\documentclass[pmlr,twocolumn,a4paper,10pt]{jmlr}
\usepackage{isipta2019}
\usepackage{mathtools}
\usepackage[american]{babel}
\usepackage{booktabs} 
\usepackage{tikz} 
\usepackage{enumitem}

\usepackage{nicefrac}
\usepackage{url}

\usepackage{stackengine}

\usepackage{graphicx}

\newcommand{\propositionref}[1]{Proposition~\ref{#1}}

\newenvironment{proofof}[1]{\par\noindent{\bfseries\upshape Proof of #1\ }}{\jmlrQED}

\newcommand{\closedaboveopenbelow}[1]{%
	\smash{\ensurestackMath{\stackengine{0.5pt}{#1}{
			\begin{tikzpicture}
				\draw circle(0.8pt);
				\draw (0.8pt, 0pt) -- (8 pt, 0pt);
			\end{tikzpicture}
		}{O}{c}{F}{F}{S}}}
	\vphantom{#1}
}

\newcommand{\norm}[1]{\left\lVert#1\right\rVert}
\newcommand{\abs}[1]{\left\lvert#1\right\rvert}

\newcommand{\nats}{\mathbb{N}}
\newcommand{\natswith}{\mathbb{N}_0}
\newcommand{\reals}{\mathbb{R}}

\newcommand{\realsext}{\overline{\reals}}
\newcommand{\realsextabove}{\closedaboveopenbelow{\reals}}

\newcommand{\states}{\mathcal{X}}
\newcommand{\gambles}{\mathcal{L}(\mathcal{X})}

\newcommand{\gamblesextabove}{\closedaboveopenbelow{\mathcal{L}}(\mathcal{X})}

\newcommand{\hprod}{\cdot}  

\newcommand{\lexpvovk}{\smash{\underline{\mathbb{E}}_{\mathcal{T}}^\mathrm{V}}}
\newcommand{\uexpvovk}{\smash{\overline{\mathbb{E}}_{\mathcal{T}}^\mathrm{V}}}

\newcommand{\lexpirr}{\smash{\underline{\mathbb{E}}_{\mathcal{T}}^\mathrm{I}}}
\newcommand{\uexpirr}{\smash{\overline{\mathbb{E}}_{\mathcal{T}}^\mathrm{I}}}

\newcommand{\lexphom}{\smash{\underline{\mathbb{E}}_{\mathcal{T}}^\mathrm{H}}}
\newcommand{\uexphom}{\smash{\overline{\mathbb{E}}_{\mathcal{T}}^\mathrm{H}}}

\newcommand{\expprec}{\smash{\mathbb{E}_{P}}}

\newcommand{\transmatset}{\mathcal{T}}

\newcommand{\imchom}{\smash{\mathcal{P}_\transmatset^\mathrm{H}}}
\newcommand{\imcirr}{\smash{\mathcal{P}_\transmatset^\mathrm{I}}}

\newcommand{\hittimeprec}{\smash{h_A^P}}
\newcommand{\hitprobprec}{\smash{p_A^P}}

\jmlrworkshop{ISIPTA 2019}

\title[Hitting Times and Probabilities for Imprecise Markov Chains]{
  Hitting Times and Probabilities for Imprecise Markov Chains
}

\author{
  \Name{Thomas {Krak}}\Email{Thomas.Krak@UGent.be}\\
  \Name{Natan {T'Joens}}\Email{Natan.TJoens@UGent.be}\\
  \Name{Jasper {De Bock}}\Email{Jasper.DeBock@UGent.be}\\
  \addr ELIS -- FLip, Ghent University, Belgium
}

\begin{document}
\maketitle

\begin{abstract}
We consider the problem of characterising 
expected hitting times and hitting probabilities for imprecise Markov chains. To this end, we consider three distinct ways in which imprecise Markov chains have been defined in the literature: as sets of homogeneous Markov chains, as sets of more general stochastic processes, and as game-theoretic probability models. Our first contribution is that all these different types of imprecise Markov chains have the same lower and upper expected hitting times, and similarly the hitting probabilities are the same for these three types. Moreover, we provide a characterisation of these quantities that directly generalises a similar characterisation for precise, homogeneous Markov chains.
\end{abstract}

\section{Introduction}\label{sec:intro}

Markov chains are mathematical models that probabilistically describe the uncertain behaviour of a dynamical system~\cite{norris:markovchains}. We here consider Markov chains that can only be in a finite number of states, and that can only change state at discrete steps in time.
An important class of inferences for Markov chains are the so called \emph{expected hitting times} and \emph{hitting probabilities} for some subset $A$ of the set of all states $\states$ that the system can be in. Informally, their aim is to answer the questions ``How long will it take until the system enters a state in $A$?'' and ``What is the probability of ever visiting a state in $A$?'', respectively. Under some regularity conditions, closed-form solutions to these questions are available in the literature~\cite{norris:markovchains, grinstead2012introduction}.

A generalisation of Markov chains that also incorporates (higher order) uncertainty about one's knowledge of the model description itself are \emph{imprecise Markov chains}~\cite{kozine2002interval, campos2003computing, hartfiel2006markov, skulj2006finite, decooman:2008:trees, decooman:2009:markovlimit, decooman:2010:markovepistemic, itip:stochasticprocesses, lopatatzidis2016robust}. Their theoretical foundations are based on the theory of imprecise probabilities~\cite{walley1991statistical, augustin:2014}, and they allow one to incorporate uncertainties about the numerical model parameters as well as about structural assumptions, like history independence---the canonical Markov property---and time homogeneity. 

However, the generalisation of Markov chains to their imprecise counterpart is not unambiguous~\cite{itip:stochasticprocesses}. There are various ways in which this might be done, and they can lead to different conclusions for particular inferences of interest. 

On the one hand we have what might be called the ``sensitivity analysis'' interpretation of an imprecise Markov chain. Here, one's model essentially constitutes an entire \emph{set} of stochastic processes that are all compatible with one's assessments about the system's uncertain behaviour. But there are multiple versions of this interpretation, depending on which models one chooses to include in this set; for instance, do we only include all (time-homogeneous) Markov chains that are compatible with our assessments~\cite{kozine2002interval, campos2003computing}, or do we also include more general stochastic processes~\cite{decooman:2010:markovepistemic,itip:stochasticprocesses}? Each choice has its own merits, depending on the particular situation. Regardless of the choice that one makes here, inferences for this ``sensitivity analysis'' interpretation always consist in computing tight lower and upper bounds on inferences for all the models that are included in the chosen set~\cite{itip:stochasticprocesses}.

An entirely different formalisation of imprecise Markov chains is based on the game-theoretic probability framework that was popularised by Shafer and Vovk~\cite{shafer2001probability}. These models are not necessarily given an interpretation in terms of compatible ``precise'' models; rather, this theory of stochastic processes is based on rational betting behaviour in repeated games with uncertain outcomes, and naturally leads to imprecise probabilistic models~\cite{decooman:2008:trees, lopatatzidis2016robust, decooman2016impreciseprocesses}. The correspondence between this framework and the ``sensitivity analysis'' interpretation of imprecise Markov chains was first explored in~\cite{decooman:2008:trees, destercke2008relating}.

In this present work, we consider the inference problems of computing lower and upper expected hitting times and hitting probabilities for an imprecise Markov chain---regardless of the specific interpretation that one chooses for these models. In fact, the first of our main results is that these inferences are the same for all of the different types of imprecise Markov chains discussed above. Our second main result is an exact generalisation to the imprecise setting, of a well-known characterisation of these inferences for precise, time-homogeneous Markov chains. 

To the best of our knowledge, this problem has never been considered in the literature at this level of generality. The most closely related work that we are aware of is that of Lopatatzidis \emph{et al.}~\cite{lopatatzidis2015calculating, lopatatzidis2017computing_birthdeath}, who prove similar properties for imprecise Markov chains that have the structure of birth-death chains. Moreover, De Cooman \emph{et al.}~\cite{decooman2016impreciseprocesses} previously derived a non-linear system describing expected hitting times, that is similar to our characterisation stated in~\corollaryref{cor:imprecise_hitting_time_is_minimal_system_solution}.

Throughout, some of the lengthier proofs as well as proofs of technical lemmas have been deferred to the appendix.

\section{Preliminaries}\label{sec:prelim}

Throughout, $\nats$ denotes the natural numbers, and we let $\natswith\coloneqq \nats\cup\{0\}$. $\reals$ denotes the real numbers, we define $\realsextabove\coloneqq \reals\cup\{+\infty\}$, that is, the reals that are closed above, and we let $\overline{\reals}\coloneqq \reals\cup\{-\infty, +\infty\}$. The sets $\realsextabove$ and $\overline{\reals}$ are endowed with the usual (order) topology, and we adopt the convention that $0\cdot+\infty = 0 = 0\cdot -\infty$.

We use $\states$ to denote the finite non-empty set of states that the Markov chain can be in. Without loss of generality, it can be identified with the set $\states=\{1,\ldots,k\}$ for some $k\in\nats$. We use $\gambles$ to denote the set of real-valued functions on $\states$. The set $\gamblesextabove$ contains all functions on $\states$ that take values in $\realsextabove$. Since $\states$ is finite, any $f$ in $\gambles$ or $\gamblesextabove$ can be identified with a vector in $\reals^k$ or ${(\realsextabove)}^k$, respectively. The set $\gambles$ is endowed with the supremum norm, i.e. $\norm{f}\coloneqq \sup_{x\in\states}\abs{f(x)}$, and the corresponding norm topology. $\gamblesextabove$ receives the topology of pointwise convergence. 

For any $A\subset\states$, we consider the indicator $\mathbb{I}_A$ of $A$, defined as $\mathbb{I}_A(x)\coloneqq 1$ if $x\in A$ and $\mathbb{I}_A(x)\coloneqq 0$ otherwise. Constant functions on $\states$ are simply denoted by their constant values.
Finally, point-wise multiplication of two functions (i.e. vectors) is denoted by $f\hprod g$; for example, the term $\mathbb{I}_{A^c}\hprod Th_A^P$ in Equation~\eqref{lemma:eq:homogeneous_hitting_system} further on denotes the pointwise multiplication of the functions $\mathbb{I}_{A^c}$ and $Th_A^P$.

\subsection{``Measure-Theoretic'' Imprecise Markov Chains}\label{subsec:measure_impr_markov}

In order to discuss the various types of imprecise Markov chains that arise from the ``sensitivity analysis'' interpretation of imprecise probabilities, we need a formalisation of general (non-Markovian) stochastic processes. We briefly give the measure theoretic account of this formalisation.

In this framework, the unknown---that is, uncertain---realisation of the stochastic process is a \emph{path}, which is a function $\omega:\mathbb{N}_0\to\mathcal{X}$. We collect all paths in the set $\Omega$. This set $\Omega$ is endowed with a $\sigma$-algebra $\mathcal{F}$\footnote{Specifically, we assume that $\mathcal{F}$ is the $\sigma$-algebra generated by the cylinder sets; this guarantees that all functions that we consider are measurable.} and augmented to a probability space $(\Omega, \mathcal{F}, P)$ with a probability measure $P$. A stochastic process is then a family $\{X_n\}_{n\in\natswith}$ of random variables on this probability space, such that $X_n: \omega\mapsto \omega(n)$ for all $n\in\natswith$. However, for ease of notation and terminology, we will often refer to the measure $P$ as the stochastic process; different processes then correspond to different measures on the same measurable space $(\Omega,\mathcal{F})$.

A \emph{Markov chain} is a stochastic process that satisfies the \emph{Markov condition}, which is a conditional independence relation between the random states that the process obtains. In particular, a process $P$ is said to be a Markov chain if
\begin{align*}
 P(X_{n+1}=x_{n+1}\,\vert\,X_{0:n}=x_{0:n}) = P(X_{n+1}=x_{n+1}\,\vert\,X_{n}=x_{n})\,,
\end{align*}
for all $x_0,\ldots,x_n,x_{n+1}\in\states$ and all $n\in\natswith$, where we let $X_{0:n}\coloneqq(X_0,\dots,X_n)$ and similarly for $x_{0:n}$. A Markov chain is called \emph{homogeneous} if, for all $x,y\in\mathcal{X}$ and all $n\in\mathbb{N}_0$,
\begin{equation}
P(X_{n+1}=y\,\vert\,X_n=x) = P(X_{1}=y\,\vert\,X_0=x)\,.
\end{equation}
Any homogeneous Markov chain $P$ is uniquely characterised---up to its initial distribution $P(X_0)$---by a \emph{transition matrix}. A transition matrix $T$ is simply an $\lvert\mathcal{X}\rvert\times \lvert\mathcal{X}\rvert$ matrix that is row-stochastic, meaning that for all $x\in\states$, $\sum_{y\in\states}T(x,y)=1$ and $T(x,y)\geq 0$ for all $y\in\states$. Such a transition matrix identifies a homogeneous Markov chain $P$ (up to its initial distribution) that satisfies
\begin{equation}\label{eq:def_transmat}
P(X_{n+1}=y\,\vert\,X_n=x) = T(x,y) 
\text{ for all $x,y\in\mathcal{X}$, $n\in\natswith$.}
\end{equation}
Moreover, a transition matrix $T$ can also be interpreted as a linear operator that maps $\gambles$ into $\gambles$, because we have identified $\gambles$ with $\reals^{\lvert\states\rvert}$. For any $f\in\gambles$ and $x\in\states$ it then holds that
\begin{align*}
\expprec\bigl[f(X_{n+1})\,\vert\,X_n=x\bigr]
&=\sum_{y\in\states}P(X_{n+1}=y\,\vert\, X_n=x)f(y)\\
&=\sum_{y\in\states}T(x,y)f(y)=\left[Tf\right](x),
\end{align*}
so we see that $T$ encodes the conditional expectation operator for 1 time step corresponding to a process $P$ that satisfies~\eqref{eq:def_transmat}. Moreover, $T$ can be uniquely extended to an operator on $\gamblesextabove$, due to the convention that $0\cdot+\infty = 0$.


We now move on to the characterisation of \emph{imprecise} Markov chains. In all cases that we consider here, these are characterised by a \emph{set} $\transmatset$ of transition matrices. In the remainder, we will assume that $\transmatset$ is non-empty, closed
, convex, and that it has separately specified rows. 
This last condition means that if, for all $x\in\states$, we select any element $T_x\in\mathcal{T}$, there must be some $T\in\mathcal{T}$ such that $T(x,\cdot)=T_x(x,\cdot)$ for all $x\in\states$; 
see e.g.~\cite[Definition 11.6]{itip:stochasticprocesses} for further discussion.

An \emph{imprecise} Markov chain is now a set of stochastic processes that are in a specific sense ``compatible'' with the transition matrices in $\mathcal{T}$. However, there are various ways how we might construct such a set, which all lead to different types of imprecise Markov chains. 

Arguably the simplest imprecise Markov chain is the set $\imchom$, which is the set of all homogeneous Markov chains whose characterising transition matrix $T$ is included in $\transmatset$~\cite{kozine2002interval, campos2003computing}. Its corresponding lower and upper expectation operators are defined respectively as
\begin{equation*}
\lexphom\bigl[\cdot\,\vert\,\cdot\bigr] \coloneqq \inf_{P\in\imchom}\expprec\bigl[\cdot\,\vert\,\cdot\bigr]
\,\,\,\text{and}\,\,\,
\uexphom\bigl[\cdot\,\vert\,\cdot\bigr] \coloneqq \sup_{P\in\imchom}\expprec\bigl[\cdot\,\vert\,\cdot\bigr]\,,
\end{equation*}
where, for both $\smash{\lexphom[\cdot\,\vert\,\cdot]}$ and $\smash{\uexphom[\cdot\,\vert\,\cdot]}$, the first argument takes functions of the form $f:\Omega\to\realsextabove$, and the second is a conditioning event $X_{0:n}=x_{0:n}$ with $n\in\natswith$.\footnote{We omit a technical discussion about the required measurability and integrability properties of such $f$, and use this definition provided that $\expprec[f\,\vert\,X_{0:n}]\coloneqq\int_\Omega f(\omega)\,\mathrm{d}P(\omega\,\vert\,X_{0:n})$ is well-defined; see e.g.~\cite{tao2011introduction} for when this is the case. For our present purposes, it suffices to know that the functions that will be of interest in this work are all non-negative and measurable, making their expectation well-defined.}

Perhaps a less obvious choice is the imprecise Markov chain $\imcirr$, which is the set of \emph{all} (potentially non-Markov, non-homogeneous) stochastic processes for which for all $n\in\natswith$ and all $x_{0},\ldots,x_{n}\in\states$:
\begin{equation}\label{eq:compatible_process}
\bigl(\exists T\in\mathcal{T}\bigr)\bigl(\forall y\in\states\bigr) : P(X_{n+1}=y\,\vert\,X_{0:n}=x_{0:n}) = T(x_n,y).
\end{equation}
The associated lower and upper expectation operators are
\begin{equation*}
\lexpirr\bigl[\cdot\,\vert\,\cdot\bigr] \coloneqq \inf_{P\in\imcirr}\expprec\bigl[\cdot\,\vert\,\cdot\bigr]
\,\,\,\text{and}\,\,\,
\uexpirr\bigl[\cdot\,\vert\,\cdot\bigr] \coloneqq \sup_{P\in\imcirr}\expprec\bigl[\cdot\,\vert\,\cdot\bigr]\,,\vspace{-5pt}
\end{equation*}
whose domain we take to be the same as that of $\lexphom[\cdot\,\vert\,\cdot]$ and $\uexphom[\cdot\,\vert\,\cdot]$.
 This type of Markov chain is often considered in the literature~\cite{decooman:2010:markovepistemic, itip:stochasticprocesses}, and is called an imprecise Markov chain under \emph{epistemic irrelevance}. 
 


Next, it will be useful to consider the dual representation(s) of the set $\transmatset$, given by the \emph{lower} (resp. \emph{upper}) \emph{transition operator} $\underline{T}$ (resp. $\overline{T}$). 
For either domain $\gambles$ or $\gamblesextabove$, these are (non-linear) operators that map these function spaces into themselves; they are respectively defined for any element $f\in\gamblesextabove$ and any $x\in\states$ as
\begin{equation}\label{eq:def:lower_trans}
\bigl[\underline{T}\,f\bigr](x) \coloneqq \smash{\inf_{T\in\mathcal{T}}} \bigl[Tf\bigr](x)~\text{and}~\bigl[\overline{T}\,f\bigr](x) \coloneqq \smash{\sup_{T\in\mathcal{T}}} \bigl[Tf\bigr](x)\,.
\end{equation}
Under the stated conditions on $\mathcal{T}$, these operators satisfy the following useful properties:
\begin{lemma}\label{lemma:trans_op_continuous}
For all $f\in\gamblesextabove$, there exist $T,S\in\mathcal{T}$ such that\vspace{-5pt}
\begin{equation}\label{eq:lower_trans_reached}
Tf=\underline{T}\,f\quad\text{and}\quad Sf=\overline{T}f.
\end{equation}
Moreover, $\underline{T}$ and $\overline{T}$ are continuous operators on $\gambles$, and are continuous on $\gamblesextabove$ with respect to non-decreasing sequences.
\end{lemma}
The usefulness of these operators stems from the fact that---similar to transition matrices for homogeneous Markov chains---they encode the (1 time step) lower and upper expectation operators for $\imchom$ and $\imcirr$. That is,
\begin{align}
\left[\overline{T}\,f\right](x_n) &= \uexphom\bigl[f(X_{n+1})\,\vert\,X_{0:n}=x_{0:n}\bigr]\notag \\
&= \uexpirr\bigl[f(X_{n+1})\,\vert\,X_{0:n}=x_{0:n}\bigr],\label{eq:MarkovJasper}
\end{align}
for all $f\in\gambles$, all $n\in\natswith$ and all $x_0,\ldots,x_n\in\states$; and similarly for the lower expectations and $\underline{T}$. Observe that the left hand side in this expression does not depend on the states $x_{0:n-1}$, which can be interpreted as saying that the (lower and upper) expectations of $\imchom$ and $\imcirr$ satisfy an \emph{imprecise Markov property}. This explains in particular why we call $\imcirr$ an ``imprecise Markov chain'', while it consists of processes which in general do not themselves satisfy the Markov property. Moreover, despite the above, it is worth noting that equality of $\lexphom[f(X_{n+m})\,\vert\,X_{0:n}]$ and $\lexpirr[f(X_{n+m})\,\vert\,X_{0:n}]$ does not in general hold when $m>1$; see e.g.~\cite[Example 10]{lopatatzidis2016robust}.

Finally, we note that in our definition and notation of the imprecise Markov chains above, we paid no further attention to the initial models $P(X_0)$ of their elements $P$. If this were to be of interest, we could specify an imprecise initial model. That is, we could consider a non-empty set $\mathcal{M}$ of probability mass functions on $\states$, and then include $P$ in $\imchom$ or $\imcirr$ if and only if, in addition to its compatibility with $\transmatset$ as discussed above, it holds that $P(X_0)\in\mathcal{M}$. 

However, we purposely restricted the domains of the corresponding lower and upper expectation operators to conditioning events of the form $X_{0:n}=x_{0:n}$, as this will suffice for all our results. Therefore, as one can easily see, these lower and upper expectations are invariant under any particular choice of such an initial model $\mathcal{M}$, which is why we have omitted any further reference to it for ease of notation and clarity of exposition.

\subsection{Hitting Times and Probabilities}

We next introduce the two inferences that are of interest in this work. The first of these is the expected \emph{hitting time} of a set of states $A\subset\states$. The hitting time $H_A:\Omega\to\mathbb{N}_{0}\cup\{+\infty\}$ for this set $A$ is a function 
defined for all $\omega\in\Omega$ as
$H_A(\omega) := \inf\left\{ t\in\mathbb{N}_{0}\,:\, \omega(t)\in A\right\}$. 
The vector of expected hitting times $\hittimeprec\in\gamblesextabove$ for a given stochastic process $P$, conditional on the starting state $X_0$, is defined for all $x\in\mathcal{X}$ as
\begin{equation*}
\hittimeprec(x) := \expprec\bigl[H_A\,\vert\,X_0=x\bigr] := \int_{\Omega}H_A(\omega)\,\mathrm{d}P(\omega\,\vert\,\omega(0)=x)\,.
\end{equation*}
Thus, $\hittimeprec(x)$ is the expected number of steps before the process $P$ reaches any element of $A$, starting from $x$.

For the imprecise Markov chains $\imchom$ and $\imcirr$, the \emph{lower} expected hitting times are defined respectively as
\begin{equation}\label{eq:lower_hitting_measure_hom}
\lexphom\bigl[H_A\,\vert\,X_0=x\bigr] := \inf_{P\in\imchom}\hittimeprec(x)\,,
\vspace{-9pt}
\end{equation}
and
\begin{equation}\label{eq:lower_hitting_measure_irr}
\lexpirr\bigl[H_A\,\vert\,X_0=x\bigr] := \inf_{P\in\imcirr}\hittimeprec(x)\,,
\end{equation}
with the corresponding \emph{upper} expected hitting times defined analogously with suprema.

The second inference that we are after is the vector of conditional \emph{hitting probabilities} $p_A^P\in\gambles$: the probabilities that the process will eventually visit an element of $A$. An explicit way of encoding this inference problem uses the function $G_A:\Omega\to \{0,1\}$, defined for all $\omega\in\Omega$ as
\begin{equation}
G_A(\omega) := \sup\{ \mathbb{I}_A\bigl(\omega(t)\bigr)\,:\,t\in\natswith\}\,.
\end{equation}
Thus, $G_A$ takes the value one on a path $\omega$ if this path at some point in time passes through any of the states in $A$; otherwise it takes the value zero. Therefore, clearly, for any stochastic process $P$, the hitting probability is given by
\begin{equation*}
p_A^P(x) \coloneqq \expprec\bigl[G_A\,\vert\,X_0=x\bigr] := \int_{\Omega}G_A(\omega)\,\mathrm{d}P(\omega\,\vert\,\omega(0)=x)\,.
\end{equation*}
Correspondingly, the lower hitting probability for the imprecise Markov chain $\imchom$ is given by
\begin{equation}
\lexphom\bigl[G_A\,\vert\,X_0=x\bigr] := \inf_{P\in\imchom} p_A^P(x)\,,
\end{equation} 
and similarly for the upper probability and for $\imcirr$.

\subsection{``Game-Theoretic'' Imprecise Markov Chains}\label{sec:Game}

In Section~\ref{subsec:measure_impr_markov} we introduced (imprecise) Markov chains using their ``measure-theoretic'' formalisation. An entirely different mathematical framework for describing stochastic processes---and imprecise Markov chains in particular---is the ``game-theoretic'' framework popularised by Shafer and Vovk~\cite{shafer2001probability}. For an in-depth treatise on this formalism, we refer the interested reader to \cite{shafer2001probability, itip:gametheory, natan:game_theory}. Explicit discussions about the connection to the measure-theoretic framework can be found in references~\cite{decooman:2008:trees,destercke2008relating, lopatatzidis2016robust}. 

For our present purposes, we restrict attention to a discussion of some essential properties of the corresponding (lower or upper) expectation operators. To this end, it suffices to think of such a game-theoretic model as simply a different characterisation of the uncertain behaviour of the dynamical system of interest. And, although this characterisation is different from the measure-theoretic one, it still leads to the same inferences in a large number of cases; in fact, it is one of the aims of this present paper to show that the expected hitting times and hitting probabilities are the same for these two different characterisations. 

The operators that we will consider in this section are functionals on functions on paths $\omega\in\Omega$. We will need a slight notational digression to introduce these domains. We let $\mathcal{L}(\Omega)$ be the set of all functions on $\Omega$ that take values in $\reals$. The domains $\closedaboveopenbelow{\mathcal{L}}(\Omega)$ and $\overline{\mathcal{L}}(\Omega)$ contain the functions taking values in $\realsextabove$ and $\overline{\reals}$, respectively. We also need the concept of an $n$-measurable function: this is a function on $\Omega$ whose value $f(\omega)$ only depends on the states $X_0$ to $X_n$. For any $n\in\nats$, we let $\mathcal{L}_n(\Omega)$ denote the set of all $n$-measurable functions taking values in $\reals$. The sets $\closedaboveopenbelow{\mathcal{L}}_n(\Omega)$ and $\overline{\mathcal{L}}_n(\Omega)$ contain the $n$-measurable functions taking values in $\realsextabove$ and $\overline{\reals}$, respectively.

A \emph{game-theoretic upper expectation operator} is now a specific $\overline{\reals}$-valued functional  $\smash{\overline{\mathbb{E}}^\mathrm{V}}[\cdot\,\vert\,\cdot]$~\cite[Definition 2]{natan:game_theory}, where the first argument takes values in $\overline{\mathcal{L}}(\Omega)$ and the second is an event of the form $X_{0:n}=x_{0:n}$. 

To specify such a game-theoretic upper expectation operator, one needs to
provide a family of operators $\overline{Q}_s\colon \gamblesextabove \to \realsextabove{}$ indexed by \emph{situations} $s \in\mathbb{S}\cup\{\Box\}$, with $\mathbb{S}\coloneqq\{x_{0:n} \in \states{}^{n+1} \colon n \in \natswith{}\}$. Furthermore, for every situation $s\in \mathbb{S}\cup\{\Box\}$, $\overline{Q}_s$ should satisfy the following axioms:
\begin{enumerate}[ref={E\arabic*},label={ E\arabic*}.,series=sepcoherence]
\item \label{coherence: const is const} 
$\overline{Q}_s(c)=c$ for all $c\in\reals{}$; 
\item \label{coherence: sublinearity} 
$\overline{Q}_s(f+g)\leq\overline{Q}_s(f) + \overline{Q}_s(g)$ for all $f,g\in\gamblesextabove$;
\item \label{coherence: homog for ext lambda}
$\overline{Q}_s(\lambda f)=\lambda \overline{Q}_s(f)$ for all positive $\lambda\in\realsextabove{}$ and all non-negative $f\in\gamblesextabove$;
\item \label{coherence: monotonicity}
if\/ $f \leq g$ then $\overline{Q}_s(f)\leq\overline{Q}_s(g)$ for all $f, g\in\gamblesextabove$.
\end{enumerate}
Crucially, every such family leads to a unique corresponding game-theoretic upper expectation operator $\smash{\overline{\mathbb{E}}^\mathrm{V}}[\cdot\,\vert\,\cdot]$. Unfortunately, however, explaining how this works requires quite some technical machinary, including the notion of a supermartingale. Since we believe this would be too much of a digression, we prefer to refer the interested reader to appendix~\ref{appendix: game-theoretic}, and
 here content ourselves with describing some of its properties.


A first important property is that every $\overline{Q}_s$ can be interpreted as a \emph{local uncertainty model} associated with the situation $s$. In particular, for $s=x_{0:n}$, we have that
\begin{equation}\label{eq:coincideswithlocalJasper}
\smash{\overline{\mathbb{E}}^\mathrm{V}}[f(X_{n+1})\,\vert\,X_{0:n}=x_{0:n}] = \overline{Q}_{x_{0:n}}(f)\,,
\end{equation}
for any $f\in\gamblesextabove$.
Similarly, the operator $\overline{Q}_{\Box}$ describes the uncertainty about the initial state. Note however that, analogous to our discussion for measure-theoretic imprecise Markov chains, we have restricted the domain of $\smash{\overline{\mathbb{E}}^\mathrm{V}}[\cdot\,\vert\,\cdot]$ to be conditional on $X_{0:n}$. Here too, as we explain in Appendix~\ref{appendix: game-theoretic}, this implies that the initial model---$\overline{Q}_{\Box}$, in this case---has no effect on our operator. For ease of notation, we will therefore make no further reference to it, and will henceforth specify $\smash{\overline{\mathbb{E}}^\mathrm{V}}[\cdot\,\vert\,\cdot]$ by providing a family of operators $\bigl\{\overline{Q}_s\bigr\}_{s\in\mathbb{S}}$, without $\overline{Q}_{\Box}$.


 We next remark that the axioms \ref{coherence: const is const}--\ref{coherence: monotonicity} can be recognised as being analogous to familiar properties of coherent lower previsions~\cite{walley1991statistical, itip:lowerprevisions}. The following result essentially shows that the upper expectation operator $\smash{\overline{\mathbb{E}}^\mathrm{V}}[\cdot\,\vert\,\cdot]$ induced by \emph{local} models that satisfy these properties, inherits these properties on the \emph{global} domain $\overline{\mathcal{L}}(\Omega)$. It also provides some properties for the conjugate game-theoretic lower expectation operator, defined as $\smash{\underline{\mathbb{E}}^\mathrm{V}}[\cdot\,\vert\,\cdot]\coloneqq -\smash{\overline{\mathbb{E}}^\mathrm{V}}[-\cdot\,\vert\,\cdot]$.

\begin{proposition}{\cite[Proposition 13]{natan:game_theory}}\label{prop:vovk_satisfies_coherence}
	Let $\overline{\mathbb{E}}^\mathrm{V}[\cdot\,\vert\,\cdot]$ be a game-theoretic upper expectation operator. Then for all $f,g\in\overline{\mathcal{L}}(\Omega)$, all $\lambda\in\reals$ with $\lambda\geq 0$, and all $n\in\natswith$:
	\begin{enumerate}
		\item $\overline{\mathbb{E}}^\mathrm{V}[f+g\,\vert\, X_{0:n}] \leq \overline{\mathbb{E}}^\mathrm{V}[f\,\vert\, X_{0:n}] + \overline{\mathbb{E}}^\mathrm{V}[g\,\vert\, X_{0:n}]$
		\footnote{
			If $f+g$ and $\overline{\mathbb{E}}^\mathrm{V}[f\,\vert\, X_{0:n}] + \overline{\mathbb{E}}^\mathrm{V}[g\,\vert\, X_{0:n}]$ 
			are well-defined; the ambiguity of $\infty + -\infty$ makes formalising this property a bit cumbersome.}
		\item $\overline{\mathbb{E}}^\mathrm{V}[\lambda f\,\vert\,X_{0:n}] = \lambda\overline{\mathbb{E}}^\mathrm{V}[f\,\vert\,X_{0:n}]$
		\item $f\leq g\,\Rightarrow \, \overline{\mathbb{E}}^\mathrm{V}[f\,\vert\,X_{0:n}] \leq \overline{\mathbb{E}}^\mathrm{V}[g\,\vert\,X_{0:n}]$
	\end{enumerate}
and, moreover, 
\begin{enumerate}
	\setcounter{enumi}{3}
	\item for all $x_0,\ldots,x_n\in\states$,
	\begin{align*}
	\inf_{\omega\in\Gamma(x_{0:n})}f(\omega) &\leq \underline{\mathbb{E}}^\mathrm{V}[f\,\vert\,X_{0:n}=x_{0:n}] \\
	&\leq
	 \overline{\mathbb{E}}^\mathrm{V}[f\,\vert\,X_{0:n}=x_{0:n}] \leq \sup_{\omega\in\Gamma(x_{0:n})}f(\omega)
	\end{align*}  with $\Gamma(x_{0:n}) \coloneqq \bigl\{ \omega\in\Omega\,\big\vert\, \forall t\in \{0,\ldots,n\}\,:\, \omega(t)=x_t \bigr\}$
	\item $\overline{\mathbb{E}}^\mathrm{V}[f + \mu\,\vert\,X_{0:n}] = \overline{\mathbb{E}}^\mathrm{V}[f\,\vert\,X_{0:n}]+\mu$ for all $\mu\in\reals$.
\end{enumerate}
\end{proposition}

With the general framework of game-theoretic upper expectation operators in place, we now move on to discussing two specific kinds of such operators, that will be particularly important in the remainder of this work. The first are those that correspond to a precise stochastic process $P$ in the measure-theoretic sense.
\begin{proposition}\label{prop: vovk local precise}
Let $P$ be a stochastic process as in Section~\ref{subsec:measure_impr_markov}, and consider the family $\bigl\{\overline{Q}_s\bigr\}_{s\in\mathbb{S}}$ defined for all $f\in\gamblesextabove$, all $n\in\natswith$, and all $x_0,\ldots,x_n\in\states$ as
\begin{equation}\label{Equation: precise local models }
\overline{Q}_{x_{0:n}}(f) \coloneqq \expprec\bigl[f(X_{n+1})\,\vert\,X_{0:n}=x_{0:n}\bigr].
\end{equation}
Then the operators in $\bigl\{\overline{Q}_s\bigr\}_{s\in\mathbb{S}}$ satisfy \ref{coherence: const is const}--\ref{coherence: monotonicity}, and therefore determine a unique corresponding game-theoretic upper expectation operator.
\end{proposition}
We will denote this game-theoretic upper expectation operator as $\smash{\overline{\mathbb{E}}_P^\mathrm{V}[\cdot\,\vert\,\cdot]}$, and the conjugate game-theoretic lower expectation operator as $\smash{\underline{\mathbb{E}}_P^\mathrm{V}[\cdot\,\vert\,\cdot]}$.

Our next result establishes that these game-theoretic operators agree with the measure-theoretic expectation $\expprec$ on all $n$-measurable real-valued functions.


\begin{proposition}\label{prop:vovk_precise_equal_measure_precise_on_n_measurable}
	Let $P$ be a stochastic process as in Section~\ref{subsec:measure_impr_markov} and let $\smash{\underline{\mathbb{E}}_P^\mathrm{V}}[\cdot\,\vert\,\cdot]$ and $\smash{\overline{\mathbb{E}}_P^\mathrm{V}}[\cdot\,\vert\,\cdot]$ be its game-theoretic lower and upper expectation operators. Then for all $n\in\natswith$ and all $f_m\in\mathcal{L}_m(\Omega)$ with $m\in\nats$, it holds that
	\begin{align*}
	\underline{\mathbb{E}}_P^\mathrm{V}[f_m\,\vert\,X_{0:n}] &= \overline{\mathbb{E}}_P^\mathrm{V}[f_m\,\vert\,X_{0:n}] = \expprec[f_m\,\vert\,X_{0:n}]\,.
	\end{align*}
\end{proposition}

The second type of game-theoretic expectation operator in which we are interested, is that corresponding to an imprecise Markov chain.
\begin{proposition}\label{prop: vovk local Markov}
	Let $\mathcal{T}$ be a non-empty, closed, and convex set of transition matrices that has separately specified rows, let $\overline{T}$ be the corresponding upper transition operator as in Section~\ref{subsec:measure_impr_markov}, and consider the family $\bigl\{\overline{Q}_s\bigr\}_{s\in\mathbb{S}}$  defined for all $f\in\gamblesextabove$, all $n\in\natswith$, and all $x_0,\ldots,x_n\in\states$ as
	\begin{equation}\label{eq:def_vovk_lexp}
	\overline{Q}_{x_{0:n}}(f) \coloneqq \bigl[\overline{T}f\bigr](x_n).
	\end{equation}
	Then the operators in $\bigl\{\overline{Q}_s\bigr\}_{s\in\mathbb{S}}$ satisfy \ref{coherence: const is const}--\ref{coherence: monotonicity}, and therefore determine a unique corresponding game-theoretic upper expectation operator.
\end{proposition}
We will denote this game-theoretic upper expectation operator as $\smash{\overline{\mathbb{E}}_\mathcal{T}^\mathrm{V}[\cdot\,\vert\,\cdot]}$, and the conjugate game-theoretic lower expectation operator as $\smash{\underline{\mathbb{E}}_\mathcal{T}^\mathrm{V}[\cdot\,\vert\,\cdot]}$.

Since the right-hand side of Equation~\eqref{eq:def_vovk_lexp} does not depend on $x_{0:(n-1)}$, it follows from Equation~\eqref{eq:coincideswithlocalJasper} that the induced game-theoretic upper expectation operator satisfies an imprecise Markov property that is entirely similar to that in Equation~\eqref{eq:MarkovJasper}. 
It is for that reason that we call $\smash{\overline{\mathbb{E}}_\mathcal{T}^\mathrm{V}[\cdot\,\vert\,\cdot]}$ and $\lexpvovk[\cdot\,\vert\,\cdot]$ the upper and lower expectation operator of a ``game-theoretic imprecise Markov chain''.

The following property shows that the operator $\lexpvovk[\cdot\,\vert\,\cdot]$ provides a lower bound for the operators $\underline{\mathbb{E}}_P^\mathrm{V}[\cdot\,\vert\,\cdot]$ whose characterising measure-theoretic process $P$ is compatible with $\mathcal{T}$. Similarly, $\uexpvovk[\cdot\,\vert\,\cdot]$ provides an upper bound.
\begin{proposition}\label{prop:vovk_imprecise_dominates_compatible_precise_vovk}
	For all $f\in\overline{\mathcal{L}}(\Omega)$, all $n\in\natswith$ and all $x_0,\ldots,x_n\in\states$, we have that\vspace{-5pt}
	\begin{equation*}
	\lexpvovk[f\,\vert\,X_{0:n}=x_{0:n}] \leq \inf_{P\in\imcirr} \underline{\mathbb{E}}_P^\mathrm{V}[f\,\vert\,X_{0:n}=x_{0:n}]\vspace{-5pt}
	\end{equation*}
	and\vspace{-6pt}
	\begin{equation*}
	\sup_{P\in\imcirr} \overline{\mathbb{E}}_P^\mathrm{V}[f\,\vert\,X_{0:n}=x_{0:n}] \leq \uexpvovk[f\,\vert\,X_{0:n}=x_{0:n}].
	\vspace{-2pt}
	\end{equation*}
\end{proposition}

Finally, we will need the following continuity property; here and in what follows, we consider $\closedaboveopenbelow{\mathcal{L}}(\Omega)$ to be endowed with the topology of pointwise convergence:
\begin{proposition}\label{prop:vovk_monotone_continuity}
	Consider a non-decreasing sequence $\{f_m\}_{m\in\nats}$ in $\mathcal{L}(\Omega)$ such that $f_m\in\mathcal{L}_m(\Omega)$ for all $m\in\nats$ and $\lim_{m\to+\infty} f_m = f\in\closedaboveopenbelow{\mathcal{L}}(\Omega)$. Then for all $n\in\natswith$ and all $x_0,\ldots,x_n\in\states$ it holds that\vspace{-5pt}
	\begin{equation*}
	\underline{\mathbb{E}}^\mathrm{V}[f\,\vert\,X_{0:n}=x_{0:n}] = \lim_{m\to+\infty} \underline{\mathbb{E}}^\mathrm{V}[f_m\,\vert\,X_{0:n}=x_{0:n}]\vspace{-6pt}
	\end{equation*}
	and\vspace{-6pt}
	\begin{equation*}
	\overline{\mathbb{E}}^\mathrm{V}[f\,\vert\,X_{0:n}=x_{0:n}] = \lim_{m\to+\infty} \overline{\mathbb{E}}^\mathrm{V}[f_m\,\vert\,X_{0:n}=x_{0:n}].
	\vspace{-8pt}
	\end{equation*}
\end{proposition}

\section{Characterisation and Invariance}\label{sec:invariance}

With the various definitions of imprecise Markov chains in place, we now move on to characterising their (lower and upper) expected hitting times and probabilities, and showing that these are the same for all of the different types of models that we discussed above. We start this discussion with our result for the hitting times, in Section~\ref{subsec:expected_hit_times}. The results for the hitting probabilities are largely analogous from a technical point of view, and are presented in Section~\ref{subsec:hit_probs}.

\subsection{Lower and Upper Expected Hitting Times}\label{subsec:expected_hit_times}

The starting point for our results in this section is the following well-known characterisation of the expected hitting times of a (precise) homogeneous Markov chain:
\begin{lemma}[\cite{norris:markovchains} {Theorem 1.3.5}]\label{lemma:homogen_hitting_time_is_minimal_system_solution}
Consider a homogeneous Markov chain $P$ with transition matrix $T$. Its vector of expected hitting times $\hittimeprec \in\gamblesextabove$ is the minimal non-negative solution to
\begin{equation}\label{lemma:eq:homogeneous_hitting_system}
\hittimeprec = \mathbb{I}_{A^c} + \mathbb{I}_{A^c}\hprod T\hittimeprec\,,
\end{equation}
where $A^c=\mathcal{X}\setminus A$, and minimality means that $\hittimeprec(x)\leq g(x)$ for all $x\in\mathcal{X}$, for any non-negative $g\in\gamblesextabove$ that also satisfies \eqref{lemma:eq:homogeneous_hitting_system}.
\end{lemma}

Inspired by this result, we introduce a recursive scheme that essentially iterates an imprecise version of Equation~\eqref{lemma:eq:homogeneous_hitting_system}. To this end, let $\smash{\underline{h}_A^{(0)}\coloneqq \overline{h}_A^{(0)}\coloneqq \mathbb{I}_{A^c}}$ and, for all $n\in\natswith$, define
\begin{equation}\label{thm:limit_is_lower_exp:eq:iteration_low}
\underline{h}_A^{(n+1)} \coloneqq \mathbb{I}_{A^c} + \mathbb{I}_{A^c}\hprod \underline{T}\,\underline{h}_A^{(n)}
\end{equation}
and
\begin{equation}\label{thm:limit_is_lower_exp:eq:iteration_up}
\overline{h}_A^{(n+1)} \coloneqq \mathbb{I}_{A^c} + \mathbb{I}_{A^c}\hprod\overline{T}\,\overline{h}_A^{(n)}.
\end{equation}
We will see in~\lemmaref{lemma:iterate_hit_time_is_restricted_exp} below that these functions can be given a clear interpretation. To this end, for all $n\in\natswith$, let $\smash{H_A^{(n)}}:\Omega\to\{0,\ldots,n+1\}$ be defined for all $\omega\in\Omega$ as
\begin{equation*}
H_A^{(n)}(\omega) \coloneqq \left\{\begin{array}{ll}
H_A(\omega) & \text{if $H_A(\omega)\leq n$, and} \\
n+1 & \text{otherwise.}
\end{array}\right.
\end{equation*}
Thus, $H_A^{(n)}(\omega)$ is the number of steps until $A$ was visited on the path $\omega$, provided that this happened in at most $n$ steps; otherwise its value is fixed to be $n+1$. The aforementioned interpretation now goes as follows:
\begin{lemma}\label{lemma:iterate_hit_time_is_restricted_exp}
For all $n\in\natswith$ it holds that
\begin{equation*}
\underline{h}_A^{(n)}=\lexpvovk\bigl[H_A^{(n)}\,\vert\,X_0\bigr]\quad\text{and}\quad \overline{h}_A^{(n)}=\uexpvovk\bigl[H_A^{(n)}\,\vert\,X_0\bigr]\,.
\end{equation*}
\end{lemma}
Moreover, it clearly holds that $\lim_{n\to+\infty}H_A^{(n)}=H_A$.
The next result tells us that the equalities in~\lemmaref{lemma:iterate_hit_time_is_restricted_exp} continue to hold as we pass to this limit; therefore, we can use the above recursive scheme to compute the (lower and upper) expected hitting times for a game-theoretic imprecise Markov chain:
\begin{proposition}\label{prop:limit_is_lower_exp}
$\lexpvovk\bigl[H_A\,\vert\,X_0\bigr] = \underline{h}_A^* \coloneqq \lim_{n\to+\infty} \smash{\underline{h}_A^{(n)}}$ 
and
$\uexpvovk\bigl[H_A\,\vert\,X_0\bigr] = \overline{h}_A^* \coloneqq \lim_{n\to+\infty} \overline{h}_A^{(n)}$. 
\end{proposition}

\begin{proof}
Each $\smash{H_A^{(n)}}$ is $n$-measurable and the sequence $\smash{H_A^{(n)}}$ is non-decreasing. Therefore, using~\lemmaref{lemma:iterate_hit_time_is_restricted_exp} and \propositionref{prop:vovk_monotone_continuity}, the limit $\underline{h}_A^*$ exists and equals $\lexpvovk\bigl[H_A\,\vert\,X_0\bigr]$. 
The proof for $\smash{\overline{h}_A^*}$ is completely analogous.
\end{proof}

In a similar manner, we can use these functions $H_A^{(n)}$ to establish that the game-theoretic hitting times corresponding to a (precise) stochastic process $P$, agree with the measure-theoretic expected hitting times of this process; this property allows us to relate the quantities obtained under the two different frameworks that we are using.
\begin{lemma}\label{lemma:precise_vovk_hit_time_equal_precise_measure_hit_time}
Let $P$ be any measure-theoretic stochastic process. Then
$\underline{\mathbb{E}}_P^\mathrm{V}[H_A\,\vert\,X_0] = \overline{\mathbb{E}}_P^\mathrm{V}[H_A\,\vert\,X_0] = \expprec[H_A\,\vert\,X_0]$.
\end{lemma}
\begin{proof}
Note that each $\smash{H_A^{(n)}}$ is $n$-measurable, that the sequence $\smash{H_A^{(n)}}$ is non-decreasing and non-negative, and that $\lim_{n\to+\infty} \smash{H_A^{(n)}}=H_A$. Hence, using~\propositionref{prop:vovk_monotone_continuity}, \propositionref{prop:vovk_precise_equal_measure_precise_on_n_measurable}, and the continuity of $\mathbb{E}_P[\cdot\,\vert\,\cdot]$ with respect to pointwise converging non-decreasing non-negative sequences (Lebesgue's monotone convergence theorem), we find that
\begin{align*}
\underline{\mathbb{E}}_P^\mathrm{V}[H_A\,\vert\,X_0] &= \lim_{n\to+\infty} \underline{\mathbb{E}}_P^\mathrm{V}[H_A^{(n)}\,\vert\,X_0] \\
 &= \lim_{n\to+\infty} \mathbb{E}_P[H_A^{(n)}\,\vert\,X_0] = \expprec[H_A\,\vert\,X_0]\,.
\end{align*}
The proof for $\overline{\mathbb{E}}_P^\mathrm{V}[H_A\,\vert\,X_0]$ is completely analogous.
\end{proof}

We now need one more property before we can state our first main result. Since the sequence $\smash{H_A^{(n)}}$ is non-decreasing, it follows from~\propositionref{prop:vovk_satisfies_coherence} together with~\lemmaref{lemma:iterate_hit_time_is_restricted_exp} that the sequences $\smash{\underline{h}_A^{(n)}}$ and $\smash{\overline{h}_A^{(n)}}$ are also non-decreasing. Hence, using the continuity of $\underline{T}$ with respect to non-decreasing sequences in $\gamblesextabove$---see~\lemmaref{lemma:trans_op_continuous}---we find that
\begin{equation}\label{eq:limit_is_fixed_point}
\underline{h}_A^* 
= \lim_{n\to+\infty} \big(\mathbb{I}_{A^c}+\mathbb{I}_{A^c}\hprod\underline{T}\, \underline{h}_A^{(n)}\big) = \mathbb{I}_{A^c}+\mathbb{I}_{A^c}\hprod\underline{T}\,\underline{h}_A^*\,.
\end{equation}
So, $\underline{h}_A^*$ is a fixed-point of the iterative scheme~\eqref{thm:limit_is_lower_exp:eq:iteration_low}. Similarly, $\overline{h}_A^*$ is a fixed-point of~\eqref{thm:limit_is_lower_exp:eq:iteration_up}.

By combining Equation~\eqref{eq:limit_is_fixed_point} with the properties of game-theoretic expectation operators and the known characterisation for precise, homogeneous Markov chains in~\lemmaref{lemma:homogen_hitting_time_is_minimal_system_solution}, we can now derive the following remarkable consequence; it states that the (lower and upper) expected hitting time for any type of imprecise Markov chain is obtained by a homogeneous Markov chain that is compatible with it. Consequently, the (lower and upper) expected hitting times are the same for all types of imprecise Markov chains!

\begin{theorem}\label{thm:lower_hitting_time_reach_and_equal}
There exists a $P\in\mathcal{P}_\mathcal{T}^\mathrm{H}$ such that
\begin{equation}\label{cor:lower_hitting_time_homogen:eq:is_equal_low}
\lexpvovk\bigl[H_A\,\vert\,X_0\bigr] = \lexpirr\bigl[H_A\,\vert\,X_0\bigr] = \lexphom\bigl[H_A\,\vert\,X_0\bigr] = \expprec\bigl[H_A\,\vert\,X_0\bigr]\,.
\end{equation}
Moreover, there exists a $P\in\mathcal{P}_\mathcal{T}^\mathrm{H}$, such that
\begin{equation}\label{cor:lower_hitting_time_homogen:eq:is_equal_up}
\uexpvovk\bigl[H_A\,\vert\,X_0\bigr] = \uexpirr\bigl[H_A\,\vert\,X_0\bigr] = \uexphom\bigl[H_A\,\vert\,X_0\bigr] = \expprec\bigl[H_A\,\vert\,X_0\bigr]\,.
\end{equation}
\end{theorem}
\begin{proof}
From the fixed-point property~\eqref{eq:limit_is_fixed_point} and the reachability property~\eqref{eq:lower_trans_reached}, we find a $T\in\mathcal{T}$ such that
\begin{equation}\label{cor:lower_hitting_time_homogen:eq:precise_system}
\underline{h}_A^* = \mathbb{I}_{A^c} + \mathbb{I}_{A^c}\hprod T\,\underline{h}_A^*\,.
\end{equation}
Using~\eqref{eq:def_transmat}, we find a homogeneous Markov chain $P$ with transition matrix $T$, and clearly $P\in\imchom$. It remains to show that~\eqref{cor:lower_hitting_time_homogen:eq:is_equal_low} holds for this $P$.

Since $\underline{h}_A^*$ satisfies~\eqref{cor:lower_hitting_time_homogen:eq:precise_system}, it clearly is a solution to~\eqref{lemma:eq:homogeneous_hitting_system}. Hence, by~\lemmaref{lemma:homogen_hitting_time_is_minimal_system_solution} and~\propositionref{prop:limit_is_lower_exp}, it holds that
\begin{equation*}
\expprec\bigl[H_A\,\vert\,X_0\bigr] \leq \underline{h}_A^* = \lexpvovk\bigl[H_A\,\vert\,X_0\bigr]\,.
\end{equation*}
Conversely, we infer from~\propositionref{prop:vovk_imprecise_dominates_compatible_precise_vovk} and~\lemmaref{lemma:precise_vovk_hit_time_equal_precise_measure_hit_time} that
\begin{align*}
\lexpvovk\bigl[H_A\,\vert\,X_0\bigr] &\leq \inf_{Q\in\imcirr} \underline{\mathbb{E}}_Q^\mathrm{V}\bigl[H_A\,\vert\,X_0\bigr] \\
 &= \inf_{Q\in\imcirr} \mathbb{E}_Q\bigl[H_A\,\vert\,X_0\bigr] \\
 &= \lexpirr\bigl[H_A\,\vert\,X_0\bigr] \\
 & \leq \lexphom\bigl[H_A\,\vert\,X_0\bigr] \leq \expprec\bigl[H_A\,\vert\,X_0\bigr]\,,
\end{align*}
where the last two inequalities hold since $P\in\imchom\subseteq\imcirr$.

The proof for the upper expected hitting time is far more tedious; it can be found in the appendix.
\end{proof}


We want to stress how powerful this result is: no matter what kind of imprecise generalisation of a Markov chain one wishes to use, the corresponding expected hitting time will always be the same (provided the regularity conditions of the set $\mathcal{T}$ are satisfied). This is not only powerful from a theoretical point of view; numerically, it allows one to use algorithms for computing (lower and upper) expectations of, say, a game-theoretic model, even when the model that one is using is a set of homogeneous Markov chains.

We conclude this section with the following characterisation of the lower and upper expected hitting times of an arbitrary imprecise Markov chain; note that this is a direct generalisation of~\lemmaref{lemma:homogen_hitting_time_is_minimal_system_solution}.
\begin{corollary}\label{cor:imprecise_hitting_time_is_minimal_system_solution}
Consider an imprecise Markov chain with lower transition operator $\underline{T}$ and upper transition operator $\overline{T}$. Its vector of lower expected hitting times $\underline{h}_A\in\gamblesextabove$ is the minimal non-negative solution to
\begin{equation}\label{cor:imprecise_hitting_time_is_minimal_system_solution:eq:system}
\underline{h}_A = \mathbb{I}_{A^c} + \mathbb{I}_{A^c}\hprod\underline{T}\,\underline{h}_A\,,
\end{equation}
and its vector of upper expected hitting times $\overline{h}_A\in\gamblesextabove$ is the minimal non-negative solution to
\begin{equation}\label{cor:imprecise_hitting_time_is_minimal_system_solution:eq:system_up}
\overline{h}_A = \mathbb{I}_{A^c} + \mathbb{I}_{A^c}\hprod\overline{T}\,\overline{h}_A\,.
\end{equation}
\end{corollary}
\begin{proof}
Due to \theoremref{thm:lower_hitting_time_reach_and_equal}, the lower expected hitting time is the same for every type of imprecise Markov chain; let $\underline{h}_A:=\lexphom\bigl[H_A\,\vert\,X_0\bigr]=\lexpvovk\bigl[H_A\,\vert\,X_0\bigr]$ be this lower expected hitting time. That $\underline{h}_A$ satisfies~\eqref{cor:imprecise_hitting_time_is_minimal_system_solution:eq:system} is immediate from~\propositionref{prop:limit_is_lower_exp} and~\eqref{eq:limit_is_fixed_point}. That it is non-negative follows from the non-negativity of $H_A$.

It remains to show that it is the minimal solution. So, let $\underline{g}_A\in\gamblesextabove$ be any non-negative solution of~\eqref{cor:imprecise_hitting_time_is_minimal_system_solution:eq:system}, and suppose \emph{ex absurdo} that $\underline{g}_A(x)<\underline{h}_A(x)$ for some $x\in\mathcal{X}$.

From~\eqref{cor:imprecise_hitting_time_is_minimal_system_solution:eq:system} and~\eqref{eq:lower_trans_reached}, we find a $T\in\mathcal{T}$ such that
\begin{equation}\label{cor:lower_hitting_time_is_minimal_system:eq:precise_system}
\underline{g}_A = \mathbb{I}_{A^c} + \mathbb{I}_{A^c}\hprod T\,\underline{g}_A\,.
\end{equation}
Using~\eqref{eq:def_transmat}, we find a homogeneous Markov chain $P$ with transition matrix $T$, and clearly $P\in\imchom$. By~\lemmaref{lemma:homogen_hitting_time_is_minimal_system_solution} and~\eqref{cor:lower_hitting_time_is_minimal_system:eq:precise_system} we conclude that
\begin{equation*}
h_A^P(x) \leq \underline{g}_A(x) < \underline{h}_A(x) = \lexphom\bigl[H_A\,\vert\,X_0=x\bigr]\,,
\end{equation*}
which yields a contradiction using~\eqref{eq:lower_hitting_measure_hom}.

The proof of the corresponding statement for the upper expected hitting time is a bit more tedious and can be found in the appendix.
\end{proof}


\subsection{Lower and Upper Hitting Probabilities}\label{subsec:hit_probs}

We now move on to study the (lower and upper) hitting probabilities of an imprecise Markov chain. Both the discussion and the technical results largely mirror that of the hitting times in Section~\ref{subsec:expected_hit_times}. We again start with a well-known characterisation for (precise) homogeneous Markov chains:

\begin{lemma}[\cite{norris:markovchains} {Theorem 1.3.2}]\label{lemma:homogen_hitting_prob_is_minimal_system_solution}
Consider a homogeneous Markov chain $P$ with transition matrix $T$. Its vector of hitting probabilities $\hitprobprec \in\gambles$ is the minimal non-negative solution to
\begin{equation}\label{lemma:eq:homogeneous_hitting_prob_system}
\hitprobprec = \mathbb{I}_{A} + \mathbb{I}_{A^c}\hprod T\hitprobprec\,.
\end{equation}
\end{lemma}

Once more, we proceed by defining a recursive scheme that is inspired by this characterisation: we let $\smash{\underline{p}_A^{(0)}} \coloneqq \smash{\overline{p}_A^{(0)}}\coloneqq \mathbb{I}_A$ and, for all $n\in\natswith$, we define
\begin{equation}\label{thm:hit_prob_is_limit:eq:recursion}
\underline{p}_A^{(n+1)} := \mathbb{I}_A + \mathbb{I}_{A^c}\hprod \underline{T}\,\underline{p}_A^{(n)}\vspace{-4pt}
\end{equation}
and
\begin{equation}\label{thm:hit_prob_is_limit:eq:recursion_up}
\overline{p}_A^{(n+1)} := \mathbb{I}_A + \mathbb{I}_{A^c}\hprod \overline{T}\,\overline{p}_A^{(n)}.
\vspace{2pt}
\end{equation}
In order to give these functions a clear interpretation, we require some auxiliary functions. For all $n\in\natswith$, we let $G_A^{(n)}:\Omega\to\{0,1\}$ be defined for all $\omega\in\Omega$ as
\begin{equation*}
G_A^{(n)}(\omega) := \sup\bigl\{ \mathbb{I}_A\bigl(\omega(t)\bigr)\,:\,t\in\{0,\ldots,n\} \bigr\}\,.
\end{equation*}
Thus $G_A^{(n)}$ takes the value one on $\omega$ if $\omega$ visits $A$ in the first $n$ steps; otherwise it takes the value zero. The aforementioned interpretation now goes as follows:
\begin{lemma}\label{lemma:recursion_lower_prob_is_truncated_function}
For all $n\in\natswith$ it holds that
\begin{equation*}
\underline{p}_A^{(n)} = \lexpvovk\bigl[ G_A^{(n)}\,\vert\,X_0 \bigr]
\quad\text{and}\quad
\overline{p}_A^{(n)} = \uexpvovk\bigl[ G_A^{(n)}\,\vert\,X_0 \bigr]\,.
\end{equation*}
\end{lemma}
Moreover, we again clearly have that $\lim_{n\to+\infty} G_A^{(n)} = G_A$. As the following result tells us, the equalities in~\lemmaref{lemma:recursion_lower_prob_is_truncated_function} continue to hold as we pass to this limit; so, the above recursive scheme can be used to compute the (lower and upper) hitting probabilities for a game-theoretic imprecise Markov chain:
\begin{proposition}\label{prop:hit_prob_is_limit}
$\lexpvovk\big[G_A\,\vert\,X_0\bigr] = \underline{p}_A^* \coloneqq \lim_{n\to+\infty} \underline{p}_A^{(n)}$ 
and
$\uexpvovk\big[G_A\,\vert\,X_0\bigr] = \overline{p}_A^* \coloneqq \lim_{n\to+\infty} \overline{p}_A^{(n)}$. 
\end{proposition}
\begin{proof}
First, we remark that each $\smash{G_A^{(n)}}$ is $n$-measurable and that the sequence $\smash{G_A^{(n)}}$ is non-decreasing. Therefore, using~\lemmaref{lemma:recursion_lower_prob_is_truncated_function} and~\propositionref{prop:vovk_monotone_continuity}, the limit $\underline{p}_A^{*}$ exists and equals $\lexpvovk\bigl[G_A\,\vert\,X_0\bigr]$. 
The proof for $\overline{p}_A^*$ is completely analogous.
\end{proof}

We can also use these functions $G_A^{(n)}$ to establish a connection between the game-theoretic hitting probabilities corresponding to a (precise) stochastic process $P$, and the measure-theoretic hitting probabilities of this process:
\begin{lemma}\label{lemma:precise_vovk_hit_prob_equal_precise_measure_hit_prob}
Let $P$ be any measure-theoretic stochastic process. Then
$\underline{\mathbb{E}}_P^\mathrm{V}[G_A\,\vert\,X_0] = \overline{\mathbb{E}}_P^\mathrm{V}[G_A\,\vert\,X_0] = \expprec[G_A\,\vert\,X_0]$.
\end{lemma}
\begin{proof}
Completely analogous to the proof of~\lemmaref{lemma:precise_vovk_hit_time_equal_precise_measure_hit_time}.
\end{proof}

Since each $\underline{p}_A^{(n)},\underline{p}_A^*\in\gambles$ is bounded---this follows from the boundedness of $G_A$ and each $\smash{G_A^{(n)}}$, together with~\lemmaref{lemma:recursion_lower_prob_is_truncated_function},~\propositionref{prop:hit_prob_is_limit}, and ~\propositionref{prop:vovk_satisfies_coherence}---the continuity of $\underline{T}$ on $\gambles$ immediately yields the fixed-point property of the iterative scheme~\eqref{thm:hit_prob_is_limit:eq:recursion}:
\begin{equation}\label{eq:fixed_point_hit_prob}
\underline{p}_A^* = \lim_{n\to+\infty} \big(\mathbb{I}_A + \mathbb{I}_{A^c}\hprod \underline{T}\,\underline{p}_A^{(n)}\big) = \mathbb{I}_A + \mathbb{I}_{A^c}\hprod \underline{T}\,\underline{p}_A^{*}\,.
\end{equation}
Similarly, $\overline{p}_A^*$ is a fixed-point of the scheme~\eqref{thm:hit_prob_is_limit:eq:recursion_up}. We again conclude that the lower and upper hitting probabilities are the same for all types of imprecise Markov chains:
\begin{theorem}\label{thm:imprecise_hit_probs_equal}
There exists a $P\in\imchom$ such that
\begin{equation}
\lexpvovk\bigl[G_A\,\vert\,X_0\bigr] = \lexpirr\bigl[G_A\,\vert\,X_0\bigr] = \lexphom\bigl[G_A\,\vert\,X_0\bigr] = \expprec\bigl[G_A\,\vert\,X_0\bigr]\,. 
\end{equation}
Moreover, it holds that\footnote{But note that the upper hitting probability is not necessarily reached by any $P\in\imchom$.}
\begin{equation}
\uexpvovk\bigl[G_A\,\vert\,X_0\bigr] = \uexpirr\bigl[G_A\,\vert\,X_0\bigr] = \uexphom\bigl[G_A\,\vert\,X_0\bigr]\,. 
\end{equation}
\end{theorem}
\begin{proof}
The proof for the lower hitting probability is completely analogous to the proof for the lower expected hitting time in~\theoremref{thm:lower_hitting_time_reach_and_equal}, only relying on ~\lemmaref{lemma:homogen_hitting_prob_is_minimal_system_solution} instead of~\lemmaref{lemma:homogen_hitting_time_is_minimal_system_solution}; on the fixed-point property~\eqref{eq:fixed_point_hit_prob} instead of~\eqref{eq:limit_is_fixed_point}; on~\propositionref{prop:hit_prob_is_limit} instead of~\propositionref{prop:limit_is_lower_exp}; and on~\lemmaref{lemma:precise_vovk_hit_prob_equal_precise_measure_hit_prob} instead of~\lemmaref{lemma:precise_vovk_hit_time_equal_precise_measure_hit_time}. The proof for the upper hitting probability is again far more tedious and can be found in  the appendix.
\end{proof}

We close with a characterisation of the lower and upper hitting probabilities for an arbitrary imprecise Markov chain, that directly generalises~\lemmaref{lemma:homogen_hitting_prob_is_minimal_system_solution}.
\begin{corollary}\label{cor:imprecise_hitting_prob_is_minimal_system_solution}
Consider an imprecise Markov chain with lower transition operator $\underline{T}$ and upper transition operator $\overline{T}$. Its vector of lower hitting probabilities $\underline{p}_A\in\gambles$ is the minimal non-negative solution to
\begin{equation}\label{cor:imprecise_hitting_prob_is_minimal_system_solution:eq:system}
\underline{p}_A = \mathbb{I}_A + \mathbb{I}_{A^c}\hprod \underline{T}\,\underline{p}_A\,,
\end{equation}
and its vector of upper hitting probabilities $\overline{p}_A\in\gambles$ is the minimal non-negative solution to
\begin{equation}\label{cor:imprecise_hitting_prob_is_minimal_system_solution:eq:system_up}
\overline{p}_A = \mathbb{I}_A + \mathbb{I}_{A^c}\hprod \overline{T}\,\overline{p}_A\,.
\end{equation}
\end{corollary}
\begin{proof}
The proof for $\underline{p}_A$ is completely analogous to the proof for $\underline{h}_A$ in~\corollaryref{cor:imprecise_hitting_time_is_minimal_system_solution}, only relying on~\theoremref{thm:imprecise_hit_probs_equal} instead of~\theoremref{thm:lower_hitting_time_reach_and_equal}; on~\propositionref{prop:hit_prob_is_limit} instead of~\propositionref{prop:limit_is_lower_exp}; on the fixed-point property~\eqref{eq:fixed_point_hit_prob} instead of~\eqref{eq:limit_is_fixed_point}; and on~\lemmaref{lemma:homogen_hitting_prob_is_minimal_system_solution} instead of~\lemmaref{lemma:homogen_hitting_time_is_minimal_system_solution}.

The proof for $\overline{p}_A$ is again different, 
 and is given in the appendix.
\end{proof}

\section{Summary and Discussion}

We have studied lower and upper expected hitting times and probabilities for imprecise Markov chains. To this end, we considered three different ways in which an imprecise Markov chain might be defined: as a set of precise, homogeneous Markov chains; as a set of precise but general (\emph{non-Markovian}) stochastic processes; and as a game-theoretic model with imprecise local models. We have shown that these quantities of interest are the same for all these types of imprecise Markov chains. Moreover, we have presented characterisations of these quantities that are direct generalisations of their well-known counterparts for precise homogeneous Markov chains.

One unexplored line of research would be to investigate the connections of these results to the theory of Markov Decision Processes (MDPs)~\cite{feinberg2012handbook}. In an MDP, the aim is to choose, at each point in time $n\in\natswith$, an \emph{action} $a_n$ from an admissible action set $\mathbb{A}_n(x_n)$ that determines the transition probabilities $P(X_{n+1}=x_{n+1}\,\vert\,X_n=x_n)$. If we interpret the choice of these actions in our current context as a selection of $T\in\transmatset$, then the connection between MDPs and the theory of imprecise Markov chains becomes intuitively clear. It is worth mentioning that this connection has been known for a while---see, e.g., the introductions of~\cite{troffaes2013model, krak2017imprecise}---yet an important semantic difference has always been the goal with which actions are selected. In an imprecise Markov chain, we optimise over $\transmatset$ in order to compute bounds on inferential quantities of interest; the goal is the quantification of uncertainty. In contrast, in an MDP, the intended outputs are typically the optimal actions themselves, which are selected to optimise a given utility function that is typically interpreted as an operational reward. However, as it pertains to the results in this current work, in Corollaries~\ref{cor:imprecise_hitting_time_is_minimal_system_solution} and~\ref{cor:imprecise_hitting_prob_is_minimal_system_solution} the characterising equations that we have derived are very reminiscent of the equations of optimality that one often encounters in the theory of MDPs, and it would be very interesting to see if this connection could be formalised.



Finally, we hope in future work to derive efficient algorithms for numerically computing the inferences that we have discussed, and aim to also extend our results to the setting of imprecise continuous-time Markov chains~\cite{skulj2015efficient, krak2017imprecise}.

\appendix
\section{Introduction to Game-theoretic Upper Expectations}\label{appendix: game-theoretic}
For readers that would like to have a better understanding of game-theoretic upper expectations, and how they are derived from the local models $\overline{Q}_s$, this appendix provides a brief introduction to the game-theoretic probability framework of Shafer and Vovk \cite{shafer2001probability,Vovk2019finance}.

Rather than using transition probabilities to describe the uncertain behaviour of a process, they assume that, for every situation $s \in \mathbb{S}\cup\{\Box\}$, we are given an operator $\overline{Q}_s\colon \gamblesextabove \to \realsextabove{}$ that satisfies \ref{coherence: const is const}--\ref{coherence: monotonicity}.
As we have already mentioned in Section~\ref{sec:Game}, these operators can be thought of as the local uncertainty models of the process.
Quite conveniently, as we have for example seen in Propositions~\ref{prop: vovk local precise} and~\ref{prop: vovk local Markov}, such a local model can be a linear expectation operator corresponding to a probability mass function on $\states$, or an upper envelope of a set of such linear expectation operators (provided the corresponding set of probability mass functions is closed and convex).
In the game-theoretic framework, however, these local models will typically be interpreted as representing the bets that a subject is willing to offer to others.
To do this, a function $f(X_{n+1})$, with $f \in \gamblesextabove{}$ and $n \in \natswith{}$, is regarded as a bet that yields a (possibly negative) uncertain reward $f(x)$ if $X_{n+1} = x$.\footnote{Note that the reward associated with these bets may also be equal to $+\infty$. In that case, it is not immediately clear how we should interpret these bets. This topic is for instance discussed in \cite{Tjoens2019NaturalExtensionISIPTA}.}
The adopted interpretation for the local model $\overline{Q}_{x_{0:n}}$ is then that conditional on the fact that he observed $X_{0:n} = x_{0:n}$, the subject is willing to offer the bet $f(X_{n+1})$, for any $f \in \gamblesextabove{}$ such that $\overline{Q}_{x_{0:n}}(f) \leq 0$.
Axioms \ref{coherence: const is const}--\ref{coherence: monotonicity} can then be regarded as constraints on what it means to offer bets rationally.
In the same way, the operator $\overline{Q}_{\Box}(\cdot) \colon \gamblesextabove \to \realsextabove{}$ represents bets on the initial state $X_0$ that the subject is willing to offer.

The idea is now to combine all these local bets to obtain a global uncertainty model ${\overline{\mathbb{E}}^\mathrm{V}}[\,\cdot\,\vert\,\cdot]$ that extends the information that is gathered in the local models.
This is achieved using the concept of a \emph{supermartingale}.

Formally, we define a supermartingale $\mathcal{M}$ to be an extended real-valued map on $\mathbb{S}\cup\{\Box\}$ that is uniformly bounded below, i.e. there is a real $c$ such that $\mathcal{M}(s) \geq c$ for all $s \in \mathbb{S}\cup\{\Box\}$, and that satisfies $\overline{Q}_s(\mathcal{M}(s \, \cdot) \, ) \leq \mathcal{M}(s)$ for all $s \in \mathbb{S}\cup\{\Box\}$.
Here, we used $\mathcal{M}(s \, \cdot)$ to denote the function in $\gamblesextabove$ that takes the value $\mathcal{M}(s x)$ for each $x \in \states{}$; note that $\Box$ can be interpreted as the empty string, and thus when $s=\Box$ it holds that $s x=\Box x=x$ for all $x\in\states$.
Indeed, $\mathcal{M}$ is uniformly bounded below, so $\mathcal{M}(s \, \cdot)$ will only take values in $\realsextabove{}$.
The key property here is that a supermartingale $\mathcal{M}$ should satisfy $\overline{Q}_s(\mathcal{M}(s \, \cdot) \, ) \leq \mathcal{M}(s)$ for all $s \in \mathbb{S}\cup\{\Box\}$, which essentially states that $\mathcal{M}$ represents a possible way to take the subject up on the bets that he is offering.
Indeed, if, for the sake of simplicity, we assume that $\mathcal{M}(s)$ is real, then it follows from the constant additivity of $\overline{Q}_s$---see \cite[Proposition 1]{natan:game_theory}---that $\overline{Q}_s(\mathcal{M}(s \, \cdot) \, ) \leq \mathcal{M}(s)$ is equivalent to $\overline{Q}_s(\mathcal{M}(s \, \cdot) - \mathcal{M}(s)) \leq 0$, which implies that our subject is willing to offer the bet $\mathcal{M}(s \, \cdot) - \mathcal{M}(s)$.
In this way, it becomes clear that $\mathcal{M}$ describes a possible evolution of a person's capital when he is gambling according to the bets offered by our subject.
For a given choice of local models $\{ \overline{Q}_s \}_{s \in \mathbb{S}}$ and $\overline{Q}_\Box$, we will use $\smash{\overline{\mathbb{M}}_{b}}$ to denote the corresponding set of all such supermartingales.

Consider now any $f \in \overline{\mathcal{L}}(\Omega)$ and $s=x_{0:n}\in\mathbb{S}$. The global (game-theoretic) upper expectation of $f$ conditional on $s$ is then defined by
\begin{align*}
\smash{\overline{\mathbb{E}}^\mathrm{V}}[\,f\,\vert\, s] \coloneqq \inf\{ \mathcal{M}(s) \colon \mathcal{M} \in \overline{\mathbb{M}}_b \text{ and } \liminf \mathcal{M} \geq_s f \},
\end{align*}
where $\liminf \mathcal{M} \geq_s f$ is taken to mean that for every path $\omega=z_0z_1\ldots z_n\ldots\in\Omega$ such that $z_{0:n}=x_{0:n}$, $\liminf_{m\to+\infty}\mathcal{M}(z_{0:m})\geq f(\omega)$. 



An intuitive meaning can be given to these upper expectations if we interpret $f \in \overline{\mathcal{L}}(\Omega)$ as an uncertain reward that depends on the path $\omega \in \Omega$ that is taken by the process. In particular, the  game-theoretic upper expectation $\smash{\overline{\mathbb{E}}^\mathrm{V}}[\,f\,\vert\, s \,]$ can then be interpreted as the infimum starting capital $\mathcal{M}(s)$ that is needed in order to guarantee that, by starting in the situation $s=x_{0:n}$ and then gambling against the subject in an appropriate (and allowed) way, we can be \emph{sure} to (eventually) end up with a capital that is larger than the reward that is associated with $f$, in the sense that this will be true for every path $\omega=z_0z_1\ldots z_n\ldots\in\Omega$ such that $z_{0:n}=x_{0:n}$.
In other words, if we are in the situation $X_{0:n} = x_{0:n}$, then any capital $\alpha$ larger than $\smash{\overline{\mathbb{E}}^\mathrm{V}}[\,f\,\vert\, s \,]$ should be worth more to us than the uncertain reward $f$, because we can bet with $\alpha$ to (eventually) obtain a reward that is guaranteed to be higher than $f$. We can therefore regard $\smash{\overline{\mathbb{E}}^\mathrm{V}}[\,f\,\vert\, s \,]$ as a lower bound on these capitals $\alpha$.

As can be expected from this interpretation, the upper expectation $\smash{\overline{\mathbb{E}}^\mathrm{V}}[\,f\,\vert\, s \,]$ does not depend on the chosen initial model $\overline{Q}_\Box$, thereby justifying our claim in Section~\ref{sec:Game}. The following result formalizes this.

\begin{proposition}
For any $f \in \gamblesextabove{}$ and $s \in \mathbb{S}$, the upper expectation ${\overline{\mathbb{E}}^\mathrm{V}}[\,f\,\vert\, s \,]$ does not depend on the choice of $\overline{Q}_\Box$.
\end{proposition}

Finally, we want to add that the global game-theoretic upper expectation operator $\smash{\overline{\mathbb{E}}^\mathrm{V}}[\, \cdot \,\vert\, \cdot \,]$ can also be characterised in a completely different way, without the use of game-theoretic concepts such as supermartingales.
Indeed, in \cite{Tjoens2019NaturalExtensionISIPTA}, it is shown that this operator is the most conservative---so least informative---upper expectation that is consistent with the local uncertainty models and satisfies a number of basic rationality axioms. 

For more information on the subject of game-theoretic probabilities, we refer the interested reader to the textbooks of Shafer and Vovk \cite{shafer2001probability,Vovk2019finance}.

\acks{%
The work in this paper was partially supported by H2020-MSCA-ITN-2016 UTOPIAE, grant agreement 722734.
We are very grateful to three anonymous reviewers for their helpful suggestions and constructive comments.
}

\bibliography{isipta2019-template}


\clearpage
\section{Proofs of Statements in \sectionref{sec:prelim}}\label{appendix:prelim_proofs}

We already introduced transition matrices in~\sectionref{subsec:measure_impr_markov}. Since we will be working with such transition matrices a lot, we start by repeating the definition explicitly.
\begin{definition}
An $\lvert\states\rvert \times \lvert\states\rvert$ matrix $T$ is called a \emph{transition matrix} if, for all $x\in\states$, it satisfies
\begin{enumerate}
\item $T(x,y)\geq 0$ for all $y\in\states$, and
\item $\sum_{y\in\states} T(x,y)=1$\,.
\end{enumerate}
\end{definition}
In other words, a matrix $T$ is a transition matrix if, for all $x\in\states$, its $x$-th row $T(x,\cdot)$ is a probability mass function on $\states$. We remarked in~\sectionref{subsec:measure_impr_markov} that the convention $0\cdot+\infty=0$ uniquely determines the extension of $T$ from $\gambles$ to $\gamblesextabove$. To see this, we first recall that for any $f\in\gambles$ and any $x\in\states$, it holds that
\begin{equation*}
\bigl[Tf\bigr](x) \coloneqq \sum_{y\in\states} T(x,y)f(y)\,.
\end{equation*}
The extension to $\gamblesextabove$ is then indeed straightforward; for any $f\in\gamblesextabove$ and any $x\in\states$,
\begin{align*}
\bigl[Tf\bigr](x) \coloneqq \sum_{y\in\states} T(x,y)f(y)
 &= \sum_{y\in\states : f(y)\in\reals} T(x,y)f(y) \\
  &\quad + \sum_{y\in\states:f(y)=+\infty} T(x,y)f(y)
\end{align*}
from which it follows that
\begin{equation*}
\bigl[Tf\bigr](x) = +\infty\,,
\end{equation*}
if $T(x,y)>0$ and $f(y)=+\infty$ for some $y\in\states$, and
\begin{equation*}
\bigl[Tf\bigr](x) = \sum_{y\in\states} T(x,y)f(y) = \sum_{y\in\states : f(y)\in\reals} T(x,y)f(y)\,,
\end{equation*}
otherwise, where we used the convention $0\cdot+\infty=0$ in the last step.

Hence, for any $T$ and any $x\in\states$, $\bigl[T(\cdot)\bigr](x)$ maps $\gamblesextabove$ into $\realsextabove$. Therefore, the infimum and supremum in the definitions of $\underline{T}$ and $\overline{T}$ in~\eqref{eq:def:lower_trans} exist. Moreover, as claimed in~\sectionref{subsec:measure_impr_markov}, these operators therefore map the spaces $\gambles$ and $\gamblesextabove$ into themselves; this is a property that we will often use without mentioning it explicitly.
\begin{proposition}
For any $f\in\gambles$ it holds that $\underline{T}\,f \in \gambles$ and $\overline{T}f\in\gambles$.
Moreover, for any $f\in\gamblesextabove$ it holds that $\underline{T}\,f \in \gamblesextabove$ and $\overline{T}f\in\gamblesextabove$.
\end{proposition} 
\begin{proof}
First, fix any $f\in\gambles$. Then, because $\lvert\states\rvert$ is finite, $f$ is uniformly bounded both below and above, meaning that $\min_{x\in\states}f(x)\in\reals$ and $\max_{x\in\states}f(x)\in\reals$. 

Moreover, for any $T\in\transmatset$ and any $x\in\states$ it holds that
\begin{equation*}
\bigl[Tf\bigr](x) = \sum_{y\in\states} T(x,y)f(y) \geq \min_{y\in\states} f(y)\,,
\end{equation*}
and 
\begin{equation*}
\bigl[Tf\bigr](x) = \sum_{y\in\states} T(x,y)f(y) \leq \max_{y\in\states} f(y)\,,
\end{equation*}
because $T$ is row-stochastic. Hence, for any $x\in\states$,
\begin{equation*}
\bigl[\underline{T}\,f\bigr](x) = \inf_{T\in\transmatset} \bigl[Tf\bigr](x) \geq \min_{y\in\states}f(y) \in\reals\,,
\end{equation*}
and 
\begin{equation*}
\bigl[\underline{T}\,f\bigr](x) = \inf_{T\in\transmatset} \bigl[Tf\bigr](x) \leq \max_{y\in\states}f(y) \in\reals\,.
\end{equation*}
This means that $\bigl[\underline{T}\,f\bigr](x)\in\reals$ and, because $x\in\states$ is arbitrary, this implies that $\underline{T}\,f\in\gambles$. Repeating this argument, \emph{mutatis mutandis}, will show that also $\overline{T}\,f\in\gambles$.

We now move on to prove the second statement, so, fix any $f\in\gamblesextabove$. Then, again because $\lvert\states\rvert$ is finite, $f$ is uniformly bounded below, meaning that there is some $c\in\reals$ such that $c\leq \min_{x\in\states}f(x)$. Then again for every $T\in\transmatset$ and any $x\in\states$ it holds that
\begin{equation*}
\bigl[Tf\bigr](x) = \sum_{y\in\states} T(x,y)f(y) \geq \min_{y\in\states}f(y)\,,
\end{equation*}
because $T$ is row-stochastic. Hence it follows that
\begin{equation*}
\bigl[\underline{T}\,f\bigr](x) = \inf_{T\in\transmatset} \bigl[Tf\bigr](x) \geq \min_{y\in\states}f(y) \geq c \in\reals
\end{equation*}
and 
\begin{equation*}
\bigl[\overline{T}\,f\bigr](x) = \sup_{T\in\transmatset} \bigl[Tf\bigr](x) \geq \min_{y\in\states}f(y) \geq c\in\reals\,.
\end{equation*}
Because $x\in\states$ is arbitrary, this implies that both $\underline{T}\,f$ and $\overline{T}\,f$ are uniformly bounded below by $c\in\reals$, and hence $\underline{T}\,f\in\gamblesextabove$ and $\overline{T}\,f\in\gamblesextabove$.
\end{proof}

In the remainder of this appendix, we implicitly use the assumption that the set $\mathcal{T}$ of transition matrices is non-empty, closed, convex, and that it has separately specified rows. We remark that closure here is understood in the operator norm topology when viewing the elements $T\in\mathcal{T}$ as operators on $\gambles$; thus, for any operator $M:\gambles\to\gambles$, we define
\begin{equation*}
\norm{M} \coloneqq \sup\Bigl\{ \norm{Mf}\,:\,f\in\gambles, \norm{f}\leq 1 \Bigr\}\,,
\end{equation*}
and $\mathcal{T}$ is assumed to be closed under this norm. We remark that a sufficient condition for this to be a norm is that $M$ is non-negatively homogeneous, meaning that $M(\lambda f)=\lambda Mf$ for any $\lambda\geq 0$ and $f\in\gambles$; this property is clearly satisfied when $M$ is a matrix, or a lower or upper transition operator---see~\lemmaref{lemma:lower_trans_basic_properties} below---so this suffices for our purposes.

Note that because $\states$ is finite, convergence of a sequence $\{M_n\}_{n\in\nats}$ of $\lvert\states\rvert\times\lvert\states\rvert$ matrices is equivalent to the element-wise convergence of $M_n(x,y)$ for all $x,y\in\states$. Moreover, note that convergence in norm implies pointwise convergence; thus if $\lim_{n\to+\infty}M_n=M$ then for any $f\in\gambles$ it also holds that $\lim_{n\to+\infty}M_nf=Mf$.

We next state a number of properties of the lower and upper transition operators corresponding to $\mathcal{T}$, which will be useful throughout this appendix. We start with an important interpretation: for any set $B\subseteq \states$, any $n\in\nats$, and any $x\in\states$, the quantity $\bigl[\underline{T}\mathbb{I}_B\bigr](x)$ can be interpreted as the \emph{lower probability} of the process being in some state $y\in B$ at time $n+1$, given that it started in state $x$ at time $n$:
\begin{equation*}
\bigl[\underline{T}\mathbb{I}_B\bigr](x) = \underline{P}\bigl(X_{n+1}\in B\,\big\vert\,X_n=x\bigr)\,.
\end{equation*}
Similarly, $\bigl[\overline{T}\mathbb{I}_B\bigr](x)$ represents the corresponding \emph{upper probability}. Throughout this appendix, we often use this interpretation to attempt to provide some intuition about the technical results that we give. 

Next, we provide some technical properties about these operators. Note that these are well known in the literature and merely stated here for convenience. We emphasise that, in the remainder of this appendix, we assume these properties to be known and will often use them without explicit mention.
\begin{lemma}\label{lemma:lower_trans_basic_properties}
For all $f,g\in\gambles$ and all $\lambda,\mu\in\reals$ with $\lambda\geq 0$, it holds that
\begin{enumerate}
\item $\underline{T}\,f \geq \min_{y\in\states} f(y)$ and $\overline{T}\,f \leq \max_{y\in\states} f(y)$
\item $f\leq g \Rightarrow \underline{T}\,f \leq \underline{T}\,g$ and $\overline{T}\,f \leq \overline{T}\,g$
\item $\underline{T}\,f + \underline{T}\,g \leq \underline{T}(f+g)$ and $\overline{T}(f+g)\leq\overline{T}\,f + \overline{T}\,g$
\item $\underline{T}(\lambda f) = \lambda\underline{T}\,f$ and $\overline{T}(\lambda f) = \lambda\overline{T}\,f$
\item $\underline{T}(f + \mu) = \underline{T}\,f + \mu$ and $\overline{T}(f + \mu) = \overline{T}\,f + \mu$
\item $\abs{\underline{T}\,f - \underline{T}\,g} \leq \overline{T}\abs{f-g}$ and $\abs{\overline{T}\,f - \overline{T}g} \leq \overline{T}\abs{f-g}$
\end{enumerate}
\end{lemma}
\begin{proof}
Fix any $x\in\states$. Then for any $T\in\mathcal{T}$, we can consider the map $\bigl[T(\cdot)\bigr](x):\gambles\to\reals$, defined for $f\in\gambles$ as
\begin{equation*}
\bigl[Tf\bigr](x) = \sum_{y\in\states} T(x,y)f(y)\,.
\end{equation*}
Because $T(x,\cdot)$ is a probability mass function, this means that we can interpret $\bigl[T(\cdot)\bigr](x):\gambles\to\reals$ as a \emph{linear prevision}---essentially an expectation operator, albeit with a slightly different interpretation---on $\gambles$; see~\cite[Section 2.8]{walley1991statistical} for details.

Now, by definition, for any $f\in\gambles$, it holds that
\begin{equation*}
\bigl[\underline{T}\,f\bigr](x) = \inf_{T\in\mathcal{T}} \bigl[Tf\bigr](x)\quad\text{and}\quad \bigl[\overline{T}\,f\bigr](x) = \sup_{T\in\mathcal{T}} \bigl[Tf\bigr](x)\,.
\end{equation*}
So, the map $\bigl[\underline{T}(\cdot)\bigr](x):\gambles\to\reals$ is a lower envelope of the linear previsions $\bigl[T(\cdot)\bigr](x)$ for which $T\in\mathcal{T}$. Similarly, the map $\bigl[\overline{T}(\cdot)\bigr](x):\gambles\to\reals$ is the upper envelope of these linear previsions. Therefore, by~\cite[Theorem 3.3.3]{walley1991statistical}, $\bigl[\underline{T}(\cdot)\bigr](x)$ is a \emph{coherent lower prevision} on $\gambles$. Similarly, $\bigl[\overline{T}(\cdot)\bigr](x)$ is a \emph{coherent upper prevision} on $\gambles$. 

The stated properties now follow (pointwise) from~\cite[Section 2.6]{walley1991statistical}.
\end{proof}

Moreover, we will need some of the above properties of $\underline{T}$ and $\overline{T}$, but generalised to their \emph{compositions} $\underline{T}^n$ and $\overline{T}^n$, with $n\in\nats$. The reason that this is possible is, essentially, that these compositions are themselves also lower and upper transition operators (but corresponding to different sets of transition matrices). Again, these properties are well-known and ubiquitous in the imprecise Markov chain literature, so we merely state them here for convenience and will subsequently often use them without explicit reference.
\begin{lemma}\label{lemma:lower_trans_composite_basic_properties}
For all $f,g\in\gambles$, all $\lambda,\mu\in\reals$ with $\lambda\geq 0$, and all $n\in\nats$, it holds that
\begin{enumerate}
	\item $\underline{T}^nf\geq \min_{y\in\states} f(y)$ and $\overline{T}^nf\leq \max_{y\in\states} f(y)$
	\item $f\leq g \Rightarrow \underline{T}^nf \leq \underline{T}^ng$ and $\overline{T}^nf \leq \overline{T}^ng$
	\item $\underline{T}^nf + \underline{T}^ng \leq \underline{T}^n(f+g)$ and $\overline{T}^n(f+g)\leq\overline{T}^nf + \overline{T}^ng$
	\item $\underline{T}^n(\lambda f) = \lambda\underline{T}^nf$ and $\overline{T}^n(\lambda f) = \lambda\overline{T}^nf$
	\item $\underline{T}^n(f + \mu) = \underline{T}^nf + \mu$ and $\overline{T}^n(f + \mu) = \overline{T}^nf + \mu$
\end{enumerate}
\end{lemma}
\begin{proof}
All of these follow straightforwardly from~\lemmaref{lemma:lower_trans_basic_properties} and induction on $n$, potentially combined with the monotonicity property (2) of~\lemmaref{lemma:lower_trans_basic_properties}. We illustrate this only for the first two properties about $\underline{T}$.

We want to show that $\underline{T}^nf\geq \min_{y\in\states} f(y)$. By~\lemmaref{lemma:lower_trans_basic_properties}, it holds that $\underline{T}f\geq \min_{y\in\states} f(y)$, so the statement is true for $n=1$. Now suppose the statement is true for $n\in\nats$. Then also
\begin{align*}
\underline{T}^{n+1}f = \underline{T}\bigl(\underline{T}^nf\bigr) &\geq \underline{T}\bigl( \min_{y\in\states}f(y) \bigr) \geq \min_{y\in\states}f(y)\,,
\end{align*}
where the first inequality used the induction hypothesis together with the monotonicity of $\underline{T}$, and the second inequality used the monotonicity of $\underline{T}$.

The second property of $\underline{T}$ can be proved in a similar manner; we want to prove that $f\leq g$ implies $\underline{T}^nf\leq \underline{T}^ng$. By~\lemmaref{lemma:lower_trans_basic_properties} this is true for $n=1$. Now suppose the statement is true for $n\in\nats$. Then, if $f\leq g$, also
\begin{align*}
\underline{T}^{n+1}f = \underline{T}\bigl(\underline{T}^nf\bigr) \leq \underline{T}\bigl(\underline{T}^ng\bigr) = \underline{T}^{n+1}g\,,
\end{align*}
where the inequality used the induction hypothesis together with the monotonicity of $\underline{T}$. 
\end{proof}

We also note that the lower and upper transition operators are \emph{conjugate} when viewed as operators on $\gambles$:
\begin{lemma}\label{lemma:transop_conjugate_gambles}
For any $f\in\gambles$, it holds that
\begin{equation*}
\underline{T}\,f = -\overline{T}(-f).
\end{equation*}
\end{lemma}
\begin{proof}
This follows immediately from the definition~\eqref{eq:def:lower_trans}.
\end{proof}
\begin{corollary}\label{lemma:transop_composite_conjugate_gambles}
	For any $f\in\gambles$ and $n\in\nats$, it holds that
	\begin{equation*}
	\underline{T}^n\,f = -\overline{T}^n(-f).
	\end{equation*}
\end{corollary}
\begin{proof}
	This is immediate from~\lemmaref{lemma:transop_conjugate_gambles} and induction on $n$.
\end{proof}

Moreover, $\overline{T}^n$ dominates $\underline{T}^n$ for any $n\in\nats$:
\begin{lemma}\label{lemma:composition_upper_dominates_composition_lower}
For any $f\in\gambles$ and $n\in\nats$, it holds that $\underline{T}^nf \leq \overline{T}^nf$.
\end{lemma}
\begin{proof}
We give a proof by induction. The fact that $\underline{T}f \leq \overline{T}f$ for any $f\in\gambles$ is immediate from the definitions of $\underline{T}$ and $\overline{T}$. 

So now assume that the result is true for some $n\in\nats$; we will show that it is also true for $n+1$. Specifically, for any $f\in\gambles$, it holds that
\begin{align*}
\underline{T}^{n+1}f = \underline{T}\underline{T}^nf \leq  \underline{T}\overline{T}^nf \leq \overline{T}\overline{T}^nf = \overline{T}^{n+1}f\,,
\end{align*}
where the first inequality used the induction hypothesis together with the monotonicity of $\underline{T}$, and the second inequality used the fact that $\underline{T}g \leq \overline{T}g$ for any $g\in\gambles$. 

\end{proof}

We need the following technical lemma:
\begin{lemma}\label{lemma:pointwise_convergence_transmats}
Consider a sequence $\{T_n\}_{n\in\nats}$ of transition matrices such that $\lim_{n\to+\infty} T_n = T$ and fix any $f\in\gamblesextabove$ and $x\in\states$. If $\bigl[T_nf\bigr](x) < +\infty$ for all $n\in\nats$, then $\lim_{n\to+\infty} \bigl[T_nf\bigr](x) = \bigl[Tf\bigr](x)$.
\end{lemma}
\begin{proof}
The condition that every step along the sequence is real-valued implies that, for all $n\in\nats$,
\begin{equation*}
f(y)=+\infty\quad\Rightarrow\quad T_n(x,y)=0\,,
\end{equation*}
which yields $T(x,y)=\lim_{n\to+\infty}T_n(x,y)=0$ if $f(y)=+\infty$.
Now let $g\in\gambles$ be such that $g(y) = f(y)$ if $f(y)<+\infty$, and $g(y)=0$, otherwise. 

Then $\bigl[Tf\bigr](x)=\bigl[Tg\bigr](x)$ and $\bigl[T_nf\bigr](x)=\bigl[T_ng\bigr](x)$ for all $n\in\nats$, so it suffices to prove that $\lim_{n\to+\infty} \bigl[T_ng\bigr](x) = \bigl[Tg\bigr](x)$. This property holds because $g\in\gambles$ and  $\lim_{n\to+\infty} T_n=T$ in the norm topology.
\end{proof}

We will next prove the various claims in~\lemmaref{lemma:trans_op_continuous} separately. We first prove the reachability properties.
\begin{lemma}\label{lemma:upper_trans_reached}
Fix any $f\in\gamblesextabove$. Then there is some $T\in\mathcal{T}$ such that $Tf=\overline{T}f$.
\end{lemma}
\begin{proof}
We first prove this pointwise per state, so fix any $x\in\states$. Now we consider two cases. First, suppose that $\bigl[\overline{T}f\bigr](x)<+\infty$. Then, for every $\epsilon>0$, we can find some $T_\epsilon\in\mathcal{T}$ such that
\begin{equation*}
\bigl[\overline{T}f\bigr](x) < \bigl[T_\epsilon f\bigr](x) + \epsilon\,.
\end{equation*}
Consider any positive sequence $\{\epsilon_n\}_{n\in\nats}$ in $\reals$ such that $\lim_{n\to+\infty} \epsilon_n=0$ and the corresponding sequence $T_{\epsilon_n}$ in $\mathcal{T}$. Because the set $\mathcal{T}$ is closed and bounded and finite dimensional, it is compact. Hence we get the existence of a subsequence $\{T_{\epsilon_{n_j}}\}_{j\in\nats}$ such that
\begin{equation*}
\lim_{j\to+\infty} T_{\epsilon_{n_j}} =: T_{\epsilon_{n_*}} \in\mathcal{T}. 
\end{equation*}
Because $\bigl[\overline{T}f\bigr](x)<+\infty$, it holds for any $j\in\nats$ that
\begin{equation*}
\bigl[T_{\epsilon_{n_j}} f\bigr](x) \leq \bigl[\overline{T}f\bigr](x) < +\infty\,.
\end{equation*}
Hence, we can apply~\lemmaref{lemma:pointwise_convergence_transmats} to find
\begin{align*}
\bigl[\overline{T}f\bigr](x) &\leq \lim_{j\to+\infty} \left( \bigl[T_{\epsilon_{n_j}}f\bigr](x) + \epsilon_{n_j}\right) \\
 &= \lim_{j\to+\infty} \bigl[T_{\epsilon_{n_j}}f\bigr](x) = \bigl[T_{\epsilon_{n_*}}f\bigr](x)\,.
\end{align*}
Conversely, $T_{\epsilon_{n_*}}\in\mathcal{T}$ directly implies
\begin{equation*}
\bigl[T_{\epsilon_{n_*}}f\bigr](x) \leq \bigl[\overline{T}f\bigr](x)\,,
\end{equation*}
so we conclude that $\bigl[T_{\epsilon_{n_*}}f\bigr](x)=\bigl[\overline{T}f\bigr](x)$.

For the other case, suppose that $\bigl[\overline{T}f\bigr](x)=+\infty$. This means that, for any $c\in \reals$, there is some $T_c\in\mathcal{T}$ such that $\bigl[T_cf\bigr](x)>c$. So consider any sequence $\{c_n\}_{n\in\nats}$ such that $\lim_{n\to+\infty} c_n=+\infty$ and the corresponding sequence $T_{c_n}\in\mathcal{T}$. If there is some $n\in\nats$ such that $\bigl[T_{c_n}f\bigr](x)=+\infty = \bigl[\overline{T}f\bigr](x)$ then we are done with this second case.

So suppose that $\bigl[T_{c_n}f\bigr](x)<+\infty$ for all $n\in\nats$. Because $\mathcal{T}$ is compact, we get a convergent subsequence $\{T_{c_{n_j}}\}_{j\in\nats}$ with $\lim_{j\to+\infty}T_{c_{n_j}} =: T_{c_{n_*}}\in\mathcal{T}$. Moreover, by~\lemmaref{lemma:pointwise_convergence_transmats}, we have
\begin{equation*}
\bigl[T_{c_{n_*}}f\bigr](x) = \lim_{j\to+\infty} \bigl[T_{c_{n_j}}f\bigr](x) \geq \lim_{j\to+\infty} c_{n_j} = +\infty \,.
\end{equation*}
This concludes the proof for the second case.

To summarise the above, for any $x\in\states$ we can find some $T_x\in\mathcal{T}$ such that $\bigl[T_xf\bigr](x) = \bigl[\overline{T}f\bigr](x)$. Now consider the matrix $T$, defined for all $x\in\states$ as $T(x,\cdot)\coloneqq T_x(x,\cdot)$. Then $T\in\mathcal{T}$ because $\mathcal{T}$ has separately specified rows. Moreover, for any $x\in\states$,
\begin{align*}
\bigl[Tf\bigr](x) &= \sum_{y\in\states} T(x,y)f(y) \\
 &= \sum_{y\in\states} T_x(x,y)f(y) = \bigl[T_xf\bigr](x) = \bigl[\overline{T}f\bigr](x)\,,
\end{align*}
and so $Tf=\overline{T}f$.
\end{proof}

\begin{lemma}\label{lemma:lower_trans_reached}
Fix any $f\in\gamblesextabove$. Then there is some $T\in\mathcal{T}$ such that $Tf = \underline{T}\,f$.
\end{lemma}
\begin{proof}
We first prove this pointwise per state, so fix any $x\in\states$. Now we consider two cases. First, suppose that $\bigl[\underline{T}f\bigr](x)=+\infty$. Then for any $T\in\mathcal{T}$ it holds that
\begin{equation*}
+\infty = \bigl[\underline{T}\,f\bigr](x) \leq \bigl[Tf\bigr](x)\,,
\end{equation*}
which implies that also $\bigl[Tf\bigr](x)=+\infty$. This concludes the first case.

For the other case, suppose that $\bigl[\underline{T}\,f\bigr](x)<+\infty$. Then, for all $\epsilon>0$, we can find some $T_\epsilon\in\mathcal{T}$ such that
\begin{equation*}
\bigl[T_\epsilon\,f\bigr](x) - \epsilon < \bigl[\underline{T}\,f\bigr](x)\,.
\end{equation*}
Consider any positive sequence $\{\epsilon_n\}_{n\in\nats}$ in $\reals$ such that $\lim_{n\to+\infty} \epsilon_n=0$ and the corresponding sequence $T_{\epsilon_n}$ in $\mathcal{T}$. Because the set $\mathcal{T}$ is closed and bounded and finite dimensional, it is compact. Hence we get the existence of a subsequence $\{T_{\epsilon_{n_j}}\}_{j\in\nats}$ such that
\begin{equation*}
\lim_{j\to+\infty} T_{\epsilon_{n_j}} =: T_{\epsilon_{n_*}} \in\mathcal{T}. 
\end{equation*}
Because $\bigl[\underline{T}f\bigr](x)<+\infty$, it holds for any $j\in\nats$ that
\begin{equation*}
\bigl[T_{\epsilon_{n_j}} f\bigr](x) - \epsilon_{n_j} < \bigl[\underline{T}f\bigr](x) < +\infty\,,
\end{equation*}
from which it follows that also $\bigl[T_{\epsilon_{n_j}} f\bigr](x)<+\infty$.

Hence, we can apply~\lemmaref{lemma:pointwise_convergence_transmats} to find
\begin{align*}
\bigl[\underline{T}f\bigr](x) &\geq \lim_{j\to+\infty} \left( \bigl[T_{\epsilon_{n_j}}f\bigr](x) - \epsilon_{n_j} \right) \\
 &= \lim_{j\to+\infty} \bigl[T_{\epsilon_{n_j}}f\bigr](x)  =  \bigl[T_{\epsilon_{n_*}}f\bigr](x)\,.
\end{align*}
Conversely, $T_{\epsilon_{n_*}}\in\mathcal{T}$ directly implies
\begin{equation*}
\bigl[T_{\epsilon_{n_*}}f\bigr](x) \geq \bigl[\underline{T}f\bigr](x)\,,
\end{equation*}
so we conclude that $\bigl[T_{\epsilon_{n_*}}f\bigr](x)=\bigl[\underline{T}\,f\bigr](x)$. This concludes the second case.

To summarise the above, for any $x\in\states$ we can find some $T_x\in\mathcal{T}$ such that $\bigl[T_xf\bigr](x) = \bigl[\underline{T}\,f\bigr](x)$. Now consider the matrix $T$, defined for all $x\in\states$ as $T(x,\cdot)\coloneqq T_x(x,\cdot)$. Then $T\in\mathcal{T}$ because $\mathcal{T}$ has separately specified rows. Moreover, for any $x\in\states$,
\begin{align*}
\bigl[Tf\bigr](x) &= \sum_{y\in\states} T(x,y)f(y) \\
 &= \sum_{y\in\states} T_x(x,y)f(y) = \bigl[T_xf\bigr](x) = \bigl[\underline{T}\,f\bigr](x)\,,
\end{align*}
and so $Tf=\underline{T}\,f$.
\end{proof}

We next need a number of additional technical properties of lower and upper transition operators, before we can prove the remaining statements in~\lemmaref{lemma:trans_op_continuous}.
\begin{lemma}\label{lemma:trans_op_monotone}
Consider any $f,g\in\gamblesextabove$ such that $f\leq g$. Then it holds that
\begin{equation*}
\underline{T}\,f\leq \underline{T}\,g\quad\text{and}\quad \overline{T}f\leq \overline{T}g\,.
\end{equation*}
\end{lemma}
\begin{proof}
Using~\lemmaref{lemma:upper_trans_reached}, we find $T\in\mathcal{T}$ such that $Tf=\overline{T}f$. Now fix any $x\in\states$. Then because $T\in\mathcal{T}$ it holds that
\begin{align*}
\bigl[\overline{T}f\bigr](x) &= \bigl[Tf\bigr](x) \\
 &= \sum_{y\in\states} T(x,y)f(y) \\
 &\leq \sum_{y\in\states} T(x,y)g(y) \\
 &= \bigl[Tg\bigr](x) \leq \bigl[\overline{T}g\bigr](x)\,,
\end{align*}
where the first inequality used the fact that $T(x,y)\geq 0$ for all $y\in\states$. This implies $\overline{T}\,f\leq\overline{T}\,g$ because $x$ was arbitrary.

We next prove the statement for $\underline{T}$. Using~\lemmaref{lemma:lower_trans_reached}, we find a $T\in\mathcal{T}$ such that $Tg=\underline{T}\,g$. Then
\begin{align*}
\bigl[\underline{T}g\bigr](x) &= \bigl[Tg\bigr](x) \\
 &= \sum_{y\in\states} T(x,y)g(y) \\
 &\geq \sum_{y\in\states} T(x,y)f(y) \\
 &= \bigl[Tf\bigr](x) \geq \bigl[\underline{T}f\bigr](x)\,,
\end{align*}
where the first inequality used the fact that $T(x,y)\geq 0$ for all $y\in\states$.
This implies $\underline{T}\,f\leq\underline{T}\,g$ because $x$ was arbitrary.

\end{proof}

\begin{lemma}\label{lemma:const_add_extgambles}
For any $f\in\gamblesextabove{}$ and any $\mu\in\reals$, it holds that
\begin{align*}
\underline{T}(f + \mu) = \underline{T}\,f + \mu \quad \text{ and } \quad \overline{T}(f + \mu) = \overline{T}\,f + \mu.
\end{align*}
\end{lemma}
\begin{proof}
Trivially from the definition of $\underline{T}$ and $\overline{T}$.
\end{proof}

\begin{lemma}\label{lemma:nneg_homogen_extgambles}
For any non-negative $f\in\gamblesextabove{}$ and any non-negative $\lambda\in\realsextabove{}$, it holds that
$\underline{T}(\lambda f) = \lambda \underline{T}\,f$.
\end{lemma}
\begin{proof}
Fix any $x\in\states$. Note that for any $T\in\transmatset$ we have that
\begin{equation}\label{lemma:nneg_homogen_extgambles:eq:linear}
\begin{split}
\bigl[T(\lambda f)\bigr](x) &= \sum_{y\in\states} T(x,y)\lambda f(y) \\
 &= \lambda \sum_{y\in\states}T(x,y)f(y) = \lambda \bigl[Tf\bigr](x)\,.
\end{split}
\end{equation}	
Observe that the non-negativity of $f$ suffices to ensure that everything is well-defined; in particular the sums do not contain terms that equal $(-\infty)$.


Using~\lemmaref{lemma:lower_trans_reached}, we can find $T,S\in\transmatset$ such that
\begin{equation*}
\bigl[T(\lambda f)\bigr](x) = \bigl[\underline{T}(\lambda f)\bigr](x)\,\,\text{and}\,\, \bigl[Sf\bigr](x) = \bigl[\underline{T}f\bigr](x)\,.
\end{equation*}
The choice of $T\in\mathcal{T}$ then implies that
\begin{equation*}
\lambda \bigl[Tf\bigr] = \bigl[T(\lambda f)\bigr] \leq  \bigl[S(\lambda f)\bigr] = \lambda \bigl[Sf\bigr]\,,
\end{equation*}
where we used~\eqref{lemma:nneg_homogen_extgambles:eq:linear} for the two equalities. Conversely, the choice of $S\in\transmatset$ implies that $\bigl[Sf\bigr](x) \leq \bigl[Tf\bigr](x)$ which, by multiplying both sides with $\lambda$, yields
\begin{equation*}
\lambda \bigl[Sf\bigr](x) \leq \lambda \bigl[Tf\bigr](x)\,,
\end{equation*}
using the non-negativity of $\lambda$.
Hence we can conclude that
\begin{align*}
\bigl[\underline{T}(\lambda f)\bigr](x) = \bigl[T(\lambda f)\bigr](x) &= \lambda \bigl[Tf\bigr](x) \\
  &= \lambda \bigl[Sf\bigr](x) = \lambda \bigl[\underline{T} f\bigr](x)\,.
\end{align*}
Since $x\in\states$ was arbitrary, the claim follows.
	


\end{proof}

\begin{lemma}\label{lemma:norm_properties_transop}
For any $f,g\in\gambles$, it holds that
\begin{equation*}
\norm{\underline{T}\,f - \underline{T}\,g} \leq \norm{f-g}\quad\text{and}\quad\norm{\overline{T}\,f - \overline{T}\,g} \leq \norm{f-g}\,.
\end{equation*}
\end{lemma}
\begin{proof}
Fix any $x\in\states$. Then, using~\lemmaref{lemma:lower_trans_basic_properties}, we get
\begin{align*}
\abs{\bigl[\underline{T}\,f\bigr](x) - \bigl[\underline{T}\,g\bigr](x)} &\leq \bigl[\overline{T}\abs{f-g}\bigr](x) \\
 &\leq \max_{y\in\states} \abs{f(y)-g(y)} = \norm{f-g}\,.
\end{align*}
Because $x$ was arbitrary we find that
\begin{equation*}
\norm{\underline{T}\,f - \underline{T}\,g} = \max_{x\in\states} \abs{\bigl[\underline{T}\,f\bigr](x) - \bigl[\underline{T}\,g\bigr](x)} \leq \norm{f-g}\,.
\end{equation*}
The proof for $\overline{T}$ is completely analogous.
\end{proof}

We can now prove that $\underline{T}$ and $\overline{T}$ are continuous on $\gambles$:
\begin{lemma}\label{lemma:transop_cont_gambles}
Consider any sequence $\{f_n\}_{n\in\nats}$ in $\gambles$ such that $\lim_{n\to+\infty} f_n=f\in\gambles$. Then
\begin{equation*}
\lim_{n\to+\infty} \underline{T}\,f_n = \underline{T}\,f\quad\text{and}\quad \lim_{n\to+\infty} \overline{T}f_n = \overline{T}f\,.
\end{equation*}
\end{lemma}
\begin{proof}
Since $\gambles$ has the norm topology, it holds that $\lim_{n\to+\infty}\norm{f_n-f} = 0$. So then, using~\lemmaref{lemma:norm_properties_transop}, it holds that
\begin{equation*}
\lim_{n\to+\infty} \norm{\underline{T}\,f_n - \underline{T}\,f} \leq \lim_{n\to+\infty} \norm{f_n - f} = 0\,,
\end{equation*}
from which we conclude that $\lim_{n\to+\infty} \underline{T}\,f_n = \underline{T}\,f$. The proof for $\overline{T}$ is completely analogous.
\end{proof}

It remains to prove that $\underline{T}$ and $\overline{T}$ are continuous on $\gamblesextabove$ with respect to non-decreasing sequences:
\begin{lemma}\label{lemma:transop_cont_ext_gambles_nondec}
Consider any non-decreasing sequence $\{f_n\}_{n\in\nats}$ in $\gamblesextabove$ such that $\lim_{n\to+\infty} f_n=f\in\gamblesextabove$. Then
\begin{equation*}
\lim_{n\to+\infty} \underline{T}\,f_n = \underline{T}\,f\quad\text{and}\quad \lim_{n\to+\infty} \overline{T}f_n = \overline{T}f\,.
\end{equation*}
\end{lemma}
\begin{proof}
We first remark that because $f_1 \in \gamblesextabove$ and $\states{}$ is finite, $f_1$ is bounded below.
Hence, because of this, we can use Lemma~\ref{lemma:const_add_extgambles}, and assume without loss of generality that $f_1$, and therefore also $f$ and all $f_n$, are non-negative; essentially we ``shift'' all functions to be non-negative, prove that the convergence holds, and then shift them back afterwards. 

We now start by proving the claim for $\overline{T}$, and prove the claim pointwise per state. So, fix any $x\in\states$. We now consider two cases. First, suppose that $\bigl[\overline{T}f\bigr](x)<+\infty$. By~\lemmaref{lemma:trans_op_monotone} and the fact that the sequence $f_n$ is non-decreasing, this implies that also $\bigl[\overline{T}f_n\bigr](x)<+\infty$ for all $n\in\nats$. We can then find functions $g_n,g\in\gambles$ such that for all $n\in\nats$, $g_n(y)=f_n(y)$ and $g(y)=f(y)$ if $f(y)<+\infty$, and $g_n(y)=g(y)=0$ otherwise. Then clearly $\lim_{n\to+\infty} g_n=g$, and, as we will show next,
\begin{equation*}
\bigl[\overline{T}f\bigr](x) = \bigl[\overline{T}g\bigr](x)\quad\text{and}\quad \bigl[\overline{T}f_n\bigr](x) = \bigl[\overline{T}g_n\bigr](x)\,,
\end{equation*}
for all $n\in\nats$. To see this, note that $\bigl[\overline{T}f\bigr](x)<+\infty$ implies, for all $T\in\mathcal{T}$, that $T(x,y)=0$ for all $y\in\states$ such that $f(y)=+\infty$. Hence, $\bigl[Tf\bigr](x) = \bigl[Tg\bigr](x)$ and $\bigl[Tf_n\bigr](x) = \bigl[Tg_n\bigr](x)$ for all $T\in\mathcal{T}$, from which the claim follows.

Hence, it suffices to note that $\lim_{n\to+\infty}\bigl[\overline{T}g_n\bigr](x)=\bigl[\overline{T}g\bigr](x)$, which follows immediately from~\lemmaref{lemma:transop_cont_gambles}.

For the other case, suppose that $\bigl[\overline{T}f\bigr](x)=+\infty$. Using~\lemmaref{lemma:upper_trans_reached}, we find a $T\in\mathcal{T}$ such that $\bigl[Tf\bigr](x)=\bigl[\overline{T}f\bigr](x)=+\infty$. This implies that there is some $y\in\states$ such that $T(x,y)>0$ and $f(y)=+\infty$. Because $\lim_{n\to+\infty} f_n(y)=f(y)=+\infty$ and the sequence $\{f_n\}_{n\in\nats}$ is  non-negative and non-decreasing, this implies that
\begin{equation*}
\lim_{n\to+\infty} \bigl[Tf_n\bigr](x) \geq \lim_{n\to+\infty} T(x,y)f_n(y) = +\infty\,. 
\end{equation*}
Because $T\in\mathcal{T}$ it holds for all $n\in\nats$ that $\bigl[Tf_n\bigr](x)\leq \bigl[\overline{T}f_n\bigr](x)$, so we find that
\begin{equation*}
+\infty = \lim_{n\to+\infty} \bigl[Tf_n\bigr](x)\leq \lim_{n\to+\infty}\bigl[\overline{T}f_n\bigr](x)\,,
\end{equation*}
and hence $\lim_{n\to+\infty}\bigl[\overline{T}f_n\bigr](x)=+\infty = \bigl[\overline{T}f\bigr](x)$. This concludes the proof of the second case.

In summary, we have shown that $\lim_{n\to+\infty}\bigl[\overline{T}f_n\bigr](x) = \bigl[\overline{T}f\bigr](x)$ for all $x\in\states$. Because $\gamblesextabove$ has the topology of pointwise convergence, it follows that $\lim_{n\to+\infty} \overline{T}f_n = \overline{T} f$. This concludes the proof for $\overline{T}$.

It remains to prove the claim for $\underline{T}$. 
Note that because of Lemma \ref{lemma:trans_op_monotone} and the non-decreasing character of $\{f_n\}_{n \in \nats{}}$, we have that $ \underline{T}\,f_n \leq \underline{T}\, f$ for all $n \in \nats{}$, and therefore $\lim_{n\to+\infty} \underline{T}\,f_n$ exists and $\lim_{n\to+\infty} \underline{T}\,f_n \leq \underline{T}\, f$.
So it remains to prove $\lim_{n\to+\infty} \underline{T}\,f_n \geq \underline{T}\, f$.

Fix any $x\in\states$ and let $\mathcal{Y} \coloneqq \{y \in \states{} \colon f(y) = +\infty\}$. 
Then according to Lemma~\ref{lemma:lower_trans_reached} there is a $T \in \mathcal{T}$ such that
\begin{align*}
[\underline{T}\mathbb{I}_{\mathcal{Y}}](x) 
= [{T}\mathbb{I}_{\mathcal{Y}}](x). 
\end{align*}
Then, according to Lemma~\ref{lemma:nneg_homogen_extgambles}, we also have, for all non-negative $\alpha\in\realsextabove$, that
\begin{align}\label{Eq: proof of lemma:transop_cont_ext_gambles_nondec}
[\underline{T}(\alpha \mathbb{I}_{\mathcal{Y}})](x) 
&= \alpha [\underline{T} \mathbb{I}_{\mathcal{Y}}](x) 
= \alpha [{T}\mathbb{I}_{\mathcal{Y}}](x).
\end{align}
We consider two cases; first, suppose that $[{T}\mathbb{I}_{\mathcal{Y}}](x) > 0$. Using~\lemmaref{lemma:trans_op_monotone} and the non-negativity of $f$, it then holds that
\begin{align*}
[\underline{T}\,f](x) 
&\geq [\underline{T}(+\infty \, \mathbb{I}_\mathcal{Y})](x)
= +\infty [T \mathbb{I}_\mathcal{Y}](x) = +\infty,
\end{align*}
where we used that $f = +\infty\mathbb{I}_\mathcal{Y} + f\hprod\mathbb{I}_{\mathcal{Y}^c} \geq +\infty \mathbb{I}_\mathcal{Y}$.
Hence, $[\underline{T}\,f](x) = +\infty$.
Similarly, for any $n \in \nats{}$, we also have that
\begin{align*}
[\underline{T}\,f_n](x) 
&\geq [\underline{T}(f_n \hprod\mathbb{I}_\mathcal{Y})](x)  \\
&\geq [\underline{T}(\min_{y \in \mathcal{Y}} f_n(y) \mathbb{I}_\mathcal{Y})](x)  \\
&= \min_{y \in \mathcal{Y}} f_n(y) [T \mathbb{I}_\mathcal{Y}](x) ,
\end{align*}
where the first two steps follow from Lemma \ref{lemma:trans_op_monotone} and the non-negativity of $f_n$, and the third from the non-negativity of $f_n$ together with Equation \eqref{Eq: proof of lemma:transop_cont_ext_gambles_nondec}.
Since $\lim_{n \to +\infty} f_n = f$, we have that $\lim_{n \to +\infty}f_n(y) = +\infty$ for all $y \in \mathcal{Y}$, and therefore, because $\mathcal{Y}$ is finite, that
\begin{align*}
\lim_{n \to +\infty} \min_{y \in \mathcal{Y}} f_n(y) = +\infty.
\end{align*}
Hence, because $[{T}\mathbb{I}_{\mathcal{Y}}](x) > 0$ by our current assumption, we find that
\begin{align*}
&\lim_{n \to +\infty} [\underline{T}\,f_n](x) \geq \lim_{n \to +\infty}  \min_{y \in \mathcal{Y}} f_n(y) [T \mathbb{I}_\mathcal{Y}](x) \geq +\infty\,.
\end{align*}
So we indeed have that $[\underline{T}\,f](x) = +\infty = \lim_{n \to +\infty} [\underline{T}\,f_n](x)$.

For the second case, suppose that  $[{T}\mathbb{I}_{\mathcal{Y}}](x) = 0$.
Then, according to Lemma~\ref{lemma:lower_trans_reached}, there is a sequence $\{T_n\}_{n \in \nats{}}$ in $\mathcal{T}$ such that $[\underline{T}\,f_n](x) = [T_n f_n](x)$ for all $n \in \nats{}$.
Because $\mathcal{T}$ is closed and bounded and finite-dimensional, it is compact. Hence, we get a convergent subsequence $\{T_{n_j}\}_{j\in\nats}$ such that $\lim_{j\to+\infty} T_{n_j}=:T_{n_*}\in\mathcal{T}$.
Moreover, we have that
\begin{align}\label{Eq: proof of lemma:transop_cont_ext_gambles_nondec 2}
\lim_{n \to +\infty} [\underline{T}\,f_n](x) 
&= \lim_{n \to +\infty} [T_n f_n](x) \nonumber \\
&= \lim_{j \to +\infty} [T_{n_j} f_{n_j}](x) \nonumber \\
&\geq \limsup_{j \to +\infty} \, [T_{n_j} (f_{n_j}\hprod \mathbb{I}_{\mathcal{Y}^c})](x) \nonumber \\
&= \limsup_{j \to +\infty} \sum_{y \in \mathcal{Y}^c} T_{n_j}(x,y) f_{n_j}(y),
\end{align}
where we have used the non-negativity of $f_{n_j}$ in the third step.
Furthermore, because $\mathcal{Y}^c$ is finite and $\lim_{j \to +\infty}f_{n_j}(y) = f(y) \in \reals{}$ and $\lim_{j \to +\infty} T_{n_j}(x,y) = T_{n_*}(x,y) \in \reals{}$ for all $y \in \mathcal{Y}^c$, it follows from Equation~\eqref{Eq: proof of lemma:transop_cont_ext_gambles_nondec 2} that
\begin{align}\label{Eq: proof of lemma:transop_cont_ext_gambles_nondec 3}
\lim_{n \to +\infty} [\underline{T}\,f_n](x) 
&\geq  \sum_{y \in \mathcal{Y}^c} \lim_{j \to +\infty} T_{n_j}(x,y) \lim_{j \to +\infty} f_{n_j}(y) \nonumber \\
&=  \sum_{y \in \mathcal{Y}^c} T_{n_*}(x,y) f(y)\,.
\end{align}


Next, note that 
\begin{equation*}
[\underline{T}\,f](x) \leq [T f](x) = [T (f\hprod \mathbb{I}_{\mathcal{Y}^c})](x)<+\infty,
\end{equation*}
where the equality follows from $[T \mathbb{I}_{\mathcal{Y}}](x) = 0$ and our convention that $0 \cdot + \infty = 0$, and the last inequality follows from the boundedness (both above and below) of $f\hprod \mathbb{I}_{\mathcal{Y}^c}$, which is a consequence of the finiteness of $\states{}$.
Hence, because also $\lim_{n \to +\infty} [\underline{T}\,f_n](x) \leq [\underline{T}\,f](x)$, it follows that $\lim_{n \to +\infty} [\underline{T}\,f_n](x) < +\infty$, which implies that
\begin{align*}
+\infty > \lim_{n \to +\infty} [\underline{T}\,f_n](x) 
&= \lim_{n \to +\infty} [T_n f_n](x) \\
&= \lim_{j \to +\infty} [T_{n_j} f_{n_j}](x) \\
&\geq \limsup_{j \to +\infty} \, [T_{n_j} (f_{n_j}\hprod \mathbb{I}_{\mathcal{Y}})](x) \\
&= \limsup_{j \to +\infty} \sum_{y \in \mathcal{Y}} T_{n_j}(x,y) f_{n_j}(y)\,,
\end{align*}
where we have used the non-negativity of $f_{n_j}$ in the second inequality. Since $\lim_{j \to +\infty} f_{n_j}(y) = +\infty$ for all $y \in \mathcal{Y}$, this implies that $\lim_{j \to +\infty} T_{n_j}(x,y) = T_{n_*}(x,y) = 0$ for all $y \in \mathcal{Y}$.
Hence, 
using Equation~\eqref{Eq: proof of lemma:transop_cont_ext_gambles_nondec 3} we find that
\begin{align*}
\lim_{n \to +\infty} [\underline{T}\,f_n](x) 
&\geq \sum_{y\in\mathcal{Y}^c} T_{n_*}(x,y)f(y) \\
 &= \sum_{y\in\mathcal{Y}^c} T_{n_*}(x,y)f(y) + \sum_{y\in\mathcal{Y}}T_{n_*}(x,y)(+\infty) \\
&= [T_{n_*}f](x) \geq [\underline{T}\,f](x),
\end{align*}
where the last step follows from the definition of $\underline{T}$. 

In summary, we have shown that, in either case, $\lim_{n\to+\infty}\bigl[\underline{T}\,f_n\bigr](x)=\bigl[\underline{T}\,f\bigr](x)$. Because $x$ was arbitrary and because $\gamblesextabove$ has the topology of pointwise convergence, it follows that $\lim_{n\to+\infty} \underline{T}\,f_n=\underline{T}\,f$. This concludes the proof for $\underline{T}$.

\end{proof}

\begin{proofof}{\lemmaref{lemma:trans_op_continuous}}
These properties are proved separately in~\lemmaref{lemma:upper_trans_reached},~\lemmaref{lemma:lower_trans_reached},~\lemmaref{lemma:transop_cont_gambles}, and~\lemmaref{lemma:transop_cont_ext_gambles_nondec}.
\end{proofof}

\begin{lemma}\label{lemma: local game-theoretic upper exp}
	Let $\mathcal{P}$ be any closed and convex set of probability mass functions on $\states{}$, and let $\overline{\mathbb{E}}_{\mathcal{P}}$ be defined by 
	\begin{align*}
	\overline{\mathbb{E}}_{\mathcal{P}}(f) 
	\coloneqq \sup_{P \in \mathcal{P}} \mathbb{E}_{P}(f)
	\coloneqq \sup_{P \in \mathcal{P}} \sum_{x \in \states{}} f(x) P(x),
	\end{align*}
	for all $f\in\gamblesextabove$.
	Then $\overline{\mathbb{E}}_{\mathcal{P}}$ satisfies \ref{coherence: const is const}--\ref{coherence: monotonicity}.
\end{lemma}
\begin{proof}
\begin{itemize}
\item[E1.]
That $\overline{\mathbb{E}}_{\mathcal{P}}$ satisfies \ref{coherence: const is const} is trivial; it follows immediately from its definition and the fact that $\sum_{x \in \states{}} P(x) = 1$ for any probability mass function $P$ on $\states{}$.

\item[E2.]
Fix any $f$ and $g$ in $\gamblesextabove{}$.
We first show that $[f(x)+g(x)]P(x) = f(x)P(x) + g(x)P(x)$ for any probability mass function $P$ on $\states{}$ and any $x \in \states{}$.

So fix any probability mass function $P$ on $\states{}$ and any $x \in \states{}$.
If $P(x)=0$, then, by our convention that $0\cdot(+\infty) = 0$, we have that both $[f(x)+g(x)]P(x) = 0$ and $f(x)P(x) + g(x)P(x) = 0$, implying that $[f(x)+g(x)]P(x) = f(x)P(x) + g(x)P(x)$.
So we can assume that $P(x)>0$.
Then because both $f\in\gamblesextabove{}$ and $g\in\gamblesextabove{}$, we can distinguish two cases.
Either we have that $[f(x)+g(x)] \in \reals$, implying that both $f(x)\in\reals$ and $g(x)\in\reals$ and therefore directly that $[f(x)+g(x)]P(x) = f(x)P(x) + g(x)P(x)$.
Otherwise, we have that $[f(x)+g(x)]=+\infty$, implying that at least $f(x)$ or $g(x)$ equals $+\infty$. 
Hence, because $P(x)>0$, we get that both $[f(x)+g(x)]P(x)=+\infty$ and $f(x)P(x) + g(x)P(x)=+\infty$, again implying that $[f(x)+g(x)]P(x) = f(x)P(x) + g(x)P(x)$.

So we indeed have that $[f(x)+g(x)]P(x) = f(x)P(x) + g(x)P(x)$ for any probability mass function $P$ on $\states{}$ and any $x \in \states{}$.
Then, for any probability mass function $P$ on $\states{}$, we can write that
\begin{align*}
\sum_{x \in \states{}} [f(x)+g(x)]P(x) 
= \sum_{x \in \states{}} f(x)P(x) + g(x)P(x),
\end{align*}
implying, by definition, that $\mathbb{E}_{P}(f+g) = \mathbb{E}_{P}(f) + \mathbb{E}_{P}(g)$.
Now consider any set $\mathcal{P}$ of probability mass function on $\states{}$.
Then we indeed have that
\begin{align*}
\overline{\mathbb{E}}_{\mathcal{P}}(f+g) 
&= \sup_{P \in \mathcal{P}} \mathbb{E}_{P}(f+g) \\
&= \sup_{P \in \mathcal{P}} \bigl[ \mathbb{E}_{P}(f) + \mathbb{E}_{P}(g) \bigr] \\
&\leq \sup_{P \in \mathcal{P}} \mathbb{E}_{P}(f) + \sup_{P \in \mathcal{P}} \mathbb{E}_{P}(g) \\
&= \overline{\mathbb{E}}_{\mathcal{P}}(f) + \overline{\mathbb{E}}_{\mathcal{P}}(g).
\end{align*}

\item[E3.]
Fix any positive $\lambda\in\realsextabove{}$ and any non-negative $f\in\gamblesextabove{}$.
Note that for any $P\in\mathcal{P}$ we have that
\begin{equation}\label{lemma: local game-theoretic upper exp:eq:linear_is_homogeneous}
\begin{split}
\mathbb{E}_P(\lambda f) &= \sum_{x\in\states} \lambda f(x)P(x) \\
&= \lambda \sum_{x\in\states} f(x)P(x) = \lambda \mathbb{E}_P(f)\,.
\end{split}
\end{equation}	
Observe that the non-negativity of $f$ suffices to ensure that everything is well-defined; in particular the sums do not contain terms that equal $(-\infty)$.

We now consider several cases. First, if $\lambda \neq +\infty$ and $\overline{\mathbb{E}}_\mathcal{P}(f) \neq +\infty$, meaning that everything is real-valued, then we trivially have that $\overline{\mathbb{E}}_\mathcal{P}(\lambda f) = \lambda\overline{\mathbb{E}}_\mathcal{P}( f)$ using the non-negativity of $\lambda$ and the definition of the supremum.

Next, suppose that $\overline{\mathbb{E}}_\mathcal{P}(f) =0$. Because $f$ is non-negative, this implies that $\mathbb{E}_P(f)=0$ for all $P\in\mathcal{P}$, from which we get that
\begin{equation*}
\lambda \overline{\mathbb{E}}_\mathcal{P}(f) = 0 = \sup_{P\in\mathcal{P}}\lambda\mathbb{E}_P( f) = \sup_{P\in\mathcal{P}}\mathbb{E}_P(\lambda f) \,, 
\end{equation*}
regardless of the value of $\lambda$, and where we used~\eqref{lemma: local game-theoretic upper exp:eq:linear_is_homogeneous} for the final equality.

Next, suppose that $\lambda \neq +\infty$ but $\overline{\mathbb{E}}_\mathcal{P}(f) =+\infty$. Because $\lambda$ is positive this implies that also $\lambda\overline{\mathbb{E}}_\mathcal{P}(f) =+\infty$. Now fix any $c\in\reals$. Because $\overline{\mathbb{E}}_\mathcal{P}(f) =+\infty$, there is some $P\in\mathcal{P}$ such that $\mathbb{E}_P(f) \geq  \nicefrac{c}{\lambda}$, where we note that $\lambda\neq 0$ because it is positive. We therefore have
\begin{equation*}
\overline{\mathbb{E}}_\mathcal{P}(\lambda f) \geq \mathbb{E}_P(\lambda f) = \lambda\mathbb{E}_P(f) \geq c\,,
\end{equation*}
where we used~\eqref{lemma: local game-theoretic upper exp:eq:linear_is_homogeneous} for the equality.
Because $c\in\reals$ was arbitrary, this implies that $\overline{\mathbb{E}}_\mathcal{P}(\lambda f)=+\infty$. Hence also in this case the result holds. 

For the final case, suppose that $\lambda=+\infty$ and $\overline{\mathbb{E}}_\mathcal{P}(f) > 0$. Then on the one hand, the second condition implies that there exists some $P\in\mathcal{P}$ for which $\mathbb{E}_P(f)>0$, whence $\lambda\mathbb{E}_P(f)=+\infty$, and therefore, using ~\eqref{lemma: local game-theoretic upper exp:eq:linear_is_homogeneous}, also
\begin{equation*}
\overline{\mathbb{E}}_\mathcal{P}(\lambda f) \geq \mathbb{E}_P(\lambda f) = \lambda \mathbb{E}_P(f) = +\infty\,.
\end{equation*}
On the other hand, the assumptions immediately imply that $\lambda \overline{\mathbb{E}}_\mathcal{P}(f)=+\infty$.

Because this exhaustively covers all cases, we conclude that indeed $\lambda \overline{\mathbb{E}}_\mathcal{P}(f) =  \overline{\mathbb{E}}_\mathcal{P}(\lambda f)$ for all positive $\lambda\in\realsextabove$ and non-negative $f\in\gamblesextabove$.



\item[E4.]
That \ref{coherence: monotonicity} holds, follows trivially from the fact that the supremum operator is monotone and that $P(x) \geq 0$ for any probability mass function $P$ on $\states{}$ and any $x\in\states{}$.
\end{itemize}
\end{proof}

\begin{proofof}{\propositionref{prop: vovk local precise}}
For any $n\in\natswith$ and $x_0,\ldots,x_n\in\states$, the (singleton) set
\begin{equation*}
\mathcal{P} \coloneqq \Bigl\{ P(X_{n+1}\,\vert\,X_{0:n}=x_{0:n}) \Bigr\}
\end{equation*}
is a closed and convex set of probability mass functions on $\states$, and, for all $f\in\gamblesextabove$,
\begin{align*}
\overline{Q}_{x_{0:n}}(f) &= \mathbb{E}_P\bigl[f(X_{n+1})\,\vert\,X_{0:n}=x_{0:n}\bigr] \\
&= \sup_{Q \in \mathcal{P}} \mathbb{E}_{Q}\bigl[f(X_{n+1})\,\vert\,X_{0:n}=x_{0:n}\bigr]\,.
\end{align*}
Therefore, $\overline{Q}_{x_{0:n}}$ satisfies \ref{coherence: const is const}--\ref{coherence: monotonicity} due to~\lemmaref{lemma: local game-theoretic upper exp}.
\end{proofof}

Now, we will need a small notational trick to prove the remaining statements in~\sectionref{sec:prelim}. To this end,
we extend the domains of $\overline{T}$ and $\underline{T}$ to functions $f \colon \states{}^{n+1} \to \realsextabove{}$ for any $n \in \nats{}$. For any $n\in\nats$, let  $\closedaboveopenbelow{\mathcal{L}}(\states{}^{n+1})$ be the set of all such functions, and, for any $n\in\nats$, define $\underline{T} \colon \closedaboveopenbelow{\mathcal{L}}(\states{}^{n+1}) \to \closedaboveopenbelow{\mathcal{L}}(\states{}^{n})$ for all $f \in \closedaboveopenbelow{\mathcal{L}}(\states{}^{n+1})$ and all $x_{0:n-1} \in \states{}^{n}$ as
\begin{align}\label{Eq:def_ext_transition}
[\underline{T}\,f](x_{0:n-1}) \coloneqq \Bigl[\underline{T}\bigl(f(x_{0:n-1}\cdot)\bigr)\Bigr](x_{n-1}), 
\end{align}
and analogously for $\overline{T}$.
Note that, in the definition above, we used the function $f(x_{0:n-1}\cdot)\in\gamblesextabove$ that is defined as $f(x_{0:n-1}\cdot)(x_{n}) \coloneqq f(x_{0:n})$ for all $x_{n} \in \states{}$.

Our proofs will at times make use of the fact that the game-theoretic (lower or upper) expectation operators satisfy a \emph{law of iterated expectations}. 
The following makes this explicit.


\begin{proposition}\label{prop:vovk_iterated_expectation}
Consider any game-theoretic lower expectation operator $\underline{\mathbb{E}}^\mathrm{V}$.
Then for all $f\in\overline{\mathcal{L}}(\Omega)$ and all $n\in\natswith$, it holds that
\begin{align*}
 &\quad \overline{\mathbb{E}}^\mathrm{V}\bigl[f\,\vert\,X_{0:n}\bigr] 
 =  \overline{\mathbb{E}}^\mathrm{V}\Bigl[ \overline{\mathbb{E}}^\mathrm{V}\bigl[f\,\vert\,X_{0:n+1}\bigr] \,\Big\vert\, X_{0:n} \Bigr] \text{ and, } \\
  &\quad \underline{\mathbb{E}}^\mathrm{V}\bigl[f\,\vert\,X_{0:n}\bigr] 
 =  \underline{\mathbb{E}}^\mathrm{V}\Bigl[ \underline{\mathbb{E}}^\mathrm{V}\bigl[f\,\vert\,X_{0:n+1}\bigr] \,\Big\vert\, X_{0:n} \Bigr]\,.
\end{align*}
\end{proposition}
\begin{proof}
The equality for the upper expectation follows immediately from \cite[Theorem 16]{natan:game_theory}. 
The equality for the lower expectation then follows from conjugacy.
\end{proof}

We will also need the following technical lemma. 
The statement is obvious when working with measure-theoretic expectations, but perhaps less so for the game-theoretic expectation operators that we use here:
\begin{lemma}\label{lemma:vovk condition fixes value}
Consider any $n,m\in\natswith$ with $m>n$, any $m$-measurable function $f\in\mathcal{L}_m(\Omega)$ and any $x_0,\ldots,x_n\in\states$. 
Then it holds that
\begin{align*}
&
\begin{aligned}
\uexpvovk\bigl[f(X_{0:m})\,&\big\vert\,X_{0:n}=x_{0:n}\bigr] \\
&= \uexpvovk\bigl[f(x_{0:n}X_{n+1:m})\,\big\vert\,X_{0:n}=x_{0:n}\bigr] \text{ and, }
\end{aligned} \\
&
\begin{aligned}
\lexpvovk\bigl[f(X_{0:m})\,&\big\vert\,X_{0:n}=x_{0:n}\bigr] \\
&= \lexpvovk\bigl[f(x_{0:n}X_{n+1:m})\,\big\vert\,X_{0:n}=x_{0:n}\bigr]\,.
\end{aligned}
\end{align*}
\end{lemma}
\begin{proof}
We first give the proof for the upper expectation; we recall the definition of the game-theoretic upper expectation given in~\appendixref{appendix: game-theoretic}, i.e.,
\begin{align*}
&\uexpvovk\bigl[f(X_{0:m})\,\vert\,X_{0:n}=x_{0:n}\bigr] \\
 &= \inf\Bigl\{ \mathcal{M}(x_{0:n})\,\Big\vert\, \mathcal{M}\in\overline{\mathbb{M}}_b\,\text{ and }\,\liminf \mathcal{M} \geq_{x_{0:n}} f(X_{0:m}) \Bigr\}\,,
\end{align*}
where we also recall that $\liminf \mathcal{M} \geq_{x_{0:n}} f(X_{0:m})$ means that for every path $\omega=z_0z_1\cdots z_n\cdots$ in $\Omega$ such that $z_{0:n}=x_{0:n}$, it holds that $\liminf_{k\to+\infty} \mathcal{M}(z_{0:k}) \geq f(\omega)$.

Because this quantifier runs over paths $\omega\in\Omega$ whose first $n+1$ values are $x_0,\ldots,x_n$, it trivially holds that
\begin{align*}
 \liminf \mathcal{M} \geq_{x_{0:n}} f(X_{0:m}) \Leftrightarrow
\liminf \mathcal{M} \geq_{x_{0:n}} f(X_{0:m})\mathbb{I}_{\Gamma(x_{0:n})}\,,
\end{align*}
where we recall the definition $\Gamma(x_{0:n}) \coloneqq \bigl\{ \omega\in\Omega\,\big\vert\, \forall t\in \{0,\ldots,n\}\,:\, \omega(t)=x_t \bigr\}$. 

Now consider this function $f(X_{0:m})\mathbb{I}_{\Gamma(x_{0:n})}$. Then, for any $\omega\in\Omega$, it holds that
\begin{align*}
&f\bigl(X_{0:m}(\omega)\bigr)\mathbb{I}_{\Gamma(x_{0:n})}(\omega) \\
 &= \left\{\begin{array}{ll}
f\bigl(x_{0:n}X_{n+1:m}(\omega)\bigr) & \text{if $X_{0:n}(\omega)=x_{0:n}$, and} \\
0 & \text{otherwise.}
\end{array}\right.
\end{align*}
Therefore in particular, it holds that
\begin{equation*}
f(X_{0:m})\mathbb{I}_{\Gamma(x_{0:n})} = f(x_{0:n}X_{n+1:m})\mathbb{I}_{\Gamma(x_{0:n})}\,. \end{equation*}
We therefore see that
\begin{align*}
\quad&\liminf \mathcal{M} \geq_{x_{0:n}} f(X_{0:m}) \\
\Leftrightarrow
&\liminf \mathcal{M} \geq_{x_{0:n}} f(X_{0:m})\mathbb{I}_{\Gamma(x_{0:n})} \\
\Leftrightarrow
&\liminf \mathcal{M} \geq_{x_{0:n}} f(x_{0:n}X_{n+1:m})\mathbb{I}_{\Gamma(x_{0:n})} \\
\Leftrightarrow
&\liminf \mathcal{M} \geq_{x_{0:n}} f(x_{0:n}X_{n+1:m})\,,
\end{align*}
where the final step again used the fact that the quantification runs over paths whose first $n+1$ values are $x_0,\ldots,x_n$.

Applying the definition of the game-theoretic upper expectation a final time therefore yields
\begin{align*}
 &\uexpvovk\bigl[f(X_{0:m})\,\vert\,X_{0:n}=x_{0:n}\bigr] \\
  &= \uexpvovk\bigl[f(x_{0:n}X_{n+1:m})\,\vert\,X_{0:n}=x_{0:n}\bigr]\,,
\end{align*}
which concludes the proof for the upper expectation.

To prove the claim for the lower expectation, we use the fact that it is defined through conjugacy, together with the above result for the upper expectation, i.e.,
\begin{align*}
 &\hphantom{= -} \lexpvovk\bigl[f(X_{0:m})\,\big\vert\,X_{0:n}=x_{0:n}\bigr] \\
 &= -\uexpvovk\bigl[-f(X_{0:m})\,\big\vert\,X_{0:n}=x_{0:n}\bigr] \\
&= -\uexpvovk\bigl[-f(x_{0:n}X_{n+1:m})\,\big\vert\,X_{0:n}=x_{0:n}\bigr] \\
 &= \hphantom{-} \lexpvovk\bigl[f(x_{0:n}X_{n+1:m})\,\big\vert\,X_{0:n}=x_{0:n}\bigr]\,.
\end{align*}
\end{proof}

\begin{proofof}{\propositionref{prop:vovk_precise_equal_measure_precise_on_n_measurable}}
We will first prove the proposition for the upper expectation $\overline{\mathbb{E}}_P^\mathrm{V}$
.
The equality for the lower expectation will then follow from conjugacy.
Note that we can assume $m\geq n+3$ without loss of generality, because any $\ell$-measurable function is also $m$-measurable for all $m > \ell$.
Since $f_m$ is $m$-measurable, we have that $f_m = f(X_{0:m})$ for some $f \colon \states{}^{m+1} \to \reals{}$.
Now, for any $x_{0:m-1}\in\states^{m}$ we have that 
\begin{align*}
\overline{\mathbb{E}}_P^\mathrm{V}\bigl[f(X_{0:m})\,&\vert\,X_{0:m-1} = x_{0:m-1}\bigr]  \\
&= \overline{\mathbb{E}}_P^\mathrm{V}\bigl[f(x_{0:m-1}X_m)\,\vert\,X_{0:m-1} = x_{0:m-1}\bigr] \\
&= \expprec\bigl[f(x_{0:m-1}X_m)\,\vert\,X_{0:m-1} = x_{0:m-1}\bigr] \\
&= \expprec\bigl[f(X_{0:m})\,\vert\,X_{0:m-1} = x_{0:m-1}\bigr],
\end{align*}
where the first step follows from Lemma~\ref{lemma:vovk condition fixes value} and the second from Equation~\eqref{Equation: precise local models } together with \cite[Proposition 15]{natan:game_theory} which we can use because $f$ is bounded as a consequence of its real-valuedness and the fact that it only takes a finite number of values.
Hence, we have that 
\begin{align*}
\overline{\mathbb{E}}_P^\mathrm{V}\bigl[f(X_{0:m})\,\vert\,X_{0:m-1}\bigr]  
= \expprec\bigl[f(X_{0:m})\,\vert\,X_{0:m-1}\bigr],
\end{align*}
implying that 
\begin{align*}
\overline{\mathbb{E}}_P^\mathrm{V}\bigl[f(X_{0:m}) \,&\vert\,X_{0:m-2}\bigr] \\
&= \overline{\mathbb{E}}_P^\mathrm{V}\Bigl[\overline{\mathbb{E}}_P^\mathrm{V}\bigl[f(X_{0:m})\,\vert\,X_{0:m-1}\bigr]  \,\big\vert\,X_{0:m-2}\Bigr] \\
&= \overline{\mathbb{E}}_P^\mathrm{V}\Bigl[ \expprec\bigl[f(X_{0:m})\,\vert\,X_{0:m-1}\bigr]\,\big\vert\,X_{0:m-2}\Bigr],
\end{align*}
where we have used Proposition~\ref{prop:vovk_iterated_expectation} in the first step.
Now note that $g(X_{0:m-1}) \coloneqq \expprec\bigl[f(X_{0:m})\,\vert\,X_{0:m-1}\bigr]$ is an $(m-1)$-measurable function that is moreover bounded because $f$ is bounded. 
So we can repeat the argument above, where we replace $f$ by $g$ and condition on $X_{0:m-2}$, to find that
\begin{align*}
\overline{\mathbb{E}}_P^\mathrm{V} \bigl[f(X_{0:m})\,&\big\vert\,X_{0:m-2}\bigr] \\
&= \expprec\Bigl[ \expprec\bigl[f(X_{0:m})\,\vert\,X_{0:m-1}\bigr] \,\big\vert\,X_{0:m-2}\Bigr] \\
&= \expprec\bigl[ f(X_{0:m})\,\vert\,X_{0:m-2}\bigr],
\end{align*}
where the last step follows from the law of iterated expectations (applied on $\expprec$).
As such, we can iterate the argument above to find that indeed 

\begin{align*}
\overline{\mathbb{E}}_P^\mathrm{V}\bigl[f_m\,\vert\,X_{0:n}\bigr] = \expprec\bigl[f(X_{0:m})\,\vert\,X_{0:n}\bigr] = \expprec\bigl[f_m\,\vert\,X_{0:n}\bigr].
\end{align*}
The equality for the game-theoretic lower expectation then follows from conjugacy and the fact that $\expprec\bigl[f_m\,\vert\,X_{0:n}\bigr] = - \expprec\bigl[- f_m\,\vert\,X_{0:n}\bigr]$. 
\end{proofof}

\begin{proofof}{\propositionref{prop: vovk local Markov}}
	For any $x\in\states$, let $\mathcal{T}_x$ denote the set of $x$-th rows of $\mathcal{T}$:
	\begin{equation*}
	\mathcal{T}_x \coloneqq \Bigl\{T(x,\cdot)\,:\,T\in\mathcal{T}\Bigr \}\,.
	\end{equation*}
	Then, because each $T\in\mathcal{T}$ is row-stochastic, and because $\mathcal{T}$ is closed and convex, $\mathcal{T}_x$ is a closed and convex set of probability mass functions on $\states$.
	
	Moreover, for any $x\in\states$ and $f\in\gamblesextabove$ it holds that
	\begin{equation*}
	\bigl[\overline{T}f\bigr](x) = \sup_{T\in\mathcal{T}} \bigl[Tf\bigr](x) = \sup_{T(x,\cdot)} \sum_{y\in\states} T(x,y)f(y)\,.
	\end{equation*}
	Therefore, due to~\lemmaref{lemma: local game-theoretic upper exp}, $\bigl[\overline{T}(\cdot)\bigr](x)$ satisfies \ref{coherence: const is const}--\ref{coherence: monotonicity}. Because this is true for every $x\in\states$, it follows that, for all $n\in\natswith$ and all $x_0,\ldots,x_n\in\states$, the operators
	\begin{equation*}
	\overline{Q}_{x_{0:n}}(\cdot) = \bigl[\overline{T}(\cdot)\bigr](x_n)\,,
	\end{equation*}
	also satisfy \ref{coherence: const is const}--\ref{coherence: monotonicity}.
\end{proofof}

\begin{proposition}\label{prop: Vovk iterated T}
Consider any $\mathcal{T}$ and let $\lexpvovk[\cdot\,\vert\,\cdot]$ and $\uexpvovk[\cdot\,\vert\,\cdot]$ be the corresponding game-theoretic lower and upper expectation operators. 
Then for all $n,m\in\natswith$ with $m>n$ and all $f:\states^{m-n+1}\to\reals$, it holds that
\begin{align*}
  \lexpvovk{}[f(X_{n:m})\,\vert\,X_{0:n}] 
  &= \bigr[\underline{T}^{m-n}\, f \bigr] (X_{n})\text{ and, } \\
 \uexpvovk{}[f(X_{n:m})\,\vert\,X_{0:n}] 
  &= \bigr[\overline{T}^{m-n}\, f \bigr] (X_{n}) .
\end{align*}
\end{proposition}
\begin{proof}
We only prove this for the upper expectation; the proof for the lower expectation then follows from conjugacy.
Fix any $n,m \in \natswith{}$ such that $m>n$, and any $f:\states^{m-n+1}\to\reals$.
Note that, for all $x_{0:m-1} \in \states{}^{m}$, we have that
\begin{align}\label{Eq: proof prop: Vovk iterated T}
\uexpvovk\bigl[f(X_{n:m})\,&\vert\,X_{0:m-1}=x_{0:m-1}\bigr] \nonumber \\  
&= \uexpvovk\bigl[f(x_{n:m-1}X_{m})\,\vert\,X_{0:m-1}=x_{0:m-1}\bigr] \nonumber \\
&= \Bigl[\overline{T} \bigl(f(x_{n:m-1} \cdot)\bigr)\Bigr] (x_{m-1}) \nonumber \\
&= [\overline{T}\,f] (x_{n:m-1}),
\end{align}
where the first step follows from Lemma~\ref{lemma:vovk condition fixes value}, the second from Equation~\eqref{eq:def_vovk_lexp} together with \cite[Proposition 15]{natan:game_theory} which we can use because $f$ is bounded as a consequence of its real-valuedness and the fact that it only takes a finite number of values, and the third step from Equation~\eqref{Eq:def_ext_transition}.
If now $m = n+1$, we are done.
If not, we have that $m > n+1$. Then, using Proposition \ref{prop:vovk_iterated_expectation}, we find that
\begin{align}\label{Eq: proof prop: Vovk iterated T next step}
\uexpvovk\bigl[f(X_{n:m})\,&\vert\,X_{0:m-2}\bigr] \nonumber \\
 &= \uexpvovk\Bigl[\uexpvovk\bigl[f(X_{n:m})\,\vert\,X_{0:m-1}\bigr] \Big\vert X_{0:m-2}\Bigr]\,,
\end{align}
and, therefore, for all $x_{0:m-2} \in \states{}^{m-1}$ we have that 
\begin{align*}
 \uexpvovk\bigl[&f(X_{n:m})\,\vert\,X_{0:m-2}=x_{0:m-2}\bigr]  \\
 &= \uexpvovk\Bigl[\uexpvovk\bigl[f(X_{n:m})\,\vert\,X_{0:m-2}=x_{0:m-2},X_{m-1}\bigr] \\
  &\quad\quad\quad\quad\quad\quad\quad\quad\quad\quad\quad\quad\quad\Big\vert X_{0:m-2}=x_{0:m-2}\Bigr] \\
&= \uexpvovk\Bigl[[\overline{T}\,f](x_{n:m-2}X_{m-1}) \Big\vert X_{0:m-2}=x_{0:m-2}\Bigr] \\
&= \bigl[\overline{T}^2f\bigr] (x_{n:m-2}),
\end{align*}
where the first step follows from Equation~\eqref{Eq: proof prop: Vovk iterated T next step} and Lemma~\ref{lemma:vovk condition fixes value}, the second step from Equation \eqref{Eq: proof prop: Vovk iterated T} and the third from Equation~\eqref{eq:def_vovk_lexp} together with \cite[Proposition 15]{natan:game_theory}, which we can use once more because $\overline{T}\,f$ is bounded as a consequence of the boundedness of $f$ and~\lemmaref{lemma:lower_trans_basic_properties}, and the third step from Equation~\eqref{Eq:def_ext_transition}.
If $m = n+2$, we are done.
If not, then we can continue in this way until we find that $\uexpvovk{}[f(X_{n:m})\,\vert\,X_{0:n}=x_{0:n}] 
  = \bigr[\overline{T}^{m-n}\, f \bigr] (x_{n})$.
\end{proof}

\begin{corollary}\label{prop: vovk time shift}
Consider any $n \in \nats{}$ and any function $f:\states^{n}\to\reals$. 
Then, for any $x,y \in \states{}$, it holds that
\begin{align}\label{eq: prop: vovk time shift}
\lexpvovk\bigl[f(X_{1:n})\,\big\vert\,X_{0:1} = xy \bigr] = \lexpvovk\bigl[f(X_{0:n-1})\,\big\vert\,X_0 = y \bigr]\, \\
\uexpvovk\bigl[f(X_{1:n})\,\big\vert\,X_{0:1} = xy \bigr] = \uexpvovk\bigl[f(X_{0:n-1})\,\big\vert\,X_0 = y \bigr]\,. \nonumber
\end{align}
\end{corollary}
\begin{proof}
It follows from Proposition \ref{prop: Vovk iterated T} that, for any $x,y \in \states{}$, both sides of~\eqref{eq: prop: vovk time shift} are equal to $\bigr[\underline{T}^{n-1}\, f \bigr] (y)$. 
The same reasoning holds for the upper expectations.

\end{proof}

\begin{proofof}{\propositionref{prop:vovk_imprecise_dominates_compatible_precise_vovk}}
We first prove the inequality for the upper expectation.
Note that 
\begin{align*}
\uexpvovk[g(X_{n+1})\,\vert\,X_{0:n}=x_{0:n}] 
&= \bigl[\overline{T}g\bigr](x_n) \\
&\geq \mathbb{E}_P[g(X_{n+1})\,\vert\,X_{0:n}=x_{0:n}] \\
&= \overline{\mathbb{E}}_P^\mathrm{V}[g(X_{n+1})\,\vert\,X_{0:n}=x_{0:n}],
\end{align*}
for all $g\in\gamblesextabove$, all $n\in\natswith$, all $x_0,\ldots,x_n\in\states$ and all $P\in\imcirr$, where the last step used Equation~\eqref{Equation: precise local models } together with \cite[Proposition 15]{natan:game_theory} which is applicable here because $g$ is bounded below.
Then it follows immediately from \cite[Proposition 14]{natan:game_theory} that 
\begin{equation*}
\uexpvovk[f\,\vert\,X_{0:n}=x_{0:n}] \geq \overline{\mathbb{E}}_P^\mathrm{V}[f\,\vert\,X_{0:n}=x_{0:n}]\,,
\vspace{2pt}
\end{equation*}
for all $f\in\overline{\mathcal{L}}(\Omega)$, all $n\in\natswith$, all $x_0,\ldots,x_n\in\states$ and all $P\in\imcirr$, implying that 
\begin{equation*}
\uexpvovk[f\,\vert\,X_{0:n}=x_{0:n}] \geq \sup_{P\in\imcirr} \overline{\mathbb{E}}_P^\mathrm{V}[f\,\vert\,X_{0:n}=x_{0:n}].
\vspace{2pt}
\end{equation*}
The inequality for lower expectations then follows from conjugacy.
\end{proofof}

\begin{proofof}{\propositionref{prop:vovk_monotone_continuity}}
From the fact that each $f_m$ is $m$-measurable and $\states$ is finite we infer that $f_m$ is bounded below and above. Moreover, since the sequence is non-decreasing, $\{f_m\}_{m\in\nats}$ and $f$ are uniformly bounded below. 
Hence the claim for the upper expectation is proved by~\cite[Theorem 24]{natan:game_theory}.

Since $\{f_m\}_{m\in\nats}$ is non-decreasing, the corresponding sequence $\{-f_m\}_{m\in\nats}$ is non-increasing, and clearly each function $-f_m$ is still $m$-measurable. By~\cite[Proposition 30]{natan:game_theory} we therefore have
\begin{equation*}
\overline{\mathbb{E}}^\mathrm{V}[-f\,\vert\,X_{0:n}=x_{0:n}] = \lim_{m\to+\infty} \overline{\mathbb{E}}^\mathrm{V}[-f_m\,\vert\,X_{0:n}=x_{0:n}]\,.
\end{equation*}
Hence, the claim for the lower expectation follows from conjugacy of the lower and upper expectations.
\end{proofof}

\section{Proofs of Statements in~\sectionref{subsec:expected_hit_times}}

\begin{lemma}\label{lemma:hitting_time_approx_sequence_nonnegative}
Consider the sequences $\smash{\underline{h}_A^{(n)}}$ and $\smash{\overline{h}_A^{(n)}}$ defined in~\sectionref{subsec:expected_hit_times}. Then $\smash{0 \leq \underline{h}_A^{(n)} \leq \overline{h}_A^{(n)} \leq n+1}$ for all $n\in\natswith$.
\end{lemma}
\begin{proof}
We give a proof by induction. First we note that by definition $\smash{\underline{h}_A^{(0)}=\overline{h}_A^{(0)}=\mathbb{I}_{A^c}}$. Hence, because $0 \leq \mathbb{I}_{A^c} \leq 1$, we have established the required induction base.

Now assume that the claim is true for $n-1$, with $n\in\nats$; we will show that the claim is also true for $n$. We start by establishing the non-negativity of $\smash{\underline{h}_A^{(n)}}$, which by definition is given by
\begin{equation*}
\underline{h}_A^{(n)} = \mathbb{I}_{A^c} + \mathbb{I}_{A^c}\cdot \underline{T}\,\underline{h}_A^{(n-1)}\,.
\end{equation*}
By the induction hypothesis, $\underline{h}_A^{(n-1)}$ is non-negative, which implies that also $\smash{\underline{T}\,\underline{h}_A^{(n-1)}}$ is non-negative. Hence, because $\mathbb{I}_{A^c}$ is non-negative, also the product $\smash{\mathbb{I}_{A^c}\cdot \underline{T}\,\underline{h}_A^{(n-1)}}$ is non-negative. Noting again that $\mathbb{I}_{A^c}$ is non-negative, it follows that therefore the sum $\smash{\mathbb{I}_{A^c} + \mathbb{I}_{A^c}\cdot \underline{T}\,\underline{h}_A^{(n-1)}}$ is non-negative, which concludes the proof that $0\leq \underline{h}_A^{(n)}$.

To show that $\smash{\underline{h}_A^{(n)}\leq \overline{h}_A^{(n)}}$, we first use the induction hypothesis and the monotonicity of $\underline{T}$ to establish that
\begin{align*}
\underline{h}_A^{(n)} &= \mathbb{I}_{A^c} + \mathbb{I}_{A^c}\cdot \underline{T}\,\underline{h}_A^{(n-1)} \leq \mathbb{I}_{A^c} + \mathbb{I}_{A^c}\cdot \underline{T}\overline{h}_A^{(n-1)}\,.
\end{align*}
Moreover, it follows from the definitions of $\underline{T}$ and $\overline{T}$ that $\underline{T}f \leq \overline{T}f$ for any $f\in\gamblesextabove$, and hence we find
\begin{align*}
\underline{h}_A^{(n)} &\leq \mathbb{I}_{A^c} + \mathbb{I}_{A^c}\cdot \underline{T}\overline{h}_A^{(n-1)} \leq \mathbb{I}_{A^c} + \mathbb{I}_{A^c}\cdot \overline{T}\overline{h}_A^{(n-1)} = \overline{h}_A^{(n)}\,,
\end{align*}
which concludes the proof that $\smash{\underline{h}_A^{(n)} \leq \overline{h}_A^{(n)}}$.

It remains to show that $\smash{\overline{h}_A^{(n)} \leq n+1}$. By definition we have that
\begin{equation*}
\overline{h}_A^{(n)} = \mathbb{I}_{A^c} + \mathbb{I}_{A^c}\cdot\overline{T}\,\overline{h}_A^{(n-1)}\,.
\end{equation*}
The induction hypothesis tells us that $\smash{\overline{h}_A^{(n-1)} \leq n}$, which implies that also $\smash{\overline{T}\,\overline{h}_A^{(n-1)} \leq n}$. Because $0\leq \mathbb{I}_{A^c} \leq 1$ it holds that the product $\smash{\mathbb{I}_{A^c}\cdot\overline{T}\,\overline{h}_A^{(n-1)} \leq n}$. Now, again because $\mathbb{I}_{A^c} \leq 1$, it follows that the sum $\smash{\mathbb{I}_{A^c} + \mathbb{I}_{A^c}\cdot\overline{T}\,\overline{h}_A^{(n-1)} \leq n+1}$, which concludes the proof.
\end{proof}

\begin{proofof}{\lemmaref{lemma:iterate_hit_time_is_restricted_exp}}
We first give the proof for $\overline{h}_A^{(n)}$; the proof for $\underline{h}_A^{(n)}$ will follow an analogous reasoning. 

This proof is by induction, and we first establish the induction base. To this end, note that $\smash{H_A^{(0)}(X_0)}$ only depends on the state $X_0$, so it follows from~[4.] in~\propositionref{prop:vovk_satisfies_coherence} that, for all $x_0\in\states$, it holds that
\begin{equation*}
\uexpvovk\bigl[H_A^{(0)}(X_0)\,\big\vert\,X_0=x_0\bigr] = H_A^{(0)}(x_0)\,.
\end{equation*}
Furthermore, it follows immediately from the definitions of $H_A^{(0)}$ and $\overline{h}_A^{(0)}$ that
\begin{equation*}
H_A^{(0)}(x_0) = \mathbb{I}_{A^c}(x_0) = \overline{h}_A^{(0)}(x_0)\,,
\end{equation*}
for all $x_0\in\states$, and so we get
\begin{equation*}
\uexpvovk\bigl[H_A^{(0)}(X_0)\,\big\vert\,X_0=x_0\bigr] = \overline{h}_A^{(0)}(x_0)\,,
\end{equation*}
which is the induction base that we are after.
Note that in the above, we used a small notational trick to write $H_A^{(0)}(x_0)$ for the value of $H_A^{(0)}(X_0)$ on a path $\omega$ that satisfies $\omega(0)=x_0$.

In fact, this notational trick can be extended by noting that, for any $n\in\natswith$, the value of $H_A^{(n)}$ on any $\omega\in\Omega$ is completely determined by the values $\omega(0),\ldots,\omega(n)\in\states$. Thus, we can equivalently interpret $H_A^{(n)}$ as a function on $\states^{n+1}$. Moreover, we can then write, for any $n\in\nats$ and any $x_0,\ldots,x_n\in\states$, that
\begin{equation*}
H_A^{(n)}(x_{0:n}) = \mathbb{I}_{A^c}(x_0) + \mathbb{I}_{A^c}(x_0)H_A^{(n-1)}(x_{1:n})\,,
\end{equation*}
with $H_A^{(0)}(x_0)=\mathbb{I}_{A^c}(x_0)$. Using this observation, let us now proceed with the induction step; so, we assume that the statement is true for some $n-1$, with $n\in\nats$. We will show that the statement is also true for $n$.

First, for any $x_0\in\states$ we have that
\begin{align}\label{Eq: proof of lemma:iterate_hit_time_is_restricted_exp}
&\quad \uexpvovk\bigl[H_A^{(n)}(X_{0:n})\,\vert\,X_0=x_0\bigr]  \nonumber \\
&= \uexpvovk\bigl[\mathbb{I}_{A^c}(x_0) + \mathbb{I}_{A^c}(x_0)H_A^{(n-1)}(X_{1:n})\,\vert\,X_0=x_0\bigr] \nonumber \\
&= \mathbb{I}_{A^c}(x_0) + \mathbb{I}_{A^c}(x_0)\uexpvovk\bigl[H_A^{(n-1)}(X_{1:n})\,\vert\,X_0=x_0\bigr] \nonumber \\
&= \mathbb{I}_{A^c}(x_0) + \mathbb{I}_{A^c}(x_0)\uexpvovk\Bigl[\uexpvovk\bigl[H_A^{(n-1)}(X_{1:n})\,\vert\,X_{0:1}\bigr]\,\Big\vert\,X_0=x_0\Bigr],
\end{align}
where the first step used~\lemmaref{lemma:vovk condition fixes value}, the second step used [2.] and [5.] in Proposition~\ref{prop:vovk_satisfies_coherence}, and the last step used Proposition~\ref{prop:vovk_iterated_expectation}.
Now, according to Corollary~\ref{prop: vovk time shift}, $\uexpvovk\bigl[H_A^{(n-1)}(X_{1:n})\,\vert\,X_{0:1}\bigr]$ does not depend on the initial state $X_0$, so there is a function $g_{n-1} \colon \states{} \to \realsext{}$ such that
\begin{align*}
g_{n-1}(X_1) = \uexpvovk\bigl[H_A^{(n-1)}(X_{1:n})\,\vert\,X_{0:1}\bigr].
\end{align*}
Moreover, again using Corollary~\ref{prop: vovk time shift}, we also have that 
\begin{align*}
g_{n-1}(X_0) = \uexpvovk\bigl[H_A^{(n-1)}(X_{0:n-1})\,\vert\,X_{0}\bigr].
\end{align*}
By the induction hypothesis, we therefore have that $g_{n-1}=\overline{h}_A^{(n-1)}$. Plugging this back into Equation \eqref{Eq: proof of lemma:iterate_hit_time_is_restricted_exp}, we get that
\begin{align}\label{Equation: Eq: proof of lemma:iterate_hit_time_is_restricted_exp_2}
&\quad \uexpvovk\bigl[H_A^{(n)}(X_{0:n})\,\vert\,X_0=x_0\bigr] \nonumber \\
&= \mathbb{I}_{A^c}(x_0) + \mathbb{I}_{A^c}(x_0)\uexpvovk\Bigl[\overline{h}_A^{(n-1)}(X_1)\,\Big\vert\,X_0=x_0\Bigr] \nonumber \\
&= \mathbb{I}_{A^c}(x_0) + \mathbb{I}_{A^c}(x_0)[\overline{T} \, \overline{h}_A^{(n-1)}] (x_0) 
= \overline{h}_A^{(n)}(x_0), 
\end{align}
for all $x_0\in\states$. 
In the expression above, the second step uses Equation \eqref{eq:def_vovk_lexp} together with \cite[Proposition 15]{natan:game_theory} which we can use because $\overline{h}_A^{(n-1)}$ is non-negative (and therefore bounded below) by~\lemmaref{lemma:hitting_time_approx_sequence_nonnegative}, and the last step uses the definition of $\overline{h}_A^{(n)}$. This concludes the first part of the proof.

Analogously, we can follow the reasoning above in order to show that the statement for $\underline{h}_A^{(n)}$ holds.
Properties [2.] and [5.] in Proposition~\ref{prop:vovk_satisfies_coherence} clearly also hold for lower expectations because of conjugacy.
The only step that requires some closer attention is the second equality in Equation~\eqref{Equation: Eq: proof of lemma:iterate_hit_time_is_restricted_exp_2} because Equation~\eqref{eq:def_vovk_lexp} and \cite[Proposition 15]{natan:game_theory} are only given for upper expectations.
However, we can easily derive this through the use of conjugacy:
\begin{align*}
\lexpvovk\Bigl[\underline{h}_A^{(n-1)}(X_1)\,\Big\vert\,X_0=x_0\Bigr]
&= - \uexpvovk\Bigl[- \underline{h}_A^{(n-1)}(X_1)\,\Big\vert\,X_0=x_0\Bigr] \\
&= - [\overline{T} \, \bigl(-\underline{h}_A^{(n-1)}\bigr)] (x_0) \\
&= [\underline{T} \, \underline{h}_A^{(n-1)}] (x_0),
\end{align*}
for all $x_0 \in \states{}$.
Again, we used Equation \eqref{eq:def_vovk_lexp} together with \cite[Proposition 15]{natan:game_theory} in the second step, which was allowed since $\underline{h}_A^{(n-1)}$ is bounded above (and therefore $-\underline{h}_A^{(n-1)}$ is bounded below) as a consequence of~\lemmaref{lemma:hitting_time_approx_sequence_nonnegative}.
\end{proofof}

The remainder of this section of the appendix sets up the required results to prove the statements for upper expectations in~\theoremref{thm:lower_hitting_time_reach_and_equal} and~\corollaryref{cor:imprecise_hitting_time_is_minimal_system_solution}. 

The following generic property will be useful.
\begin{lemma}\label{lemma:events_conjugate_complement}
For all $B\subseteq \states$, with $B^c\coloneqq \states\setminus B$, it holds for all $n\in\nats$ and all $x\in\states$ that
\begin{equation*}
\bigl[\underline{T}^n\mathbb{I}_B\bigr](x) = 0\quad\Leftrightarrow\quad \bigl[\overline{T}^n\mathbb{I}_{B^c}\bigr](x) = 1\,,
\end{equation*}
and, moreover,
\begin{equation*}
\bigl[\underline{T}^n\mathbb{I}_B\bigr](x) = 1\quad\Leftrightarrow\quad \bigl[\overline{T}^n\mathbb{I}_{B^c}\bigr](x) = 0\,.
\end{equation*}
\end{lemma}
\begin{proof}
Fix $x\in\states$ and suppose that $\bigl[\underline{T}^n\mathbb{I}_B\bigr](x)=0$. Note that $\mathbb{I}_{B^c} = 1-\mathbb{I}_{B}$. Hence, using~\lemmaref{lemma:lower_trans_composite_basic_properties} and~\lemmaref{lemma:transop_composite_conjugate_gambles},
\begin{equation*}
\bigl[\overline{T}^n\mathbb{I}_{B^c}\bigr](x) = \bigl[\overline{T}^n(1-\mathbb{I}_B)\bigr](x) = 1 -\bigl[\underline{T}^n\mathbb{I}_B\bigr](x) = 1\,.
\end{equation*}
Conversely, suppose that $\bigl[\overline{T}^n\mathbb{I}_{B^c}\bigr](x) = 1$ and note that $\mathbb{I}_{B}=1-\mathbb{I}_{B^c}$. Then using~\lemmaref{lemma:lower_trans_composite_basic_properties} and~\lemmaref{lemma:transop_composite_conjugate_gambles},
\begin{equation*}
\bigl[\underline{T}^n\mathbb{I}_{B}\bigr](x) = \bigl[\underline{T}^n(1-\mathbb{I}_{B^c})\bigr](x) = 1 -\bigl[\overline{T}^n\mathbb{I}_{B^c}\bigr](x) = 0\,.
\end{equation*}
The proof of the other claim is completely analogous.
\end{proof}

The proofs of~\theoremref{thm:lower_hitting_time_reach_and_equal} and~\corollaryref{cor:imprecise_hitting_time_is_minimal_system_solution} work by first proving the statements of interest on a modification of the actual (imprecise) stochastic process, and then showing that the inferences are invariant under this modification. We will call this modification the \emph{$A$-inert modification} of the imprecise Markov chain, which is based on the following definition:
\begin{definition}\label{def:a_inert}
A non-empty set $\transmatset$ of transition matrices is called \emph{$A$-inert} if, for all $x\in A$ and all $T\in\transmatset$, it holds that $T(x,x)=1$. If $\transmatset$ is $A$-inert, we also say that its corresponding lower and upper transition operators are $A$-inert, or, similarly, that an imprecise Markov chain derived from it is $A$-inert.
\end{definition}

\begin{lemma}\label{lemma:a_inert_fixed_value}
Let $\transmatset$ be a non-empty set of transition matrices that is $A$-inert, and let $\underline{T}$ and $\overline{T}$ be its lower and upper transition operators. Then, for all $f\in\gambles$, it holds that $\bigl[\underline{T}\,f\bigr](x)=f(x)$ and $\bigl[\overline{T}\,f\bigr](x)=f(x)$ for all $x\in A$.
\end{lemma}
\begin{proof}
Fix $f\in\gambles$ and $x\in A$, and choose any $T\in\transmatset$. Then because $\transmatset$ is $A$-inert,
\begin{equation*}
\bigl[Tf\bigr](x) = \sum_{y\in\states} T(x,y)f(y) = f(x)\,.
\end{equation*}
Then from the definition of $\underline{T}$ we get
\begin{equation*}
\bigl[\underline{T}\,f\bigr](x) = \inf_{T\in\transmatset} \bigl[Tf\bigr](x) = f(x)\,,
\end{equation*}
and from the definition of $\overline{T}$ we get
\begin{equation*}
\bigl[\overline{T}\,f\bigr](x) = \sup_{T\in\transmatset} \bigl[Tf\bigr](x) = f(x)\,,
\end{equation*}
which concludes the proof.
\end{proof}

\begin{corollary}\label{cor:a_inert_on_indicators}
Suppose that $\underline{T}$ is an $A$-inert lower transition operator. Then $\bigl[\underline{T}\mathbb{I}_{\{x\}}\bigr](x)=1$ for all $x\in A$.
\end{corollary}
\begin{proof}
Apply~\lemmaref{lemma:a_inert_fixed_value} to the function $\mathbb{I}_{\{x\}}\in\gambles$.
\end{proof}
Intuitively, it can be seen that an imprecise Markov chain is $A$-inert, if all states $x\in A$ are \emph{absorbing}: once the process reaches a state $x\in A$, it remains there with (lower) probability one. The modification is now straightforward; we start with the modification of a set of transition matrices:
\begin{definition}\label{def:modified_set_a_inert}
For any non-empty set of transition matrices $\transmatset$, its \emph{$A$-inert modification} $\mathcal{S}$ is the set of transition matrices defined by
\begin{align*}
\mathcal{S} \coloneqq \Bigl\{ S\in\mathbb{T} \,\Big\vert\, &\bigl(\forall x\in A: S(x,x)=1\bigr), \\
 &\bigl(\exists T\in\mathcal{T}: \forall x\in A^c : S(x,\cdot) = T(x,\cdot)\bigr)\Bigr\}\,,
\end{align*}
where $\mathbb{T}$ is the set of all transition matrices.
\end{definition}

This $A$-inert modification satisfies the following properties:
\begin{lemma}\label{lemma:ainert_mod_set_satisfies_properties}
Let $\transmatset$ be a non-empty set of transition matrices that is closed and convex and that has separately specified rows, and let $\mathcal{S}$ be its $A$-inert modification. Then $\mathcal{S}$ is non-empty, closed, convex, $A$-inert, and has separately specified rows.

Moreover, let $\underline{T}$ and $\overline{T}$, and $\underline{S}$ and $\overline{S}$, denote the lower and upper transition operators corresponding to $\transmatset$ and $\mathcal{S}$, respectively. Then, for all $f\in\gamblesextabove$, if holds that $\bigl[\underline{T}\,f\bigr](x) = \bigl[\underline{S}\,f\bigr](x)$ and $\bigl[\overline{T}\,f\bigr](x) = \bigl[\overline{S}\,f\bigr](x)$ for all $x\in A^c$.
\end{lemma}
\begin{proof}
We first show that $\mathcal{S}$ is non-empty. To this end, choose any $T\in\transmatset$; this is possible because $\transmatset$ is non-empty. Consider the $\lvert\states\rvert\times \lvert\states\rvert$ matrix $S$, whose elements are defined, for all $x,y\in\states$, as
\begin{equation*}
S(x,y) \coloneqq \left\{\begin{array}{ll}
1 & \text{if $x=y$ and $x\in A$,} \\
0 & \text{if $x\neq y$ and $x\in A$, and} \\
T(x,y) & \text{otherwise.}
\end{array}\right.
\end{equation*}
Then for all $x\in A$ we have that $\sum_{y\in\states}S(x,y)=1$ and $S(x,y)\geq 0$ for all $y\in\states$. Moreover, for all $x\in A^c$ we have $\sum_{y\in\states}S(x,y)=\sum_{y\in\states}T(x,y)=1$ and $S(x,y)=T(x,y)\geq 0$ for all $y\in\states$, because $T$ is a transition matrix. Hence, $S$ is also a transition matrix, so an element of $\mathbb{T}$, which means that $S\in \mathcal{S}$.

We next show that $\mathcal{S}$ is closed. To this end, consider any converging sequence $\{S_n\}_{n\in\natswith}$ in $\mathcal{S}$ with $S_* \coloneqq \lim_{n\to+\infty}S_n$. First note that, for all $n\in\natswith$, $S_n(x,x)=1$ for all $x\in A$, because $S_n\in\mathcal{S}$. This implies that also $S_*(x,x)=1$ for all $x\in A$. Similarly, $S_*(x,y)=0$ for all $x\in A$ and $y\neq x$.

Next, fix an arbitrary $T\in\transmatset$. Then, for all $n\in\natswith$, because $S_n\in\mathcal{S}$, there is some $T_n\in\transmatset$ such that $S_n(x,\cdot)=T_n(x,\cdot)$ for all $x\in A^c$. Moreover, because $\transmatset$ has separately specified rows, there is some $T_n'\in\transmatset$ such that $T_n'(x,\cdot)=T(x,\cdot)$ for all $x\in A$, and $T_n'(x,\cdot)=T_n(x,\cdot)$ for all $x\in A^c$. Taking limits in $n$, for $x\in A$ we immediately get $\lim_{n\to+\infty} T_n'(x,\cdot)=\lim_{n\to+\infty}T(x,\cdot)=T(x,\cdot)$. Moreover, for $x\in A^c$ we get
\begin{equation*}
\lim_{n\to+\infty} T_n'(x,\cdot)=\lim_{n\to+\infty}T_n(x,\cdot)=\lim_{n\to+\infty}S_n(x,\cdot)=S_*(x,\cdot)\,.
\end{equation*}
Hence, there exists the limit $T_*' \coloneqq \lim_{n\to+\infty} T_n'$ and, because $\transmatset$ is closed and $T_n'\in\transmatset$ for all $n\in\natswith$, we have $T_*'\in\transmatset$. Now note that $S_*(x,\cdot)=T_*'(x,\cdot)$ for all $x\in A^c$, whence $S_*\in\mathcal{S}$.

To prove that $\mathcal{S}$ is convex, fix $S_1,S_2\in\mathcal{S}$ and $\lambda\in [0,1]$, and define $S\coloneqq \lambda S_1 + (1-\lambda) S_2$; we need to show that $S\in\mathcal{S}$. Because for $x\in A$ we have $S_1(x,x)=S_2(x,x)=1$, clearly also $S(x,x)=1$. Similarly, $S(x,y)=0$ for all $x\in A$ and $y\neq x$. Next, because $S_1,S_2\in\mathcal{S}$, there exist $T_1,T_2\in\transmatset$ such that $S_1(x,\cdot)=T_1(x,\cdot)$ and $S_2(x,\cdot)=T_2(x,\cdot)$ for all $x\in A^c$. Let $T\coloneqq \lambda T_1 + (1-\lambda) T_2$. Because $\transmatset$ is convex, it holds that $T\in\transmatset$. Moreover, for all $x\in A^c$,
\begin{align*}
S(x,\cdot) &= \lambda S_1(x,\cdot) + (1-\lambda)S_2(x,\cdot) \\
 &= \lambda T_1(x,\cdot) + (1-\lambda)T_2(x,\cdot) = T(x,\cdot)\,.
\end{align*}
Hence we have found a $T\in\transmatset$ such that $S(x,\cdot)=T(x,\cdot)$ for all $x\in A^c$, which implies $S\in\mathcal{S}$.

The fact that $\mathcal{S}$ is $A$-inert is immediate from the definition; $S(x,x)=1$ for all $S\in\mathcal{S}$.

We now prove that $\mathcal{S}$ has separately specified rows. So, for all $x\in\states$, choose $S_x\in\mathcal{S}$, and define the matrix $S$ as $S(x,\cdot)\coloneqq S_x(x,\cdot)$ for all $x\in\states$. We need to show that $S\in\mathcal{S}$. First, for any $x\in A$, because $S_x\in\mathcal{S}$, it holds that $S(x,x)=S_x(x,x)=1$. Similarly, $S_x(x,y)=0$ for all $x\in A$ and $y\neq x$. Next, for any $x\in A^c$, because $S_x\in\mathcal{S}$, there is some $T_x\in\transmatset$ such that $S_x(x,\cdot)=T_x(x,\cdot)$. Because $\transmatset$ has separately specified rows, there exists some $T\in\transmatset$ such that $T(x,\cdot)=T_x(x,\cdot)$ for all $x\in A^c$. Hence also $T(x,\cdot)=S_x(x,\cdot)=S(x,\cdot)$ for all $x\in A^c$, which implies $S\in\mathcal{S}$.

For the property about the lower transition operators, fix any $f\in\gamblesextabove$ and $x\in A^c$. Using~\lemmaref{lemma:trans_op_continuous}, we can find $T\in\transmatset$ and $S\in\mathcal{S}$ such that
\begin{equation*}
\bigl[\underline{T}\,f\bigr](x) = \bigl[T\,f\bigr](x)\quad\text{and}\quad \bigl[\underline{S}\,f\bigr](x) = \bigl[S\,f\bigr](x)\,.
\end{equation*}
Because $T\in\transmatset$, there is some corresponding $S_T\in\mathcal{S}$ such that $S_T(x,\cdot)=T(x,\cdot)$. This implies that
\begin{equation*}
\bigl[\underline{S}\,f\bigr](x) \leq \bigl[S_Tf\bigr](x) = \bigl[Tf\bigr](x) = \bigl[\underline{T}\,f\bigr](x)\,.
\end{equation*}
Similarly, because $S\in\mathcal{S}$, there is some corresponding $T_S\in\transmatset$ such that $T_S(x,\cdot)=S(x,\cdot)$. Hence 
\begin{equation*}
\bigl[\underline{T}\,f\bigr](x) \leq \bigl[T_Sf\bigr](x) = \bigl[Sf\bigr](x) = \bigl[\underline{S}\,f\bigr](x)\,,
\end{equation*}
from which we get $\bigl[\underline{T}\,f\bigr](x) = \bigl[\underline{S}\,f\bigr](x)$. 

The argument for the upper transition operators proceeds analogously, by simply changing the directions of the inequalities.
\end{proof}

We can now finally define the $A$-inert modification of an imprecise Markov chain:
\begin{definition}\label{def:modified_imprecise_markov}
For any imprecise Markov chain that is characterised by a set $\transmatset$ of transition matrices, its (unique) \emph{$A$-inert modification} is the imprecise Markov chain that is characterised by the $A$-inert modification $\mathcal{S}$ of $\transmatset$, and that has the same initial model as the original imprecise Markov chain.
\end{definition}



It is clear that the $A$-inert modification of a given imprecise Markov chain, is indeed $A$-inert. In what follows, we will want to work with the $A$-inert modification of the actual imprecise Markov chain under consideration, when proving results about the expected hitting time. The reason that we can do this, intuitively, is that this modification only changes the behaviour of the process \emph{once it has visited $A$}, while the expected hitting time is only dependent on the behaviour \emph{before} $A$ is reached. The following makes this explicit for precise, homogeneous Markov chains:
\begin{lemma}\label{lemma:hit_time_homogen_is_a_inert}
Consider a homogeneous Markov chain $P$ with transition matrix $T$. Because $T$ can be interpreted as a singleton set $\{T\}$, the $A$-inert modification $\{S\}$ is well defined. Let $Q$ be the unique homogeneous Markov chain with transition matrix $S$ and initial model $Q(X_0)=P(X_0)$. Then $h_A^P = h_A^Q$.
\end{lemma}
\begin{proof}
By~\lemmaref{lemma:homogen_hitting_time_is_minimal_system_solution}, it holds that
\begin{equation*}
h_A^P = \mathbb{I}_{A^c} + \mathbb{I}_{A^c}\cdot Th_A^P\,.
\end{equation*}
Thus for any $x\in A$, because $\mathbb{I}_{A^c}(x)=0$, it holds that
\begin{align*}
h_A^P(x) &= \mathbb{I}_{A^c}(x) + \mathbb{I}_{A^c}(x)\cdot \bigl[Th_A^P\bigr](x) \\
 &= \mathbb{I}_{A^c}(x) + \mathbb{I}_{A^c}(x)\cdot \bigl[Sh_A^Q\bigr](x) = 0\,.
\end{align*}
Next we point out that the lower (and upper) transition operator corresponding to the singleton set $\{T\}$ is simply the matrix $T$ itself; analogously, the lower (and upper) transition operator of $\{S\}$ is again the matrix $S$ itself. Hence, by~\lemmaref{lemma:ainert_mod_set_satisfies_properties}, for all $f\in\gamblesextabove$ and $x\in A^c$ it holds that $\bigl[Sf\bigr](x) = \bigl[Tf\bigr](x)$, so then also
\begin{align*}
h_A^P(x) &= \mathbb{I}_{A^c}(x) + \mathbb{I}_{A^c}(x)\cdot \bigl[Th_A^P\bigr](x) \\
 &= \mathbb{I}_{A^c}(x) + \mathbb{I}_{A^c}(x)\cdot \bigl[Sh_A^P\bigr](x)\,.
\end{align*}
We can therefore conclude that
\begin{equation*}
h_A^P = \mathbb{I}_{A^c} + \mathbb{I}_{A^c}\cdot Sh_A^P\,,
\end{equation*}
which by~\lemmaref{lemma:homogen_hitting_time_is_minimal_system_solution} implies that $h_A^Q \leq h_A^P$. But by symmetry we can modify the above argument to also find $h_A^P\leq h_A^Q$.
\end{proof}

Similarly, as the next result shows, the lower and upper expected hitting times of a game-theoretic imprecise Markov chain, are equal to those of its $A$-inert modification:
\begin{lemma}\label{lemma:hit_time_game_is_a_inert}
Consider a game-theoretic imprecise Markov chain parameterised by $\mathcal{T}$. Then its lower and upper expected hitting time are equal to those of its $A$-inert modification.
\end{lemma}
\begin{proof}
By~\propositionref{prop:limit_is_lower_exp}, it holds that $\lexpvovk[H_A\,\vert\,X_0] = \lim_{n\to+\infty} \underline{h}_A^{(n)}$ and $\uexpvovk[H_A\,\vert\,X_0] = \lim_{n\to+\infty} \overline{h}_A^{(n)}$.

Now consider the $A$-inert modification $\mathcal{S}$ of $\mathcal{T}$, and let $\underline{S}$ and $\overline{S}$ be the corresponding lower and upper transition operators, respectively.

Then consider the corresponding sequences $\underline{s}_A^{(n)}$ and $\overline{s}_A^{(n)}$, so, $\underline{s}_A^{(0)}\coloneqq \overline{s}_A^{(0)} \coloneqq \mathbb{I}_{A^c}$ and, for all $n\in\natswith$,
\begin{equation*}
\underline{s}_A^{(n+1)} \coloneqq \mathbb{I}_{A^c} + \mathbb{I}_{A^c}\cdot \underline{S}\,\underline{s}_A^{(n)}\,,
\end{equation*}
and
\begin{equation*}
\overline{s}_A^{(n+1)} \coloneqq \mathbb{I}_{A^c} + \mathbb{I}_{A^c}\cdot \overline{S}\,\overline{s}_A^{(n)}\,.
\end{equation*}
Due to~\lemmaref{lemma:ainert_mod_set_satisfies_properties}, the set $\mathcal{S}$ satisfies all the required properties on which our results about the corresponding imprecise Markov chain rely. In particular, we can apply~\propositionref{prop:limit_is_lower_exp} to find that $\underline{s}_A^*\coloneqq \lim_{n\to+\infty} \underline{s}_A^{(n)}$ and $\overline{s}_A^*\coloneqq \lim_{n\to+\infty} \overline{s}_A^{(n)}$ are the lower and upper expected hitting times of the $A$-inert model, respectively. To finish the proof it therefore suffices to show that $\underline{h}_A^{(n)}=\underline{s}_A^{(n)}$ and $\overline{h}_A^{(n)}=\overline{s}_A^{(n)}$ for all $n\in\natswith$.

We will prove these equalities by induction, and start by giving the proof for the lower expected hitting time. To establish the induction base we note that $\underline{s}_A^{(0)}=\mathbb{I}_{A^c}=\underline{h}_A^{(0)}$ by definition. So let us now assume that $\underline{s}_A^{(n)}=\underline{h}_A^{(n)}$ for some $n\in\natswith$; we will show that the equality then also holds for $n+1$. To this end, note that for all $x\in A$ it holds that
\begin{equation*}
\underline{s}_A^{(n+1)}(x) = 0 = \underline{h}_A^{(n+1)}(x)\,,
\end{equation*}
and, for all $x\in A^c$, that
\begin{align*}
\underline{s}_A^{(n+1)}(x) &= \mathbb{I}_{A^c}(x) + \mathbb{I}_{A^c}(x)\cdot \bigl[\underline{S}\,\underline{s}_A^{(n)}\bigr](x) \\
 &= \mathbb{I}_{A^c}(x) + \mathbb{I}_{A^c}(x)\cdot \bigl[\underline{T}\,\underline{s}_A^{(n)}\bigr](x)\,,
\end{align*}
where we used~\lemmaref{lemma:ainert_mod_set_satisfies_properties}. This relies on the observation that $\underline{s}_A^{(n)}\in\gamblesextabove$ because $\underline{s}_A^{(n)}$ is non-negative; the argument for that is identical to the proof of~ \lemmaref{lemma:hitting_time_approx_sequence_nonnegative}, so we omit it here for brevity. It therefore follows from the induction hypothesis that
\begin{align*}
\underline{s}_A^{(n+1)}(x) &= \mathbb{I}_{A^c}(x) + \mathbb{I}_{A^c}(x)\cdot \bigl[\underline{T}\,\underline{s}_A^{(n)}\bigr](x) \\
 &=  \mathbb{I}_{A^c}(x) + \mathbb{I}_{A^c}(x)\cdot \bigl[\underline{T}\,\underline{h}_A^{(n)}\bigr](x) = \underline{h}_A^{(n+1)}(x)\,,
\end{align*}
where the last step used the definition of $\underline{h}_A^{(n+1)}$.

Because we have exhaustively covered all possible values of $x\in\states$, we conclude that, indeed, $\underline{s}_A^{(n+1)} = \underline{h}_A^{(n+1)}$, which concludes the induction step. We therefore conclude that
\begin{equation*}
\underline{s}_A^* = \lim_{n\to+\infty} \underline{s}_A^{(n)} = \lim_{n\to+\infty} \underline{h}_A^{(n)} = \underline{s}_A^{*}\,,
\end{equation*}
which concludes the proof for the lower expected hitting time. The proof that $\overline{s}_A^*=\overline{h}_A^*$ is completely analogous.
\end{proof}

In the remainder of this section of the appendix, we will until further notice \emph{implicitly} work with the $A$-inert modification of the imprecise Markov chain under consideration. That is, we simply assume that $\transmatset$ is $A$-inert, and therefore that the corresponding $\underline{T}$ and imprecise Markov chain are, as well. This has two important and immediate consequences:
\begin{lemma}\label{lemma:inert_set_remains_inert}
Suppose that $\underline{T}$ is $A$-inert. Then for all $x\in A$ and all $n\in\nats$, it holds that $\bigl[\underline{T}^n\mathbb{I}_A\bigr](x)=1$.
\end{lemma}
\begin{proof}
We consider two cases. First, if $A=\emptyset$ then the statement is vacuously true. Second, we can assume that $A\neq \emptyset$, and then fix any $x\in A$. Note that $\mathbb{I}_A \geq \mathbb{I}_{\{x\}}$. Hence by monotonicity of $\underline{T}$, it holds that 
\begin{equation*}
\bigl[\underline{T}\mathbb{I}_A\bigr](x) \geq \bigl[\underline{T}\mathbb{I}_{\{x\}}\bigr](x) = 1\,,
\end{equation*}
where the equality used~\corollaryref{cor:a_inert_on_indicators}. Moreover, it holds that
\begin{equation*}
\bigl[\underline{T}\mathbb{I}_A\bigr](x) \leq \bigl[\overline{T}\mathbb{I}_A\bigr](x) \leq \max_{y\in\states} \mathbb{I}_A(y) = 1\,,
\end{equation*}
where the equality holds because $A$ is non-empty by assumption. So, we conclude that $\bigl[\underline{T}\mathbb{I}_A\bigr](x)=1$, which proves the claim for $n=1$.

We will next show that, for any $n\in\nats$ and any $y\in A^c$, it holds that $\bigl[\underline{T}^n\mathbb{I}_A\bigr](y) \geq \mathbb{I}_A(y)$. This is vacuously true if $A^c=\emptyset$, and otherwise it holds because
\begin{equation*}
\bigl[\underline{T}^n\mathbb{I}_A\bigr](y) \geq \min_{z\in\states} \mathbb{I}_A(z) = 0 = \mathbb{I}_A(y)\,.
\end{equation*}

We now finish by induction; so, suppose the claim is true for $n\in\nats$. Then for all $z\in A$ it holds that $\bigl[\underline{T}^n\mathbb{I}_A\bigr](z)=1=\mathbb{I}_A(z)$. Moreover, as we established above, $\bigl[\underline{T}^n\mathbb{I}_A\bigr](y)\geq \mathbb{I}_A(y)$ for all $y\in A^c$. So we conclude that $\underline{T}^n\mathbb{I}_A \geq \mathbb{I}_A$. Then by the monotonicity of $\underline{T}$,
\begin{equation*}
\bigl[\underline{T}^{n+1}\mathbb{I}_A\bigr](x) = \bigl[\underline{T}\,\underline{T}^{n}\mathbb{I}_A\bigr](x) \geq \bigl[\underline{T}\mathbb{I}_A\bigr](x) = 1\,.
\end{equation*}
Moreover, using the properties of $\underline{T}^{n+1}$ and $\overline{T}^{n+1}$ we find that
\begin{equation*}
\bigl[\underline{T}^{n+1}\mathbb{I}_A\bigr](x) \leq
\bigl[\overline{T}^{n+1}\mathbb{I}_A\bigr](x) \leq \max_{y\in\states}\mathbb{I}_A(y) = 1\,,
\end{equation*}
and hence $\bigl[\underline{T}^{n+1}\mathbb{I}_A\bigr](x)=1$.
\end{proof}

\begin{lemma}\label{lemma:hit_target_monotone}
Suppose that $\underline{T}$ is $A$-inert. Then for all $x\in A^c$, if for some $n\in\nats$ it holds that $\bigl[\underline{T}^n\mathbb{I}_A\bigr](x)>0$, then $\bigl[\underline{T}^{n+1}\mathbb{I}_A\bigr](x)>0$.
\end{lemma}
\begin{proof}
We consider two cases. First, if $A^c=\emptyset$ then the statement is vacuously true. Second, we will assume that $A^c\neq\emptyset$.
For any $x\in A$, we know from~\lemmaref{lemma:inert_set_remains_inert} that $\bigl[\underline{T}\mathbb{I}_A\bigr](x)=1=\mathbb{I}_A(x)$. Moreover, for any $y\in A^c$ it holds that
\begin{equation*}
\bigl[\underline{T}\mathbb{I}_A\bigr](y) \geq \min_{z\in\states} \mathbb{I}_A(z) = 0 = \mathbb{I}_A(y)\,,
\end{equation*}
because $A^c\neq\emptyset$ by assumption. Therefore we conclude that $\underline{T}\mathbb{I}_A \geq \mathbb{I}_A$.

Now fix $x\in A^c$, which again is possible because $A^c\neq\emptyset$ by assumption, and suppose there is some $n\in\nats$ such that $\bigl[\underline{T}^n\mathbb{I}_A\bigr](x)>0$. Then by the monotonicity of $\underline{T}^n$,
\begin{equation*}
\bigl[\underline{T}^{n+1}\mathbb{I}_A\bigr](x) = \bigl[\underline{T}^{n}\underline{T}\mathbb{I}_A\bigr](x) \geq \bigl[\underline{T}^{n}\mathbb{I}_A\bigr](x) > 0\,.
\end{equation*}
\end{proof}

We are now ready to start setting up the proofs of (the upper expectation part of)~\theoremref{thm:lower_hitting_time_reach_and_equal} and~\corollaryref{cor:imprecise_hitting_time_is_minimal_system_solution}. First, we split the set $\states$ into three different classes, as follows:
\begin{definition}
Consider the sets $\mathcal{B}, \mathcal{U},\mathcal{Z}\subseteq\states$ defined as
\begin{itemize}
\item $\mathcal{B} \coloneqq \bigl\{ x\in A^c\,\big\vert\, \forall n\in\nats\,:\,\bigl[\underline{T}^n\mathbb{I}_A\bigr](x)=0\bigr\}$
\item $\mathcal{U}\coloneqq \bigl\{ x\in(A^c\setminus \mathcal{B})\,\big\vert\, \exists n\in\nats\,:\, \bigl[\overline{T}^n\mathbb{I}_\mathcal{B}\bigr](x)>0\bigr\}$
\item $\mathcal{Z} \coloneqq \states\setminus (\mathcal{B}\cup\mathcal{U})$
\end{itemize}
\end{definition}
In the remainder of this section, we will try to provide the reader with an intuitive interpretation of these sets and the results that we prove about them. To this end, we remark on an important interpretation that is well-known in the imprecise Markov chain literature: for any set of states $\mathcal{C}\subseteq\states$, any $x\in\states$, and any $n\in\natswith$, the quantity $\bigl[\underline{T}^n\mathbb{I}_\mathcal{C}\bigr](x)$ can be interpreted as the lower probability $\underline{P}\bigl(X_n\in\mathcal{C}\,\big\vert\,X_0=x\bigr)$ of being in any of the states in $\mathcal{C}$ at time $n$, given that the process started in state $x$ at time zero. An analogous statement with upper probabilities holds when using the corresponding upper transition operator. 
It should be noted that this interpretation holds for the imprecise Markov chain $\imcirr$ but that it does not hold for $\imchom$. Regardless, we will use this interpretation in the remainder to attempt to provide our results with some intuition.

For example, starting with the subsets of $\states$ defined above, the states in $\mathcal{B}$ are exactly those from which reaching $A$ in any finite number of steps has lower probability zero. 
Moreover, $\mathcal{U}$ contains the states from which reaching $\mathcal{B}$ in a finite number of steps has positive upper probability, but which are not themselves in $\mathcal{B}$. Because states in $\mathcal{B}$ have a zero lower probability to reach $A$, also the states in $\mathcal{U}$ in a sense inhibit reaching $A$. The states in $\mathcal{Z}$ are all other states, including those in $A$. 

As we will show below, the hitting time of the imprecise Markov chain is relatively well-behaved when it starts in $\mathcal{Z}$; in fact, the upper expected hitting then remains finite and we can uniquely characterise it. Conversely, as we will show, the upper expected hitting time when starting in either $\mathcal{B}$ or $\mathcal{U}$, will diverge to $+\infty$.

We start by stating a number of technical properties of these sets, which can in a certain sense be interpreted as closure properties. The first of these can intuitively be interpreted as saying that, once the process reaches $\mathcal{B}$, it remains in $\mathcal{B}$ with upper probability one.
\begin{lemma}\label{lemma:bad_states_are_closed}
For all $x\in\mathcal{B}$, it holds that $\bigl[\overline{T}\mathbb{I}_{\mathcal{B}}\bigr](x)=1$.
\end{lemma}
\begin{proof}
To avoid trivialities, we can assume without loss of generality that $\mathcal{B}$ is non-empty.

First consider any $x\in (A^c\setminus\mathcal{B})$. Because $x\notin \mathcal{B}$, there exists some $n_x\in\nats$ such that $\bigl[\underline{T}^{n_x}\mathbb{I}_A\bigr](x)>0$. Now let $n\coloneqq \max_{x\in (A^c\setminus\mathcal{B})}n_x$. Then, by~\lemmaref{lemma:hit_target_monotone}, for all $x\in (A^c\setminus\mathcal{B})$, because $n\geq n_x$, it holds that $\bigl[\underline{T}^{n}\mathbb{I}_A\bigr](x)>0$. Moreover, by~\lemmaref{lemma:inert_set_remains_inert}, clearly also $\bigl[\underline{T}^{n}\mathbb{I}_A\bigr](x)=1>0$ for all $x\in A$. We conclude that for all $x\in \mathcal{B}^c$ it holds that
\begin{equation*}
\bigl[\underline{T}^{n}\mathbb{I}_A\bigr](x)>0\,.
\end{equation*}
Now let $c\coloneqq \min_{x\in\mathcal{B}^c}\bigl[\underline{T}^{n}\mathbb{I}_A\bigr](x)$; then $c>0$ since $\mathcal{B}^c$ is finite. Define $f\in\gambles$ for all $x\in\states$ as
\begin{equation*}
f(x) \coloneqq \left\{\begin{array}{ll}
\nicefrac{1}{c}\bigl[\underline{T}^{n}\mathbb{I}_A\bigr](x) & \text{if $x\in \mathcal{B}^c$, and} \\
0 & \text{otherwise.}
\end{array}\right.
\end{equation*} 
Clearly it holds that $f \geq \mathbb{I}_{\mathcal{B}^c}$. Moreover, since for all $x\in\mathcal{B}$ it holds by definition that $\bigl[\underline{T}^n\mathbb{I}_A\bigr](x)=0$, we find
\begin{equation*}
f = \nicefrac{1}{c}\cdot\underline{T}^n\mathbb{I}_A\,.
\end{equation*}

By~\lemmaref{lemma:events_conjugate_complement}, the claim in this lemma's statement is equivalent to the claim that $\bigl[\underline{T}\mathbb{I}_{\mathcal{B}^c}\bigr](x)=0$ for all $x\in\mathcal{B}$. First, clearly $\bigl[\underline{T}\mathbb{I}_{\mathcal{B}^c}\bigr](x) \geq \min_{y\in\states} \mathbb{I}_{\mathcal{B}^c}(y)=0$. Now suppose \emph{ex absurdo} that there is some $x\in\mathcal{B}$ for which $\bigl[\underline{T}\mathbb{I}_{\mathcal{B}^c}\bigr](x)>0$. Then
\begin{align*}
\bigl[\underline{T}^{n+1}\mathbb{I}_A\bigr](x) &= \bigl[\underline{T}\,\underline{T}^{n}\mathbb{I}_A\bigr](x) \\
 &= c\bigl[\underline{T}\,f\bigr](x) \geq c\bigl[\underline{T}\,\mathbb{I}_{\mathcal{B}^c}\bigr](x) > 0\,,
\end{align*}
which implies $x\notin\mathcal{B}$, a contradiction.
\end{proof}

The next property can be interpreted as saying that, if it is possible to reach the set $\mathcal{U}$ when starting in some state $x\in A^c\setminus \mathcal{B}$, then $x$ must be in $\mathcal{U}$.
\begin{lemma}\label{lemma:unfortunate_reachable_is_unfortunate}
For all $x\in A^c\setminus \mathcal{B}$, if $\bigl[\overline{T}\mathbb{I}_\mathcal{U}\bigr](x)>0$, then it holds that $x\in\mathcal{U}$.
\end{lemma}
\begin{proof}
Fix $x\in A^c\setminus \mathcal{B}$ and suppose that $\bigl[\overline{T}\mathbb{I}_\mathcal{U}\bigr](x)>0$. This implies the existence of some $T\in\mathcal{T}$ such that $T(x,u)>0$ for some $u\in\mathcal{U}$. Because $u\in\mathcal{U}$, there is some $n\in\nats$ such that $\bigl[\overline{T}^n\mathbb{I}_\mathcal{B}\bigr](u)>0$. Therefore
\begin{align*}
\bigl[\overline{T}^{n+1}\mathbb{I}_\mathcal{B}\bigr](x) &\geq \bigl[T\overline{T}^n\mathbb{I}_\mathcal{B}\bigr](x) \\
 &= \sum_{y\in\states} T(x,y)\bigl[\overline{T}^n\mathbb{I}_\mathcal{B}\bigr](y) \\
 &\geq T(x,u)\bigl[\overline{T}^n\mathbb{I}_\mathcal{B}\bigr](u) > 0\,,
\end{align*}
which is the required condition to get $x\in\mathcal{U}$.
\end{proof}

The next result intuitively tells us that, if an element $u\in\mathcal{U}$ can reach $\mathcal{B}$ in $n_u>1$ steps, then there is also an element $v\in\mathcal{U}$ that can reach $\mathcal{B}$ in $n_u-1$ steps:
\begin{lemma}\label{lemma:unfortunate_set_induction_steps}
For every $u\in\mathcal{U}$, let $n_u\in\nats$ be the smallest number such that $\bigl[\overline{T}^{n_u}\mathbb{I}_\mathcal{B}\bigr](u)>0$; this $n_u$ exists by the definition of $\mathcal{U}$. Then, for all $u\in\mathcal{U}$, if $n_u>1$, there is some $v\in\mathcal{U}$ such that $n_v=n_u-1$, and clearly in that case $u\neq v$.
\end{lemma}
\begin{proof}
Fix any $u\in\mathcal{U}$ and suppose that $n_u>1$. 
Now, by the definition of $n_u$ and due to~\eqref{eq:lower_trans_reached}, there is some $T\in\mathcal{T}$ such that
\begin{equation*}
\bigl[T\overline{T}^{n_u-1}\mathbb{I}_\mathcal{B}\bigr](u) = \bigl[\overline{T}^{n_u}\mathbb{I}_\mathcal{B}\bigr](u) > 0\,.
\end{equation*}
Hence we have
\begin{equation}\label{lemma:unfortunate_set_induction_steps:eq:pos}
0 < \sum_{x\in\states} T(u,x) \bigl[\overline{T}^{n_u-1}\mathbb{I}_\mathcal{B}\bigr](x)\,.
\end{equation}
Let $\mathcal{V}\subseteq \states$ be defined as
\begin{equation*}
\mathcal{V} \coloneqq \Bigl\{ x\in\states\,:\, T(u,x)\bigl[\overline{T}^{n_u-1}\mathbb{I}_\mathcal{B}\bigr](x) > 0\Bigr\}\,.
\end{equation*}
That $\mathcal{V}\neq \emptyset$ is immediate from~\eqref{lemma:unfortunate_set_induction_steps:eq:pos}. We will next show that $\mathcal{V}\subseteq \mathcal{U}$. To this end, consider any $v\in\mathcal{V}$. Then, because $\mathcal{B}\subseteq A^c$, it follows from the monotonicity of $\overline{T}$ that
\begin{equation*}
\bigl[\overline{T}^{n_u-1}\mathbb{I}_{A^c}\bigr](v) \geq \bigl[\overline{T}^{n_u-1}\mathbb{I}_\mathcal{B}\bigr](v) > 0\,,
\end{equation*}
which, using~\lemmaref{lemma:events_conjugate_complement} and~\lemmaref{lemma:inert_set_remains_inert}, implies that $v\notin A$. Hence we know that $\mathcal{V}\subseteq A^c$. Next, suppose \emph{ex absurdo} that $v\in\mathcal{B}$. Then, because $T(u,v)>0$ since $v\in\mathcal{V}$, it follows that
\begin{equation*}
\bigl[\overline{T}\mathbb{I}_\mathcal{B}\bigr](u) \geq \bigl[T\mathbb{I}_\mathcal{B}\bigr](u) \geq T(u,v) > 0\,,
\end{equation*}
which would mean that $n_u=1$, a contradiction with the assumption made at the beginning of this proof. Hence $v\notin \mathcal{B}$. Finally, since $v\in\mathcal{V}$, it holds that $\bigl[\overline{T}^{n_u-1}\mathbb{I}_\mathcal{B}\bigr](v) > 0$, which is a sufficient condition to establish $v\in\mathcal{U}$.

So, indeed, we have established that $\emptyset\neq \mathcal{V}\subseteq \mathcal{U}$. Moreover, for all $v\in\mathcal{V}$, it holds that $n_v \leq n_u-1$ and, clearly, that $T(u,v)>0$. Now fix any $v\in \mathcal{V}$. 
We already established that $n_v\leq n_u-1$, and if $n_v=n_u-1$ then we are done. So suppose \emph{ex absurdo} that $n_v < n_u-1$. Then
\begin{align*}
\bigl[\overline{T}^{n_v+1}\mathbb{I}_\mathcal{B}\bigr](u) &\geq \bigl[T \overline{T}^{n_v}\mathbb{I}_\mathcal{B}\bigr](u) \\
 &= \sum_{w\in\states} T(u,w)\bigl[\overline{T}^{n_v}\mathbb{I}_\mathcal{B}\bigr](w) \\
  &\geq T(u,v)\bigl[\overline{T}^{n_v}\mathbb{I}_\mathcal{B}\bigr](v) > 0\,,
\end{align*}
where the second inequality used the fact that, for all $w\in\states$, $T(u,w)\geq 0$ because $T$ is a transition matrix, and $\bigl[\overline{T}^{n_v}\mathbb{I}_\mathcal{B}\bigr](w)\geq \min_{x\in\states} \mathbb{I}_\mathcal{B}(x)\geq 0$.
This implies that $n_u>n_v+1$ is not the smallest number for which $\bigl[\overline{T}^{n_u}\mathbb{I}_\mathcal{B}\bigr](u)>0$, a contradiction.
\end{proof}

As a consequence, if $\mathcal{U}$ is not empty, it contains at least one element that can reach $\mathcal{B}$ in a single step:
\begin{corollary}\label{cor:unfortunate_set_induction_base}
For all $u\in\mathcal{U}$, let $n_u$ be defined as in~\lemmaref{lemma:unfortunate_set_induction_steps}. Then, if $\mathcal{U}\neq \emptyset$, there is some $u\in\mathcal{U}$ with $n_u=1$.
\end{corollary}
\begin{proof}
Suppose that $\mathcal{U}\neq \emptyset$, and choose any $u\in\mathcal{U}$. If $n_u=1$ then we are done; otherwise, use~\lemmaref{lemma:unfortunate_set_induction_steps} to find $v\in\mathcal{U}$ with $n_v=n_u-1$. Then repeat until $n_v=1$.
\end{proof}

The following property can be interpreted as saying that, if the process starts in a state $x\in\mathcal{Z}$, it will remain in $\mathcal{Z}$ with lower probability one:
\begin{lemma}\label{lemma:non_bad_or_unfortunate_states_are_closed}
For all $x\in\mathcal{Z}$ it holds that $\bigl[\underline{T}\mathbb{I}_\mathcal{Z}\bigr](x)=1$.
\end{lemma}
\begin{proof}
For $x\in A\subseteq \mathcal{Z}$, it follows from the monotonicity and $A$-inertness of $\underline{T}$ that
\begin{equation*}
1\geq \bigl[\underline{T}\mathbb{I}_\mathcal{Z}\bigr](x) \geq \bigl[\underline{T}\mathbb{I}_A\bigr](x)=1\,.
\end{equation*}

So consider any $x\in\mathcal{Z}\cap A^c$. Because $x\notin \mathcal{B}$ and $x\notin \mathcal{U}$, we have $\bigl[\overline{T}\mathbb{I}_\mathcal{B}\bigr](x)=0$. Moreover, since $x\notin \mathcal{B}$ and $x\notin \mathcal{U}$, we find $\bigl[\overline{T}\mathbb{I}_\mathcal{U}\bigr](x)\leq 0$ due to~\lemmaref{lemma:unfortunate_reachable_is_unfortunate}; because also $\bigl[\overline{T}\mathbb{I}_\mathcal{U}\bigr](x) \geq \min_{y\in\states} \mathbb{I}_{\mathcal{U}}(y)\geq 0$ this implies $\bigl[\overline{T}\mathbb{I}_\mathcal{U}\bigr](x)= 0$. Now observe that since $\mathcal{B}\cap\mathcal{U}=\emptyset$, it holds that $\mathbb{I}_{\mathcal{B}\cup\mathcal{U}} = \mathbb{I}_\mathcal{B} + \mathbb{I}_\mathcal{U}$ and therefore, by the subadditivity of $\overline{T}$,
\begin{align*}
0 \leq \min_{y\in \states} \mathbb{I}_{\mathcal{B}\cup\mathcal{U}}(y) &\leq \bigl[\overline{T} \mathbb{I}_{\mathcal{B}\cup\mathcal{U}}\bigr](x) \\
 &\leq \bigl[\overline{T}\mathbb{I}_\mathcal{B}\bigr](x) + \bigl[\overline{T}\mathbb{I}_\mathcal{U}\bigr](x) = 0\,.
\end{align*}

From $\mathcal{Z}=(\mathcal{B}\cup\mathcal{U})^c$ together with~\lemmaref{lemma:events_conjugate_complement} we now conclude that indeed $\bigl[\underline{T}\mathbb{I}_\mathcal{Z}\bigr](x)=1$.
\end{proof}

The next property generalises the above result to any finite number of steps; intuitively, since if we start in $\mathcal{Z}$ and take one step we remain in $\mathcal{Z}$---this is essentially the previous statement---we also remain in $\mathcal{Z}$ (with lower probability one) after taking $n$ steps:
\begin{corollary}\label{cor:nice_remains_nice}
For all $n\in\nats$ and all $x\in\mathcal{Z}$, it holds that $\bigl[\underline{T}^n\mathbb{I}_\mathcal{Z}\bigr](x)=1$.
\end{corollary}
\begin{proof}
By~\lemmaref{lemma:non_bad_or_unfortunate_states_are_closed} we have for all $x\in\mathcal{Z}$ that $\bigl[\underline{T}\mathbb{I}_\mathcal{Z}\bigr](x)=1=\mathbb{I}_\mathcal{Z}(x)$ . Moreover, for all $y\in\states\setminus\mathcal{Z}$ we have $\bigl[\underline{T}\mathbb{I}_\mathcal{Z}\bigr](y)\geq \min_{z\in\states} \mathbb{I}_\mathcal{Z}(z) \geq 0 = \mathbb{I}_\mathcal{Z}(y)$, from which we conclude that
\begin{equation*}
\underline{T}\mathbb{I}_\mathcal{Z} \geq \mathbb{I}_\mathcal{Z}\,.
\end{equation*}
Now fix any $n\in\nats$. Then by the monotonicity of $\underline{T}^n$, we find that for all $x\in\mathcal{Z}$,
\begin{align*}
\bigl[\underline{T}^{n+1}\mathbb{I}_\mathcal{Z}\bigr](x) &= \bigl[\underline{T}^{n}\underline{T}\mathbb{I}_\mathcal{Z}\bigr](x) \\
 &\geq \bigl[\underline{T}^{n}\mathbb{I}_\mathcal{Z}\bigr](x) \\
 &\geq \bigl[\underline{T}^{n-1}\underline{T}\mathbb{I}_\mathcal{Z}\bigr](x) 
 \geq \cdots \geq \bigl[\underline{T}\mathbb{I}_\mathcal{Z}\bigr](x) = 1\,.
\end{align*}
Conversely, it holds that
\begin{equation*}
\bigl[\underline{T}^{n+1}\mathbb{I}_\mathcal{Z}\bigr](x)\leq \bigl[\overline{T}^{n+1}\mathbb{I}_\mathcal{Z}\bigr](x) \leq \max_{y\in\states}\mathbb{I}_{\mathcal{Z}}(y)=1\,,
\end{equation*}
where we made use of~\lemmaref{lemma:composition_upper_dominates_composition_lower}.
\end{proof}

The following consequence of~\lemmaref{lemma:non_bad_or_unfortunate_states_are_closed} is important, and its interpretation can be explained using the interpretation of $\underline{T}$ as encoding the 1 time step lower expectation of an imprecise Markov chain (see~\sectionref{subsec:measure_impr_markov}). It can be interpreted as saying that, because if we start in $\mathcal{Z}$ we remain in $\mathcal{Z}$, the 1 time-step conditional expectation of a function $f$, when starting in $\mathcal{Z}$, only depends on the value of $f$ on $\mathcal{Z}$:
\begin{corollary}\label{cor:nice_sets_functions_depend_on_nice}
For all $f,g\in\gambles$ such that $f(x)=g(x)$ for all $x\in\mathcal{Z}$, it holds for all $x\in\mathcal{Z}$ that
\begin{equation*}
\bigl[\underline{T}\,f\bigr](x) = \bigl[\underline{T}\,g\bigr](x)\quad\text{and}\quad \bigl[\overline{T}f\bigr](x) = \bigl[\overline{T}g\bigr](x)\,.
\end{equation*}
\end{corollary}
\begin{proof}
Let $C\coloneqq \max_{x\in\states} \abs{f-g}(x)$; then clearly $0\leq C$.

Fix any $x\in\mathcal{Z}$. Then, using~\lemmaref{lemma:lower_trans_basic_properties},
\begin{align*}
\abs{\bigl[\underline{T}\,f\bigr](x) - \bigl[\underline{T}\,g\bigr](x)}
 &\leq \bigl[\overline{T}\abs{f-g}\bigr](x) \\
 &= \bigl[\overline{T}\bigl(\mathbb{I}_{\mathcal{Z}^c}\hprod\abs{f-g}\bigr)\bigr](x) \\
 &\leq \bigl[\overline{T}\bigl(\mathbb{I}_{\mathcal{Z}^c}C\bigr)\bigr](x) \\
 &= C\bigl[\overline{T}\bigl(\mathbb{I}_{\mathcal{Z}^c}\bigr)\bigr](x) \\
 &= 0\,,
\end{align*}
where the last step used~\lemmaref{lemma:non_bad_or_unfortunate_states_are_closed} and~\lemmaref{lemma:events_conjugate_complement}. Hence indeed $\bigl[\underline{T}\,f\bigr](x) = \bigl[\underline{T}\,g\bigr](x)$. The proof for $\overline{T}$ is completely analogous.
\end{proof}

The next result intuitively tells us something about why the set $\mathcal{Z}$ is particularly well-behaved in the context of hitting-times of the set $A$: when starting in any state $x\in\mathcal{Z}\cap A^c$, there is a strictly positive lower probability of reaching $A$ in some finite number of steps: 
\begin{lemma}\label{lemma:nice_sets_reach_A}
For all $x\in \mathcal{Z}\cap A^c$, there is some $n\in\nats$ such that $\bigl[\underline{T}^n\mathbb{I}_A\bigr](x)>0$.
\end{lemma}
\begin{proof}
First, clearly, $\bigl[\underline{T}^n\mathbb{I}_A\bigr](x) \geq \min_{y\in\states} \mathbb{I}_A(y)\geq 0$ for all $n\in\nats$ and all $x\in\states$.

Moreover, by the definition of $\mathcal{B}$, for any $x\in A^c$ it holds that if $\bigl[\underline{T}^n\mathbb{I}_A\bigr](x)=0$ for all $n\in\nats$, then $x\in\mathcal{B}$. Thus the result follows because $\mathcal{Z}\cap \mathcal{B}=\emptyset$.
\end{proof}

The next statement essentially states the converse to the previous property; in short, it states that when starting in $\mathcal{Z}\cap A^c$, we can leave $A^c$ with positive upper probability in a finite number of steps:
\begin{corollary}\label{cor:nice_set_leaves_Ac}
For all $x\in\mathcal{Z}\cap A^c$, there is some $n\in\nats$ such that $\bigl[\overline{T}^n\mathbb{I}_{A^c}\bigr](x) < 1$. Moreover, also $\bigl[\overline{T}^{m}\mathbb{I}_{A^c}\bigr](x) < 1$ for all $m\in\nats$ with $n\leq m$.
\end{corollary}
\begin{proof}
Fix $x\in\mathcal{Z}\cap A^c$. By~\lemmaref{lemma:nice_sets_reach_A}, there exists some $n\in\nats$ such that $\bigl[\underline{T}^n\mathbb{I}_A\bigr](x)>0$. By~\lemmaref{lemma:events_conjugate_complement}, this implies $\bigl[\overline{T}^n\mathbb{I}_{A^c}\bigr](x)\neq 1$. But $\bigl[\overline{T}^n\mathbb{I}_{A^c}\bigr](x) \leq \max_{y\in\states}\mathbb{I}_{A^c}(y)\leq 1$.

The second claim is by induction; so, suppose that for some $n\in\nats$ it holds that $\bigl[\overline{T}^n\mathbb{I}_{A^c}\bigr](x) < 1$.
By~\lemmaref{lemma:events_conjugate_complement}, this implies $\bigl[\underline{T}^n\mathbb{I}_{A}\bigr](x)\neq 0$. Since $\bigl[\underline{T}^n\mathbb{I}_{A}\bigr](x)\geq \min_{y\in\states}\mathbb{I}_{A}(y)\geq 0$, we conclude $\bigl[\underline{T}^n\mathbb{I}_{A}\bigr](x)>0$. Due to~\lemmaref{lemma:hit_target_monotone}, this implies $\bigl[\underline{T}^{n+1}\mathbb{I}_{A}\bigr](x)>0$ because $\underline{T}$ is $A$-inert. 

Then due to~\lemmaref{lemma:events_conjugate_complement} we have that $\bigl[\overline{T}^{n+1}\mathbb{I}_{A^c}\bigr](x)\neq 1$  and, since $\bigl[\overline{T}^{n+1}\mathbb{I}_{A^c}\bigr](x)\leq \max_{y\in\states}\mathbb{I}_{A^c}(y)\leq 1$, we find that $\bigl[\overline{T}^{n+1}\mathbb{I}_{A^c}\bigr](x)<1$.
\end{proof}

The following property combines some of the above intuitive statements. It can be interpreted as saying that, if we start in $\mathcal{Z}$, the (upper) probability of moving to any state in $A^c$ is the same as the (upper) probability of moving to a state in $\mathcal{Z}\cap A^c$; this follows because when starting in $\mathcal{Z}$, we always stay in $\mathcal{Z}$.
\begin{lemma}\label{lemma:leave_nice_set_complement_stays_nice}
For all $n\in\nats$ and all $x\in \mathcal{Z}$ it holds that $\bigl[\overline{T}^n\mathbb{I}_{\mathcal{Z}\cap A^c}\bigr](x) = \bigl[\overline{T}^n\mathbb{I}_{A^c}\bigr](x)$\,.
\end{lemma}
\begin{proof}
We give a proof by induction. First, note that clearly $\mathbb{I}_{\mathcal{Z}\cap A^c}(x) = \mathbb{I}_{A^c}(x)$ for all $x\in\mathcal{Z}$. Hence, by~\corollaryref{cor:nice_sets_functions_depend_on_nice}, we have for all $x\in\mathcal{Z}$ that
\begin{equation*}
[\overline{T}\mathbb{I}_{\mathcal{Z}\cap A^c}\bigr](x) = \bigl[\overline{T}\mathbb{I}_{A^c}\bigr](x)\,.
\end{equation*}
This provides the induction base. Now suppose the statement is true for some $n\in\nats$. Then, by~\corollaryref{cor:nice_sets_functions_depend_on_nice}, we have for all $x\in\mathcal{Z}$ that
\begin{align*}
\bigl[\overline{T}^{n+1}\mathbb{I}_{\mathcal{Z}\cap A^c}\bigr](x) &= \bigl[\overline{T}\bigl(\overline{T}^n\mathbb{I}_{\mathcal{Z}\cap A^c}\bigr)\bigr](x) \\
 &= \bigl[\overline{T}\bigl(\overline{T}^n\mathbb{I}_{A^c}\bigr)\bigr](x) = \bigl[\overline{T}^{n+1}\mathbb{I}_{A^c}\bigr](x)\,.
\end{align*}
\end{proof}

At this point, it should hopefully be intuitively clear from the above results that the behaviour of the imprecise stochastic process on the states $\mathcal{Z}$ is relatively self-contained, in the sense that the process cannot leave these states once it enters them. In order to derive an expression for the upper expected hitting time when starting in $\mathcal{Z}$, we start below by explicitly and separately describing the behaviour on this part of the domain.

To this end, let $\mathcal{L}(\mathcal{Z}\cap A^c)$ be the (vector-)space of all functions $f:(\mathcal{Z}\cap A^c) \to \reals$. For any $f\in \mathcal{L}(\mathcal{Z}\cap A^c)$, let $f^\uparrow$ denote the \emph{zero-extension} of $f$ into $\gambles$, defined for all $x\in\states$ as
\begin{equation*}
f^\uparrow(x) \coloneqq \left\{\begin{array}{ll}
f(x) & \text{if $x\in\mathcal{Z}\cap A^c$, and } \\
0 & \text{otherwise.}
\end{array}\right.
\end{equation*}
Let $\overline{Q}:\mathcal{L}(\mathcal{Z}\cap A^c)\to \mathcal{L}(\mathcal{Z}\cap A^c)$ be the non-linear operator that is defined, for all $f\in \mathcal{L}(\mathcal{Z}\cap A^c)$ and all $x\in\mathcal{Z}\cap A^c$, as
\begin{equation*}
\bigl[\overline{Q}f\bigr](x) \coloneqq \bigl[\overline{T}f^\uparrow\bigr](x)\,.
\end{equation*}
Provide $\mathcal{L}(\mathcal{Z}\cap A^c)$ with the supremum norm and any non-negatively homogeneous operator on this space with the induced operator norm; as with our comments in~\appendixref{appendix:prelim_proofs}, we note that non-negative homogeneity is a sufficient condition for this to induce a norm. Moreover, it follows from~\lemmaref{lemma:uQ_nonneg_homogen} below that in particular $\overline{Q}$ is non-negatively homogeneous.

We next state a number of properties of $\overline{Q}$, which it essentially inherits from $\overline{T}$.
\begin{lemma}\label{lemma:uQ_super_absolute}
For all $f,g\in\mathcal{L}(\mathcal{Z}\cap A^c)$ it holds that
\begin{equation*}
\abs{\overline{Q}f - \overline{Q}g} \leq \overline{Q}\abs{f-g}\,.
\end{equation*}
\end{lemma}
\begin{proof}
By the definition of $\overline{Q}$ and using~\lemmaref{lemma:lower_trans_basic_properties},
\begin{align*}
\abs{\overline{Q}f - \overline{Q}g} = \abs{\overline{T}f^\uparrow- \overline{T}g^\uparrow} &\leq \overline{T}\abs{f^\uparrow-g^\uparrow} \\
 &= \overline{T}\left(\abs{f-g}^\uparrow\right) = \overline{Q}\abs{f-g}\,.
\end{align*}
\end{proof}

\begin{lemma}\label{lemma:uQ_monotone}
For all $f,g\in\mathcal{L}(\mathcal{Z}\cap A^c)$ it holds that
\begin{equation*}
f \leq g \quad\Rightarrow\quad \overline{Q}f \leq \overline{Q}g\,.
\end{equation*}
\end{lemma}
\begin{proof}
By the definition of $\overline{Q}$ and using~\lemmaref{lemma:lower_trans_basic_properties},
\begin{align*}
\overline{Q}f = \overline{T}f^\uparrow \leq \overline{T} g^\uparrow = \overline{Q}g\,.
\end{align*}
\end{proof}

\begin{corollary}\label{cor:uQ_comp_super_abs}
For all $f\in\mathcal{L}(\mathcal{Z}\cap A^c)$ and all $n\in\nats$ it holds that $\abs{\overline{Q}^nf} \leq \overline{Q}^n\abs{f}$.
\end{corollary}
\begin{proof}
We give a proof by induction. First, let $g\in\mathcal{L}(\mathcal{Z}\cap A^c)$ be identically zero, $g=0$. Then, clearly, also $\overline{Q}g = 0$. Then for any $f\in\mathcal{L}(\mathcal{Z}\cap A^c)$, using~\lemmaref{lemma:uQ_super_absolute} we get
\begin{equation*}
\abs{\overline{Q}f} = \abs{\overline{Q}f - \overline{Q}g} \leq \overline{Q}\abs{f-g} = \overline{Q}\abs{f}\,,
\end{equation*} 
which provides the induction base. Now suppose that the statement is true for $n\in\nats$; we will show that it is also true for $n+1$. In particular, for any $f\in\mathcal{L}(\mathcal{Z}\cap A^c)$, we have
\begin{align*}
\abs{\overline{Q}^{n+1}f} = \abs{\overline{Q}\bigl(\overline{Q}^nf\bigr)} \leq \overline{Q}\abs{\overline{Q}^nf} \leq \overline{Q}\overline{Q}^n\abs{f} = \overline{Q}^{n+1}\abs{f}\,, 
\end{align*}
where the first inequality uses the argument for the induction base above, and the second inequality uses the induction hypothesis together with the monotonicity of $\overline{Q}$, which is established in~\lemmaref{lemma:uQ_monotone}.
\end{proof}

\begin{lemma}\label{lemma:uQ_nonneg_homogen}
For all $f\in\mathcal{L}(\mathcal{Z}\cap A^c)$ and all $\lambda\in\reals$ with $\lambda\geq 0$, it holds that
\begin{equation*}
\overline{Q}(\lambda f) = \lambda\overline{Q}f\,.
\end{equation*}
\end{lemma}
\begin{proof}
By the definition of $\overline{Q}$ and using~\lemmaref{lemma:lower_trans_basic_properties},
\begin{align*}
\overline{Q}(\lambda f) = \overline{T}\Bigl[(\lambda f)^\uparrow\Bigr] =  \overline{T}(\lambda f^\uparrow) = \lambda\overline{T} f^\uparrow = \lambda\overline{Q}f\,.
\end{align*}
\end{proof}

The following property allows us to reduce compositions of $\overline{Q}$ to compositions of $\overline{T}$:
\begin{lemma}\label{lemma:uQ_composition_interchange}
For all $n\in\nats$ and all $f\in \mathcal{L}(\mathcal{Z}\cap A^c)$ it holds that $\bigl[\overline{Q}^nf\bigr](x) = \bigl[\overline{T}^n(f^\uparrow)\bigr](x)$ for all $x\in\mathcal{Z}\cap A^c$.
\end{lemma}
\begin{proof}
We give a proof by induction. First, by definition, $\bigl[\overline{Q}f\bigr](x) = \bigl[\overline{T}f^\uparrow\bigr](x)$ for all $x\in\mathcal{Z}\cap A^c$. We can equivalently write this using the restriction of $\overline{T}f^\uparrow$ to $\mathcal{Z}\cap A^c$; thus,
\begin{equation*}
\overline{Q}f = \bigl[\overline{T}f^\uparrow\bigr]\big\vert_{\mathcal{Z}\cap A^c}\,.
\end{equation*}
We will now show that for all $x\in\mathcal{Z}$, it holds that
\begin{equation}\label{lemma:uQ_composition_interchange:eq:restrict_extend_onestep}
\bigl[\overline{T}f^\uparrow\bigr](x) = \left(\bigl[\overline{T}f^\uparrow\bigr]\big\vert_{\mathcal{Z}\cap A^c}\right)^\uparrow(x)\,,
\end{equation}
which is to say, if we restrict $\overline{T}f^\uparrow$ to $\mathcal{Z}\cap A^c$ and then zero-extend this restricted function back into $\gambles$, the values on $\mathcal{Z}$ remain the same. Clearly this is true for $x\in\mathcal{Z}\cap A^c$. Now note that because $A\subseteq \mathcal{Z}$ it holds that $\mathcal{Z}\setminus A^c = A$. Hence because the model is $A$-inert, we have for $x\in A$ that
\begin{equation*}
\bigl[\overline{T}f^\uparrow\bigr](x) = f^\uparrow(x) = 0 = \left(\bigl[\overline{T}f^\uparrow\bigr]\big\vert_{\mathcal{Z}\cap A^c}\right)^\uparrow(x)\,,
\end{equation*}
so we conclude that~\eqref{lemma:uQ_composition_interchange:eq:restrict_extend_onestep} indeed holds.

Now, starting from~\eqref{lemma:uQ_composition_interchange:eq:restrict_extend_onestep}, we can repeatedly apply~\corollaryref{cor:nice_sets_functions_depend_on_nice} to find that, for any $n\in\nats$,
\begin{equation}\label{lemma:uQ_composition_interchange:eq:restrict_extend}
\left[\overline{T}^n\left( \overline{T}f^\uparrow \right)\right](x) = \left[\overline{T}^n\left( \left(\bigl[\overline{T}f^\uparrow\bigr]\big\vert_{\mathcal{Z}\cap A^c}\right)^\uparrow \right)\right](x)\,.
\end{equation}

Now suppose the statement is true for some $n\in\nats$. Then for all $x\in\mathcal{Z}\cap A^c$,
\begin{align*}
\bigl[\overline{Q}^{n+1}f\bigr](x) &= \bigl[\overline{Q}^n\bigl(\overline{Q}f\bigr)\bigr](x) \\
 &= \left[\overline{T}^n\Bigl[\bigl(\overline{Q}f\bigr)^\uparrow\Bigr]\right](x) \\
 &= \left[\overline{T}^n\left[\left(\bigl[\overline{T}f^\uparrow\bigr]\big\vert_{\mathcal{Z}\cap A^c}\right)^\uparrow\right]\right](x) \\
 &= \left[\overline{T}^{n}\left(\overline{T}\bigl(f^\uparrow\bigr)\right)\right](x) \\
 &= \left[\overline{T}^{n+1}f^\uparrow\right](x)\,,
\end{align*}
where the second step used the induction hypothesis, the third step used the definition of $\overline{Q}$, and the fourth step used~\eqref{lemma:uQ_composition_interchange:eq:restrict_extend}.
\end{proof}

We need the following three properties of the norm of $\overline{Q}$ and its self-compositions:
\begin{lemma}\label{lemma:uQ_norm_factors}
For all $f\in\mathcal{L}(\mathcal{Z}\cap A^c)$ and all $n\in\nats$, it holds that
\begin{equation*}
\norm{\overline{Q}^nf} \leq \norm{\overline{Q}^n}\norm{f}\,.
\end{equation*}
\end{lemma}
\begin{proof}
Fix $f\in\mathcal{L}(\mathcal{Z}\cap A^c)$ and $n\in\nats$. First we note that, for any $\lambda\in\reals$ with $\lambda\geq 0$, by repeated application of \lemmaref{lemma:uQ_nonneg_homogen}, we have
\begin{equation*}
\overline{Q}^n(\lambda f) = \lambda \overline{Q}^nf\,.
\end{equation*}
Now define $g\in \mathcal{L}(\mathcal{Z}\cap A^c)$ as
\begin{equation*}
g \coloneqq \left\{\begin{array}{ll}
\frac{1}{\norm{f}}f & \text{if $\norm{f}\neq 0$, and} \\
1 & \text{otherwise.}
\end{array}\right.
\end{equation*}
Then clearly $\norm{g}\leq 1$ and, because $\norm{f}=0$ if and only if $f=0$, we have $\norm{f}g = f$. Therefore,
\begin{align*}
\norm{\overline{Q}^nf} = \norm{\overline{Q}^n(\norm{f}g)} &=  \norm{\norm{f}\overline{Q}^ng} \\
 &= \norm{f}\norm{\overline{Q}^ng} \leq \norm{f}\norm{\overline{Q}^n}\,,
\end{align*}
where the inequality used the fact that $\norm{g}\leq 1$ and the definition of the induced operator norm.
\end{proof}

\begin{lemma}\label{lemma:uQ_unit_norm}
For all $n\in\nats$ it holds that $\norm{\overline{Q}^n}\leq 1$.
\end{lemma}
\begin{proof}
Using~\corollaryref{cor:uQ_comp_super_abs}, for any $f\in\mathcal{L}(\mathcal{Z}\cap A^c)$ and any $x\in\mathcal{Z}\cap A^c$ it holds that
\begin{equation*}
\abs{\overline{Q}^n f}(x) \leq \Bigl[\overline{Q}^n\abs{f}\Bigr](x)\,.
\end{equation*}
Moreover, due to~\lemmaref{lemma:uQ_composition_interchange}, we have
\begin{align*}
\Bigl[\overline{Q}^n\abs{f}\Bigr](x) &= \Bigl[\overline{T}^n\bigl(\abs{f}^\uparrow\bigr)\Bigr](x) \\
 &\leq \max_{y\in\states} \abs{f}^\uparrow(y) 
 = \max_{y\in\mathcal{Z}\cap A^c} \abs{f}(y) = \norm{f}\,,
\end{align*}
where we used~\lemmaref{lemma:lower_trans_composite_basic_properties} for the inequality. Hence we have $\abs{\overline{Q}^n f}(x) \leq \norm{f}$. Because $x\in\mathcal{Z}\cap A^c$ was arbitrary, this implies $\norm{\overline{Q}^nf} \leq \norm{f}$, which in turn implies $\norm{\overline{Q}^n} \leq 1$.
\end{proof}

\begin{lemma}\label{lemma:uQ_iter_contractive}
There exists some $n\in\nats$ such that $\norm{\overline{Q}^n}<1$.
\end{lemma}
\begin{proof}
Let $1$ denote the function in $\mathcal{L}(\mathcal{Z}\cap A^c)$ with constant value 1. Then for all $f\in \mathcal{L}(\mathcal{Z}\cap A^c)$ with $\norm{f}\leq 1$,
\begin{equation}\label{lemma:uQ_iter_contractive:eq:constant_best}
f \leq \abs{f} \leq 1\,.
\end{equation}

For all $x\in\mathcal{Z}\cap A^c$, due to~\corollaryref{cor:nice_set_leaves_Ac}, there is some $n_x\in\nats$ such that $\bigl[\overline{T}^{m}\mathbb{I}_{A^c}\bigr](x)<1$ for all $m\in\nats$ with $n_x\leq m$. Let $n\coloneqq \max_{x\in\mathcal{Z}\cap A^c} n_x$. Then  $\bigl[\overline{T}^{n}\mathbb{I}_{A^c}\bigr](x)<1$ for all $x\in\mathcal{Z}\cap A^c$. Now let $C\coloneqq \max_{x\in\mathcal{Z}\cap A^c} \bigl[\overline{T}^{n}\mathbb{I}_{A^c}\bigr](x)$; then $C<1$ since $\mathcal{Z}\cap A^c$ is finite.

Now fix any $f\in \mathcal{L}(\mathcal{Z}\cap A^c)$ with $\norm{f}\leq 1$. Then we find that
\begin{align*}
\norm{\overline{Q}^nf} &= \max_{x\in\mathcal{Z}\cap A^c} \abs{ \bigl[\overline{Q}^n f\bigr](x) } \\
 &\leq \max_{x\in\mathcal{Z}\cap A^c}  \bigl[\overline{Q}^n \abs{f}\bigr](x) \\
 &\leq \max_{x\in\mathcal{Z}\cap A^c}  \bigl[\overline{Q}^n 1\bigr](x) \\
 &= \max_{x\in\mathcal{Z}\cap A^c}  \bigl[\overline{T}^n 1^\uparrow\bigr](x) \\
 &= \max_{x\in\mathcal{Z}\cap A^c}  \bigl[\overline{T}^n \mathbb{I}_{\mathcal{Z}\cap A^c}\bigr](x) \\
 &= \max_{x\in\mathcal{Z}\cap A^c}  \bigl[\overline{T}^n \mathbb{I}_{A^c}\bigr](x) \\
 &= C\,,
\end{align*}
where the first inequality uses~\corollaryref{cor:uQ_comp_super_abs}, the second inequality uses~\eqref{lemma:uQ_iter_contractive:eq:constant_best} and~\lemmaref{lemma:uQ_monotone} (repeatedly),
the second equality uses~\lemmaref{lemma:uQ_composition_interchange} and the fourth equality used~\lemmaref{lemma:leave_nice_set_complement_stays_nice}.

Then,
\begin{align*}
\norm{\overline{Q}^n} &= \sup\Bigl\{ \norm{\overline{Q}^n f}\,\Big\vert\, f\in \mathcal{L}(\mathcal{Z}\cap A^c), \norm{f} \leq 1\Bigr\} \\
 &\leq \sup\Bigl\{ C\,\Big\vert\, f\in \mathcal{L}(\mathcal{Z}\cap A^c), \norm{f} \leq 1\Bigr\} \\
 &= C < 1\,.
\end{align*}
\end{proof}

We now have all the pieces to characterise the upper expected hitting time of the imprecise stochastic process when it starts in $\mathcal{Z}$. To this end, define the map $\overline{\mathbf{Q}}:\mathcal{L}(\mathcal{Z}\cap A^c)\to \mathcal{L}(\mathcal{Z}\cap A^c)$ for all $f\in \mathcal{L}(\mathcal{Z}\cap A^c)$ as
\begin{equation*}
\overline{\mathbf{Q}}(f) \coloneqq 1 + \overline{Q}f\,,
\end{equation*}
where $1$ is the function in $\mathcal{L}(\mathcal{Z}\cap A^c)$ with constant value one. For any $n\in\nats$, let $\overline{\mathbf{Q}}^n$ denote the $n$-fold composition of $\overline{\mathbf{Q}}$ with itself.

We will show below that there is a $q\in \mathcal{L}(\mathcal{Z}\cap A^c)$ describing the upper expected hitting time when starting in $\mathcal{Z}\cap A^c$. The following proposition tells us how to find this $q$.
\begin{proposition}\label{prop:nice_set_fixed_point}
There exists a unique $q\in \mathcal{L}(\mathcal{Z}\cap A^c)$ such that $\overline{\mathbf{Q}}(q) = q$. Moreover, $\lim_{n\to+\infty} \overline{\mathbf{Q}}^n(f) = q$ for all $f\in \mathcal{L}(\mathcal{Z}\cap A^c)$.
\end{proposition}
\begin{proof}
First, fix any $f,g\in \mathcal{L}(\mathcal{Z}\cap A^c)$. Note that
\begin{equation}\label{prop:nice_set_fixed_point:eq:one_step}
\abs{\overline{\mathbf{Q}}(f) - \overline{\mathbf{Q}}(g)} = \abs{\overline{Q}f - \overline{Q}g} \leq \overline{Q}\abs{f-g}\,,
\end{equation}
using the definition of $\overline{\mathbf{Q}}$ and~\lemmaref{lemma:uQ_super_absolute}. We will now show through induction that, in fact, for all $n\in\nats$, it holds that
\begin{equation*}
\abs{\overline{\mathbf{Q}}^n(f) - \overline{\mathbf{Q}}^n(g)} \leq \overline{Q}^n\abs{f-g}\,.
\end{equation*}
We have already derived the induction base above; so now assume that this is true for some $n\in\nats$. Then,
\begin{align*}
\abs{\overline{\mathbf{Q}}^{n+1}(f) - \overline{\mathbf{Q}}^{n+1}(g)} &= \abs{\overline{\mathbf{Q}}^n\bigl(\overline{\mathbf{Q}}(f)\bigr) - \overline{\mathbf{Q}}^n\bigl(\overline{\mathbf{Q}}(g)\bigr)} \\
 &\leq \overline{Q}^n\abs{\overline{\mathbf{Q}}(f) - \overline{\mathbf{Q}}(g)} \\
 &\leq \overline{Q}^{n+1}\abs{f-g}\,,
\end{align*}
where the first inequality used the induction hypothesis, and the second inequality used~\eqref{prop:nice_set_fixed_point:eq:one_step} together with the monotonicity of $\overline{Q}^n$---this last property follows straightforwardly from the monotonicity of $\overline{Q}$ in~\lemmaref{lemma:uQ_monotone}.

Hence, for any $n\in\nats$ and any $f,g\in\mathcal{L}(\mathcal{Z}\cap A^c)$, using~\lemmaref{lemma:uQ_norm_factors}, we conclude that
\begin{equation}\label{prop:nice_set_fixed_point:eq:contraction}
\norm{\overline{\mathbf{Q}}^n(f) - \overline{\mathbf{Q}}^n(g)} \leq \norm{\overline{Q}^n\abs{f-g}} \leq \norm{\overline{Q}^n}\norm{f-g}\,.
\end{equation}
By~\lemmaref{lemma:uQ_iter_contractive}, there is some $n\in\nats$ such that $\norm{\overline{Q}^n}<1$. Using~\eqref{prop:nice_set_fixed_point:eq:contraction}, this implies that $\overline{\mathbf{Q}}^n$ is a contraction mapping on $\mathcal{L}(\mathcal{Z}\cap A^c)$. Using a generalisation of the Banach fixed point theorem~\cite[Theorem 2.4]{almezel2014topics}, $\overline{\mathbf{Q}}$ (and not only $\overline{\mathbf{Q}}^n$!) therefore has a unique fixed point $q\in\mathcal{L}(\mathcal{Z}\cap A^c)$, i.e. this is the unique point for which $\overline{\mathbf{Q}}(q) = q$. This finishes the proof of the first claim.

For the second claim, i.e. that $\lim_{m\to+\infty} \overline{\mathbf{Q}}^m(f) = q$ for all $f\in\mathcal{L}(\mathcal{Z}\cap A^c)$, first note that, since $\overline{\mathbf{Q}}^n$ is a contraction, by the Banach fixed point theorem, there is a unique $q^*\in\mathcal{L}(\mathcal{Z}\cap A^c)$ such that $\overline{\mathbf{Q}}^n(q^*)=q^*$ and, for all $f\in\mathcal{L}(\mathcal{Z}\cap A^c)$, it holds that 
\begin{equation*}
\lim_{m\to+\infty}\overline{\mathbf{Q}}^{mn}(f) = q^*\,.
\end{equation*}
Now, first note that, since $\overline{\mathbf{Q}}(q)=q$, we quickly find that $\overline{\mathbf{Q}}^n(q)=q$ and, therefore, that $q=q^*$.

Now choose any $k,m\in\nats$. Then,
\begin{align*}
\norm{\overline{\mathbf{Q}}^{mn+k}(f) - q} 
 &= \norm{\overline{\mathbf{Q}}^{mn+k}(f) - \overline{\mathbf{Q}}^k(q)} \\
 &= \norm{\overline{\mathbf{Q}}^{k}\left( \overline{\mathbf{Q}}^{mn}(f) \right) - \overline{\mathbf{Q}}^k(q)} \\
 &\leq \norm{\overline{Q}^k}\norm{ \overline{\mathbf{Q}}^{mn}(f) - q} \\
 &\leq \norm{ \overline{\mathbf{Q}}^{mn}(f) - q}\,,
\end{align*}
where the first inequality used~\eqref{prop:nice_set_fixed_point:eq:contraction} and the second inequality used~\lemmaref{lemma:uQ_unit_norm}. Therefore, and because
\begin{equation*}
\lim_{m\to+\infty} \norm{\overline{\mathbf{Q}}^{mn}(f)-q} = \lim_{m\to+\infty} \norm{\overline{\mathbf{Q}}^{mn}(f)-q^*} = 0\,,
\end{equation*}
we conclude that also $\lim_{\ell\to+\infty}\norm{\overline{\mathbf{Q}}^\ell(f)-q}=0$.
\end{proof}

The next result tells us that, indeed, this $q$ is what we are after.
\begin{lemma}\label{lemma:fixed_point_nice_sets_is_upper_hit_time}
Let $\overline{h}_A^{(n)}$ and $\overline{h}_A^*$ be as in~\propositionref{prop:limit_is_lower_exp}, and let $q$ be as in~\propositionref{prop:nice_set_fixed_point}. Then $\overline{h}_A^*(x)=q(x)$ for all $x\in\mathcal{Z}\cap A^c$.
\end{lemma}
\begin{proof}
First we will show that, for any $n\in\natswith$, it holds for all $x\in\mathcal{Z}\cap A^c$ that
\begin{equation*}
\left[\overline{\mathbf{Q}}\left( \overline{h}_A^{(n)}\big\vert_{\mathcal{Z}\cap A^c} \right)\right](x) = \overline{h}_A^{(n+1)}(x)\,.
\end{equation*}
First note that, for all $x\in\mathcal{Z}$, it holds that
\begin{equation*}
\left(\overline{h}_A^{(n)}\big\vert_{\mathcal{Z}\cap A^c}\right)^\uparrow(x) = \overline{h}_A^{(n)}(x)\,,
\end{equation*}
since $\overline{h}_A^{(n)}(x)=0$ for $x\in A\subseteq \mathcal{Z}$.

Hence we find, for all $x\in\mathcal{Z}\cap A^c$, and using~\corollaryref{cor:nice_sets_functions_depend_on_nice}, that
\begin{align*}
\left[\overline{\mathbf{Q}}\left( \overline{h}_A^{(n)}\big\vert_{\mathcal{Z}\cap A^c} \right)\right](x) &= \left[1 + \overline{Q}\left(\overline{h}_A^{(n)}\big\vert_{\mathcal{Z}\cap A^c}\right)\right](x) \\
 &= \left[1 + \overline{T}\left(\left(\overline{h}_A^{(n)}\big\vert_{\mathcal{Z}\cap A^c}\right)^\uparrow\right)\right](x) \\
 &= \left[1 + \overline{T}\,\overline{h}_A^{(n)}\right](x) \\
 &= \mathbb{I}_{A^c}(x) + \mathbb{I}_{A^c}(x)\bigl[\overline{T}\overline{h}_A^{(n)}\bigr](x) \\
 &= \overline{h}_A^{(n+1)}\,.
\end{align*}
It follows that, for all $x\in\mathcal{Z}\cap A^c$ and all $n\in\nats$,
\begin{equation*}
\left[\overline{\mathbf{Q}}^n\bigl(\overline{h}_A^{(0)}\big\vert_{\mathcal{Z}\cap A^c}\bigr)\right](x) = \overline{h}_A^{(n)}(x)\,.
\end{equation*}
Thus we get, using~\propositionref{prop:limit_is_lower_exp} and ~\propositionref{prop:nice_set_fixed_point}, for all $x\in\mathcal{Z}\cap A^c$, that
\begin{equation*}
q(x) = \lim_{n\to+\infty} \left[\overline{\mathbf{Q}}^n\bigl(\overline{h}_A^{(0)}\big\vert_{\mathcal{Z}\cap A^c}\bigr)\right](x) = \lim_{n\to+\infty} \overline{h}_A^{(n)}(x) = \overline{h}_A^*(x)\,.
\end{equation*}
\end{proof}

The next step will involve finding a precise, homogeneous Markov chain $P\in\imchom$ whose expected hitting time agrees with $q$ when starting in $\mathcal{Z}\cap A^c$. To this end, we will first construct the particular transition matrix $W\in\mathcal{T}$ that characterises this $P$. 

We will now construct the transition matrix $W$ of interest, in a row-by-row manner. First, for all $z\in\mathcal{Z}$, due to~\eqref{eq:lower_trans_reached} there is some $T_z\in\mathcal{T}$ such that
\begin{equation}\label{eq:upper_time_reaching_mat_nice}
\bigl[T_z q^\uparrow\bigr](z) = \bigl[\overline{T} q^\uparrow\bigr](z)\,,
\end{equation}
and we let the $z$-th row of $W$ be defined as
\begin{equation*}
W(z,\cdot) \coloneqq T_z(z,\cdot)\,.
\end{equation*}
Next, for all $b\in\mathcal{B}$, due to~\eqref{eq:lower_trans_reached}, there is some $T_b\in\mathcal{T}$ such that
\begin{equation}\label{eq:upper_time_reaching_mat_bad}
\bigl[T_b\mathbb{I}_\mathcal{B}\bigr](b) = \bigl[\overline{T}\mathbb{I}_\mathcal{B}\bigr](b)\,,
\end{equation}
and we let the $b$-th row of $W$ be defined as
\begin{equation*}
W(b,\cdot) \coloneqq T_b(b,\cdot)\,.
\end{equation*}
Finally, for all $u\in\mathcal{U}$, by definition, there is some $n\in\nats$ such that $\bigl[\overline{T}^n\mathbb{I}_\mathcal{B}\bigr](u)>0$. Let $n_u\in\nats$ be the smallest number for which this holds, for this $u$. Then, due to~\eqref{eq:lower_trans_reached}, there is some $T_u\in\mathcal{T}$ such that
\begin{equation}\label{eq:upper_time_reaching_mat_unfortunate}
\bigl[T_u \overline{T}^{n_u-1}\mathbb{I}_\mathcal{B}\bigr](u) = \bigl[\overline{T}^{n_u}\mathbb{I}_\mathcal{B}\bigr](u)\,,
\end{equation}
and we define the $u$-th row of $W$ as
\begin{equation*}
W(u,\cdot) \coloneqq T_u(u,\cdot)\,.
\end{equation*}
The first thing to notice is that $W\in\mathcal{T}$, because $\mathcal{T}$ has separately specified rows and $W$ was constructed from the rows of elements of $\mathcal{T}$. 

We need the following property:
\begin{lemma}\label{lemma:nice_sets_complement_matrix_system_solution}
Let $W$ be the transition matrix constructed above, and let $S$ be the $\lvert \mathcal{Z}\cap A^c\rvert\times \lvert \mathcal{Z}\cap A^c\rvert$ matrix that is defined, for all $x,y\in\mathcal{Z}\cap A^c$, as $S(x,y)\coloneqq W(x,y)$. Let $I$ denote the identity matrix. Then $(I-S)^{-1}$ exists.
\end{lemma}
\begin{proof}
The proof basically follows the structure of~\cite[Theorem 11.4]{grinstead2012introduction}, but it will take a bit of effort to get all the moving pieces in place.

The first thing to note is that, because $W\in\mathcal{T}$, and due to~\lemmaref{lemma:non_bad_or_unfortunate_states_are_closed}, we have for all $x\in\mathcal{Z}$ that
\begin{equation*}
1 = \bigl[\underline{T}\mathbb{I}_\mathcal{Z}\bigr](x) \leq \bigl[W\mathbb{I}_\mathcal{Z}\bigr](x)\,,
\end{equation*}
from which it follows that $W(x,y)=0$ for all $y\notin \mathcal{Z}$.

Next, we note that because $\transmatset$ is assumed to be $A$-inert and because $W\in\mathcal{T}$, it holds for any $f\in\gambles$ that for all $x\in A$ we have $\bigl[Wf\bigr](x) = f(x)$. This immediately implies that, moreover, for any $n\in\nats$ it holds that
\begin{equation}\label{lemma:nice_sets_complement_matrix_system_solution:repeatability}
\bigl[W^nf\bigr](x) = f(x)\,,
\end{equation}
for all $x\in A$.

We will next show that, for any $f\in\gambles$ such that $f(x)=0$ for all $x\in A$ and any $n\in\natswith$, it holds for all $x\in \mathcal{Z}\cap A^c$ that
\begin{equation}\label{lemma:nice_sets_complement_matrix_system_solution:block_matrix_iterate}
\Bigl[S^n \bigl(f\big\vert_{\mathcal{Z}\cap A^c}\bigr)\Bigr](x) = \bigl[W^n f\bigr](x)\,.
\end{equation}
We will prove this by induction. The induction base is trivial; for $n=0$ we have for all $x\in \mathcal{Z}\cap A^c$ that
\begin{equation*}
\Bigl[S^n \bigl(f\big\vert_{\mathcal{Z}\cap A^c}\bigr)\Bigr](x) = \bigl[f\big\vert_{\mathcal{Z}\cap A^c}\big](x) = f(x) = \bigl[W^n f\bigr](x)\,.
\end{equation*}
So now assume that the statement is true for some $n\in\natswith$; we will show that it is also true for $n+1$. To this end, we note that
\begin{align*}
\Bigl[S^{n+1} \bigl(f\big\vert_{\mathcal{Z}\cap A^c}\bigr)\Bigr](x) &= \sum_{y\in\mathcal{Z}\cap A^c} S(x,y)\Bigl[S^{n} \bigl(f\big\vert_{\mathcal{Z}\cap A^c}\bigr)\Bigr](y) \\
 &= \sum_{y\in\mathcal{Z}\cap A^c} S(x,y)\bigl[W^n f\bigr](y) \\
&= \sum_{y\in\mathcal{Z}\cap A^c} W(x,y)\bigl[W^n f\bigr](y) \\
&= \sum_{y\in\mathcal{Z}} W(x,y)\bigl[W^n f\bigr](y) \\
&= \sum_{y\in\states} W(x,y)\bigl[W^n f\bigr](y) \\
&= \bigl[W^{n+1}f\bigr](x)\,,
\end{align*}
where the second equality used the induction hypothesis; the third equality used the definition of $S$; the fourth equality used the assumption that $f(y)=0$ for $y\in A$ and, therefore, that $\bigl[W^n f\bigr](y)=0$ for $y\in A$ due to~\eqref{lemma:nice_sets_complement_matrix_system_solution:repeatability}; and the fifth equality used that $W(x,y)=0$ for $y\notin \mathcal{Z}$, as established above. Hence we conclude that~\eqref{lemma:nice_sets_complement_matrix_system_solution:block_matrix_iterate} indeed holds.

Next, because $W\in\mathcal{T}$, it holds for any $f\in \gambles$ that $Wf \leq \overline{T}f$. It therefore follows from a straightforward induction argument based on monotonicity that, for all $n\in\nats$ and all $f\in\gambles$, it holds that
\begin{equation}\label{lemma:nice_sets_complement_matrix_system_solution:W_dominated}
W^nf \leq \overline{T}^nf\,.
\end{equation}
Now, by~\corollaryref{cor:nice_set_leaves_Ac}, for all $x\in\mathcal{Z}\cap A^c$ there is some $n_x\in\nats$ such that $\bigl[\overline{T}^{n_x}\mathbb{I}_{A^c}\bigr](x)<1$. Let $n\coloneqq \max_{x\in\mathcal{Z}\cap A^c} n_x$. Then for all $x\in\mathcal{Z}\cap A^c$, it holds that
\begin{equation*}
\bigl[\overline{T}^n\mathbb{I}_{A^c}\bigr](x) < 1\,,
\end{equation*}
again due to~\corollaryref{cor:nice_set_leaves_Ac}. 

Now due to~\eqref{lemma:nice_sets_complement_matrix_system_solution:W_dominated} and~\lemmaref{lemma:leave_nice_set_complement_stays_nice}, for all $x\in\mathcal{Z}\cap A^c$, it holds that
\begin{equation*}
\bigl[W^n\mathbb{I}_{\mathcal{Z}\cap A^c}\bigr](x) \leq \bigl[\overline{T}^n\mathbb{I}_{\mathcal{Z}\cap A^c}\bigr](x) = \bigl[\overline{T}^n\mathbb{I}_{A^c}\bigr](x) < 1\,.
\end{equation*}
Hence, using~\eqref{lemma:nice_sets_complement_matrix_system_solution:block_matrix_iterate} and the fact that $\mathbb{I}_{\mathcal{Z}\cap A^c}\big\vert_{\mathcal{Z}\cap A^c}=1$, we have for all $x\in\mathcal{Z}\cap A^c$ that
\begin{equation*}
\bigl[S^n1\bigr](x) = \bigl[W^n\mathbb{I}_{\mathcal{Z}\cap A^c}\bigr](x) <1\,.
\end{equation*}
Let $p\in\reals$ be defined as
\begin{equation*}
p \coloneqq \max_{x\in\mathcal{Z}\cap A^c}\bigl[S^n1\bigr](x)\,.
\end{equation*}
Then $p<1$ because $\mathcal{Z}\cap A^c$ is finite. Moreover, because $W\in\transmatset$ is a transition matrix, it holds that $S(x,y)=W(x,y)\geq 0$ for all $x,y\in\mathcal{Z}\cap A^c$. Clearly this implies that also $S^n(x,y)\geq 0$ for all $x,y\in\mathcal{Z}\cap A^c$, and therefore $p\geq 0$. 

To be explicit in what follows, let $\mathbf{p}$ be the vector in $\mathcal{L}(\mathcal{Z}\cap A^c)$ with constant value $p$. Then, for all $x\in \mathcal{Z}\cap A^c$ and for any $m\in\nats$,
\begin{align*}
\bigl[S^{mn}1\bigr](x) &= \bigl[S^{(m-1)n}S^n1\bigr](x) \\
 &= \bigl[S^{(m-1)n}(S^n1)\bigr](x) \\
 &\leq \bigl[S^{(m-1)n}\mathbf{p}\bigr](x) \\
 &= p\bigl[ S^{(m-1)n}1\bigr](x) \leq \cdots \leq p^m\,,
\end{align*}
where the inequalities use that fact that both $p$ and all the entries of $S^{mn}$ (and of $S^{(m-1)n}$ and so forth) are non-negative.

It follows that because $0\leq p<1$,
\begin{equation*}
\lim_{m\to+\infty} S^{mn} = 0\,,
\end{equation*}
where the right hand side denotes the zero matrix. 

Now, suppose $f\in\mathcal{L}(\mathcal{Z}\cap A^c)$ is such that
\begin{equation*}
(I-S)f = 0\,,
\end{equation*}
or in other words, that $f=Sf$. Iterating this, we find that, with $n$ as above, it holds for any $m\in\nats$ that
\begin{equation*}
f = S^{mn}f\,.
\end{equation*}
Because $\lim_{m\to+\infty}S^{mn} = 0$, this leads us to conclude that
\begin{equation*}
f = \lim_{m\to+\infty} S^{mn} f = 0\,.
\end{equation*}
Hence the kernel of $I-S$ is trivial, 
\begin{equation*}
\mathrm{ker}(I-S)=\{0\}\,,
\end{equation*}
whence $(I-S)^{-1}$ exists.
\end{proof}

The fact that $W\in\transmatset$ implies that we can find a precise homogeneous Markov chain $P\in\imchom$ corresponding to $W$ (there may in fact be multiple such $P$ with different initial models, but this is irrelevant). Next we characterise the expected hitting time of any such $P$.
\begin{lemma}\label{lemma:upper_hit_time_reaching_process_characterisation}
Let $W\in\transmatset$ be the transition matrix constructed above, and let $P\in\imchom$ be any homogeneous Markov chain with transition matrix $W$. Then, for all $x\in\states$, it holds that
\begin{equation*}
h_A^P(x) = \left\{\begin{array}{ll}
0 & \text{if $x\in A\subseteq \mathcal{Z}$, and} \\
q(x) & \text{if $x\in\mathcal{Z}\cap A^c$, and} \\
+\infty & \text{if $x\in\mathcal{B}$, and} \\
+\infty & \text{if $x\in\mathcal{U}$\,.}
\end{array}\right.
\end{equation*} 
\end{lemma}
\begin{proof}
By~\lemmaref{lemma:homogen_hitting_time_is_minimal_system_solution}, the expected hitting time $h_A^P$ of $P$ is the minimal non-negative solution to
\begin{equation}\label{lemma:upper_hit_time_reaching_process_characterisation:eq:system}
h_A^P = \mathbb{I}_{A^c} + \mathbb{I}_{A^c}\cdot Wh_A^P\,.
\end{equation}
We now prove the claimed values of $h_A^P$ in the order that they are stated. First, for any $x\in A\subseteq \mathcal{Z}$, using~\eqref{lemma:upper_hit_time_reaching_process_characterisation:eq:system} we find
\begin{equation*}
h_A^P(x) = \mathbb{I}_{A^c}(x) + \mathbb{I}_{A^c}(x)\bigl[ Wh_A^P\bigr](x) = 0\,,
\end{equation*}
due to the presence of the indicator functions $\mathbb{I}_{A^c}$.

We now move on to the values on $\mathcal{Z}\cap A^c$, and start by establishing that these are real-valued, which we will need further on. To this end, we first infer from~\propositionref{prop:vovk_imprecise_dominates_compatible_precise_vovk} and~\lemmaref{lemma:precise_vovk_hit_time_equal_precise_measure_hit_time} that,
\begin{align*}
\uexpvovk\bigl[H_A\,\vert\,X_0\bigr] &\geq \sup_{Q\in\imcirr} \overline{\mathbb{E}}_Q^\mathrm{V}\bigl[H_A\,\vert\,X_0\bigr] \\
&= \sup_{Q\in\imcirr} \mathbb{E}_Q\bigl[H_A\,\vert\,X_0\bigr] \\
&= \uexpirr\bigl[H_A\,\vert\,X_0\bigr] \\
& \geq \uexphom\bigl[H_A\,\vert\,X_0\bigr] \geq \expprec\bigl[H_A\,\vert\,X_0\bigr]\,,
\end{align*}
where the last two inequalities hold since $P\in\imchom\subseteq\imcirr$. Moreover, by~\lemmaref{lemma:fixed_point_nice_sets_is_upper_hit_time}, the game-theoretic upper expected hitting time agrees with $q$ on $\mathcal{Z}\cap A^c$. This implies that
\begin{equation*}
h_A^P(x) \leq q(x)\quad\text{for all $x\in \mathcal{Z}\cap A^c$}\,,
\end{equation*}
and, because $q$ is real-valued due to~\lemmaref{prop:nice_set_fixed_point}, this implies in particular that $h_A^P(x)$ is also real-valued for $x\in\mathcal{Z}\cap A^c$.

Now we note that, for all $z\in \mathcal{Z}\cap A^c$, the indicators evaluate to ones; so we can rewrite
\begin{equation*}
h_A^P(z) = \mathbb{I}_{A^c}(z) + \mathbb{I}_{A^c}(z)\bigl[ Wh_A^P\bigr](z) = 1 + \bigl[Wh_A^P\bigr](z)\,.
\end{equation*}
The next thing to note is that, because $W\in\mathcal{T}$, due to~\lemmaref{lemma:non_bad_or_unfortunate_states_are_closed}, for all $z\in\mathcal{Z}\cap A^c$ it holds that
\begin{equation*}
1 = \bigl[\underline{T}\mathbb{I}_\mathcal{Z}\bigr](z) \leq \bigl[W\mathbb{I}_\mathcal{Z}\bigr](z)\,,
\end{equation*}
from which we conclude that $W(z,x)=0$ for all $x\notin \mathcal{Z}$. Thus, we get
\begin{equation*}
h_A^P(z) = 1 + \sum_{x\in \mathcal{Z}} W(z,x) h_A^P(x)\,,
\end{equation*}
for all $z\in \mathcal{Z}\cap A^c$. However, as we have already shown, it holds that $h_A^P(x)=0$ for $x\in A\subseteq\mathcal{Z}$. So in fact we have
\begin{equation}\label{lemma:upper_hit_time_reaching_process_characterisation:eq:nice_system_part}
h_A^P(z) = 1 + \sum_{x\in\mathcal{Z}\cap A^c} W(z,x)h_A^P(x)\,,
\end{equation}
for all $z\in\mathcal{Z}\cap A^c$. Now let $S$ be the $\lvert \mathcal{Z}\cap A^c\rvert\times \lvert \mathcal{Z}\cap A^c\rvert$ matrix that is defined, for all $x,y \in\mathcal{Z}\cap A^c$, as
\begin{equation*}
S(x,y) \coloneqq W(x,y)\,.
\end{equation*}
Then we can rewrite the system~\eqref{lemma:upper_hit_time_reaching_process_characterisation:eq:nice_system_part} as
\begin{equation*}
h_A^P\big\vert_{\mathcal{Z}\cap A^c} = 1 + S h_A^P\big\vert_{\mathcal{Z}\cap A^c}\,,
\end{equation*}
and, using the previously established fact that $h_A^P\big\vert_{\mathcal{Z}\cap A^c}$ is real-valued, we rewrite this as
\begin{equation*}
(I - S) h_A^P\big\vert_{\mathcal{Z}\cap A^c} = 1\,.
\end{equation*}
Moreover, we also recognise $S$ as the same matrix constructed in~\lemmaref{lemma:nice_sets_complement_matrix_system_solution}, so we know that $(I-S)^{-1}$ exists. Hence in particular, we have that
\begin{equation*}
h_A^P\big\vert_{\mathcal{Z}\cap A^c} = (I-S)^{-1}1\,,
\end{equation*}
where the $1$ is still the vector in $\mathcal{L}(\mathcal{Z}\cap A^c)$ with constant value one. It remains to relate this to the quantity $q$.

We recall from~\sectionref{subsec:expected_hit_times} the fixed-point property of $\overline{h}_A^*$:
\begin{equation*}
\overline{h}_A^* = \mathbb{I}_{A^c} + \mathbb{I}_{A^c}\cdot \overline{T}\,\overline{h}_A^*\,.
\end{equation*}
Because $\overline{h}_A^*(x)=0=q^\uparrow(x)$ for all $x\in A\subseteq \mathcal{Z}$ and because, by~\lemmaref{lemma:fixed_point_nice_sets_is_upper_hit_time}, $\overline{h}_A^*(x)=q(x)$ for all $x\in\mathcal{Z}\cap A^c$, it holds that $\overline{h}_A^*(x) = q^\uparrow(x)$ for all $x\in\mathcal{Z}$. Hence by~\corollaryref{cor:nice_sets_functions_depend_on_nice} it holds for all $x\in\mathcal{Z}$ that
\begin{equation*}
q^\uparrow(x) = \mathbb{I}_{A^c}(x) + \mathbb{I}_{A^c}(x)\bigl[\overline{T}q^\uparrow\bigr](x)\,,
\end{equation*}
and, therefore in particular for all $x\in\mathcal{Z}\cap A^c$, that
\begin{equation*}
q^\uparrow(x) = 1 + \bigl[\overline{T}q^\uparrow\bigr](x)\,.
\end{equation*}
Using the selection for the construction of $W$ in~\eqref{eq:upper_time_reaching_mat_nice}, this implies for all $x\in\mathcal{Z}\cap A^c$ that
\begin{equation*}
q^\uparrow(x) = 1 + \bigl[Wq^\uparrow\bigr](x)\,.
\end{equation*}
Because $q^\uparrow(x)=0$ for all $x\notin (\mathcal{Z}\cap A^c)$, we can rewrite this to
\begin{equation*}
q(x) = 1 + \sum_{y\in\mathcal{Z}\cap A^c} W(x,y) q(y)\,,
\end{equation*}
for all $x\in\mathcal{Z}\cap A^c$, so rewriting in matrix-vector notation, and recognising again the definition of the matrix $S$,
\begin{equation*}
q = 1 + Sq\,.
\end{equation*}
We conclude that
\begin{equation*}
q = (I-S)^{-1}1 = h_A^P\big\vert_{\mathcal{Z}\cap A^c}\,,
\end{equation*}
where the last step used the earlier derivation above.

We next consider the values on $\mathcal{B}$. We start by noting that, due to~\lemmaref{lemma:bad_states_are_closed} and the selection used for the construction of $W$ in~\eqref{eq:upper_time_reaching_mat_bad}, it holds for all $b\in\mathcal{B}$ that
\begin{equation*}
\bigl[W\mathbb{I}_\mathcal{B}\bigr](b) = \bigl[\overline{T}\mathbb{I}_\mathcal{B}\bigr](b) = 1\,,
\end{equation*}
which implies that $W(b,x)=0$ for all $x\notin \mathcal{B}$.

We again rewrite the system~\eqref{lemma:upper_hit_time_reaching_process_characterisation:eq:system}, and note that since $\mathcal{B}\subseteq A^c$, for all $b\in\mathcal{B}$ it holds that
\begin{equation*}
h_A^P(b) = 1 + \bigl[Wh_A^P\bigr](b)\,.
\end{equation*}
Now let $b \in \arg\min\{h_A^P(x)\,\vert\,x\in\mathcal{B}\}$, and suppose \emph{ex absurdo} that $C\coloneqq h_A^P(b) < +\infty$. Then,
\begin{align*}
C = h_A^P(b) &= 1 + \bigl[Wh_A^P\bigr](b) \\
 &= 1 + \sum_{x\in \states} W(b,x)h_A^P(x) \\
 &= 1 + \sum_{x\in \mathcal{B}} W(b,x)h_A^P(x) \\
 &\geq 1 + \sum_{x\in \mathcal{B}} W(b,x)C \\
 &= 1 + C\,,
\end{align*}
which, since $0\leq C<+\infty$, leads us to conclude that $0\geq 1$, a contradiction. Hence we have in fact that $h_A^P(x) = +\infty$ for all $x\in\mathcal{B}$.

We finally consider the values on $\mathcal{U}$. We again start by noting that $\mathcal{U}\subseteq A^c$ so we can rewrite the system~\eqref{lemma:upper_hit_time_reaching_process_characterisation:eq:system} such that, for all $u\in\mathcal{U}$,
\begin{equation*}
h_A^P(u) = 1 + \bigl[W h_A^P\bigr](u)\,.
\end{equation*}
To avoid trivialities, we can assume that $\mathcal{U}\neq \emptyset$. The remainder of the proof is then by backwards induction. Due to~\corollaryref{cor:unfortunate_set_induction_base}, there is some $u\in\mathcal{U}$ with $n_u=1$ such that $\bigl[\overline{T}^{n_u}\mathbb{I}_\mathcal{B}\bigr](u)>0$. Then due to the selection used for the construction of $W$ in~\eqref{eq:upper_time_reaching_mat_unfortunate}, it holds that
\begin{equation*}
\bigl[W\mathbb{I}_\mathcal{B}\bigr](u) = \bigl[\overline{T}\mathbb{I}_\mathcal{B}\bigr](u) > 0\,.
\end{equation*}
This implies that there is some $b\in\mathcal{B}$ such that $W(u,b)>0$. We have already established above that $h_A^P(x)=+\infty$ for all $x\in\mathcal{B}$. Using also the non-negativity of $h_A^P$, we get
\begin{equation*}
h_A^P(u) = 1 + \bigl[W h_A^P\bigr](u) \geq 1 + W(u,b)h_A^P(b) = +\infty\,.
\end{equation*}
This provides the induction base. 

Now suppose that $h_A^P(u)=+\infty$ for all $u\in\mathcal{U}$ with $n_u\leq n$, for some $n\in\nats$. Now consider any $v\in\mathcal{U}$ with $n_v = n+1$. Due to the selection used for the construction of $W$ in~\eqref{eq:upper_time_reaching_mat_unfortunate}, it holds that
\begin{equation*}
\bigl[W\overline{T}^{n}\mathbb{I}_\mathcal{B}\bigr](v) = \bigl[\overline{T}^{n_v}\mathbb{I}_\mathcal{B}\bigr](v) > 0\,.
\end{equation*}
This implies that, for some $u\in\states$,
\begin{equation*}
W(v,u)\bigl[\overline{T}^n\mathbb{I}_\mathcal{B}\bigr](u) > 0\,.
\end{equation*}
It follows that $W(v,u)>0$ and, as we will show next, that $u\in\mathcal{U}$ and $n_u \leq n$. To see that $u\in\mathcal{U}$, first suppose \emph{ex absurdo} that $u\in A$. By the above inequality, we have $\bigl[\overline{T}^n\mathbb{I}_\mathcal{B}\bigr](u)>0$. However, because $u\in A$ and because $\transmatset$ is $A$-inert, it follows from~\lemmaref{lemma:a_inert_fixed_value} that $\bigl[\overline{T}^n\mathbb{I}_\mathcal{B}\bigr](u) = \mathbb{I}_\mathcal{B}(u) = 0$, because $A\cap \mathcal{B}=\emptyset$. This is a contradiction, and hence $u\notin A$.

Next, suppose \emph{ex absurdo} that $u\in \mathcal{B}$. Using the fact that $W(v,u)>0$ by the above inequality, we find that
\begin{equation*}
0 < W(v,u) = \bigl[W\mathbb{I}_{\{u\}}\bigr](v) \leq \bigl[\overline{T}\mathbb{I}_{\{u\}}\bigr](v) \leq \bigl[\overline{T}\mathbb{I}_{\mathcal{B}}\bigr](v)\,,
\end{equation*}
where the second inequality uses the fact that $W\in\transmatset$ and the third inequality uses the monotonicity of $\overline{T}$ and the assumption $u\in\mathcal{B}$. This implies that $n_v=1$, a contradiction due to the choice of $v$. Hence $u\notin\mathcal{B}$. It now follows directly from the definition of $\mathcal{U}$ that $u\in\mathcal{U}$, because $\bigl[\overline{T}^n\mathbb{I}_\mathcal{B}\bigr](u)>0$. That $n_u\leq n$ is then immediate.

 Hence, by the induction hypothesis, $h_A^P(u)=+\infty$. Moreover, because $W(v,u)>0$, we find that
\begin{equation*}
h_A^P(v) = 1 + \bigl[W h_A^P\bigr](v) \geq W(v,u)h_A^P(u) = +\infty\,.
\end{equation*}
\end{proof}

We can now derive the conclusion of interest, under the assumption that the model is $A$-inert.
\begin{lemma}\label{lemma:upper_hit_time_reached_for_inert}
Suppose that $\transmatset$ is $A$-inert. Then there is some $P\in\imchom$ such that
\begin{equation*}
\uexpvovk\bigl[H_A\,\vert\,X_0\bigr] = \uexpirr\bigl[H_A\,\vert\,X_0\bigr] = \uexphom\bigl[H_A\,\vert\,X_0\bigr] = \expprec\bigl[H_A\,\vert\,X_0\bigr]\,.
\end{equation*}
\end{lemma}
\begin{proof}
Let $P\in\imchom$ be any precise homogeneous Markov chain corresponding to the transition matrix $W\in\mathcal{T}$ constructed above; this choice is only unique up to the initial model $P(X_0)$, but that initial model is irrelevant in what follows. 

Using~\propositionref{prop:vovk_imprecise_dominates_compatible_precise_vovk} and~\lemmaref{lemma:precise_vovk_hit_time_equal_precise_measure_hit_time}, it holds that
\begin{align*}
\uexpvovk\bigl[H_A\,\vert\,X_0\bigr] &\geq \sup_{Q\in\imcirr} \overline{\mathbb{E}}_Q^\mathrm{V}\bigl[H_A\,\vert\,X_0\bigr] \\
 &= \sup_{Q\in\imcirr} \mathbb{E}_Q\bigl[H_A\,\vert\,X_0\bigr] \\
 &= \uexpirr\bigl[H_A\,\vert\,X_0\bigr] \geq \uexphom\bigl[H_A\,\vert\,X_0\bigr] \geq \expprec\bigl[H_A\,\vert\,X_0\bigr]\,,
\end{align*}
where we used $P\in\imchom\subseteq \imcirr$.

Now, for all $x\in\mathcal{Z}\cap A^c$ it holds due to~\propositionref{prop:limit_is_lower_exp}, ~\lemmaref{lemma:fixed_point_nice_sets_is_upper_hit_time} and~\lemmaref{lemma:upper_hit_time_reaching_process_characterisation} that
\begin{equation*}
\uexpvovk\bigl[H_A\,\vert\,X_0=x\bigr] = \overline{h}_A^*(x) = q(x) = h_A^P(x) = \expprec\bigl[H_A\,\vert\,X_0=x\bigr]\,.
\end{equation*}
Next, for all $x\in A$ it holds that $\expprec\bigl[H_A\,\vert\,X_0=x\bigr]=0$ due to~\lemmaref{lemma:upper_hit_time_reaching_process_characterisation}. Therefore, and because $\uexpvovk\bigl[H_A\,\vert\,X_0=x\bigr] = \overline{h}_A^*(x)$ and $\overline{h}_A^*$ satisfies
\begin{equation*}
\overline{h}_A^* = \mathbb{I}_{A^c} + \mathbb{I}_{A^c}\hprod \overline{T}\,\overline{h}_A^*\,,
\end{equation*}
for $x\in A$ we get using $\mathbb{I}_{A^c}(x)=0$ that
\begin{equation*}
\uexpvovk\bigl[H_A\,\vert\,X_0=x\bigr] = \overline{h}_A^*(x) = 0 = \expprec\bigl[H_A\,\vert\,X_0=x\bigr]\,.
\end{equation*}
Finally, by~\lemmaref{lemma:upper_hit_time_reaching_process_characterisation}, we have for all $x\in\mathcal{B}\cup\mathcal{U}$ that
\begin{equation*}
 \expprec\bigl[H_A\,\vert\,X_0=x\bigr] = h_A^P(x) = +\infty\,.
\end{equation*}
Because by the above it holds that $\uexpvovk\bigl[H_A\,\vert\,X_0\bigr] \geq \expprec\bigl[H_A\,\vert\,X_0\bigr]$, this implies that
\begin{equation*}
\uexpvovk\bigl[H_A\,\vert\,X_0=x\bigr] = +\infty = \expprec\bigl[H_A\,\vert\,X_0=x\bigr]\,.
\end{equation*}
We therefore conclude that in fact
\begin{equation*}
\uexpvovk\bigl[H_A\,\vert\,X_0\bigr] = \expprec\bigl[H_A\,\vert\,X_0\bigr]\,.
\end{equation*}
A final application of our earlier inequality now gives
\begin{equation*}
\uexpvovk\bigl[H_A\,\vert\,X_0\bigr] = \uexpirr\bigl[H_A\,\vert\,X_0\bigr] = \uexphom\bigl[H_A\,\vert\,X_0\bigr] = \expprec\bigl[H_A\,\vert\,X_0\bigr]\,.
\end{equation*}
\end{proof}

In order to prove the remaining statement in~\theoremref{thm:lower_hitting_time_reach_and_equal}, it simply remains to extend the above result to general (i.e. non $A$-inert) models:
\quad\newline
\begin{proofof}{\theoremref{thm:lower_hitting_time_reach_and_equal}}
Let $\mathcal{S}$ be the $A$-inert modification of $\mathcal{T}$, and let $\overline{\mathbb{E}}_\mathcal{S}^\mathrm{V}[H_A\,\vert\,X_0]$ be the upper expected hitting time of the $A$-inert modification of the game-theoretic imprecise Markov chain corresponding to $\mathcal{T}$. Due to~\lemmaref{lemma:hit_time_game_is_a_inert} it holds that
\begin{equation*}
\overline{\mathbb{E}}_\mathcal{S}^\mathrm{V}[H_A\,\vert\,X_0] = \overline{\mathbb{E}}_\mathcal{T}^\mathrm{V}[H_A\,\vert\,X_0]\,.
\end{equation*}

Due to~\lemmaref{lemma:upper_hit_time_reached_for_inert}, we can find some $Q\in\mathcal{P}_\mathcal{S}^\mathrm{H}$ such that
\begin{equation*}
\overline{\mathbb{E}}_\mathcal{S}^\mathrm{V}[H_A\,\vert\,X_0] = \mathbb{E}_{Q}[H_A\,\vert\,X_0\bigr]\,.
\end{equation*}
Because $Q\in\mathcal{P}_\mathcal{S}^\mathrm{H}$ and $\mathcal{S}$ is the $A$-inert modification of $\mathcal{T}$, there is some $P\in\imchom$ such that $Q$ is the $A$-inert modification of $P$. Due to~\lemmaref{lemma:hit_time_homogen_is_a_inert} it then holds that
\begin{equation*}
\mathbb{E}_{Q}[H_A\,\vert\,X_0\bigr] = \mathbb{E}_{P}[H_A\,\vert\,X_0\bigr]\,,
\end{equation*}
from which we conclude that
\begin{equation*}
\overline{\mathbb{E}}_\mathcal{T}^\mathrm{V}[H_A\,\vert\,X_0] = \overline{\mathbb{E}}_\mathcal{S}^\mathrm{V}[H_A\,\vert\,X_0] = \mathbb{E}_{Q}[H_A\,\vert\,X_0\bigr] = \mathbb{E}_{P}[H_A\,\vert\,X_0\bigr]\,.
\end{equation*}
Finally, using~\propositionref{prop:vovk_imprecise_dominates_compatible_precise_vovk} and~\lemmaref{lemma:precise_vovk_hit_time_equal_precise_measure_hit_time}, it holds that
\begin{align*}
\uexpvovk\bigl[H_A\,\vert\,X_0\bigr] &\geq \sup_{Q\in\imcirr} \overline{\mathbb{E}}_Q^\mathrm{V}\bigl[H_A\,\vert\,X_0\bigr] \\
 &= \sup_{Q\in\imcirr} \mathbb{E}_Q\bigl[H_A\,\vert\,X_0\bigr] \\
 &= \uexpirr\bigl[H_A\,\vert\,X_0\bigr] \geq \uexphom\bigl[H_A\,\vert\,X_0\bigr] \geq \expprec\bigl[H_A\,\vert\,X_0\bigr]\,,
\end{align*}
where we used $P\in\imchom\subseteq \imcirr$.
\end{proofof}

\begin{proofof}{\corollaryref{cor:imprecise_hitting_time_is_minimal_system_solution}}
That the upper expected hitting time $\overline{h}_A$ indeed satisfies this non-linear system is immediate from the fixed-point property of $\overline{h}_A^*$ (see~\sectionref{subsec:expected_hit_times}) and~\propositionref{prop:limit_is_lower_exp}. That it is non-negative follows from the non-negativity of $H_A$. It remains to show that it is the minimal solution.

To this end, consider any non-negative $g\in\gamblesextabove$ that satisfies
\begin{equation*}
g = \mathbb{I}_{A^c} + \mathbb{I}_{A^c}\hprod \overline{T}g\,.
\end{equation*}
We will show that $g\geq \overline{h}_A$.

First note that, for all $x\in A^c$, it holds that
\begin{equation*}
g(x) = \mathbb{I}_{A^c}(x) + \mathbb{I}_{A^c}(x)\bigl[ \overline{T}g\bigr](x) = 1 + \bigl[ \overline{T}g\bigr](x) \geq 1\,,
\end{equation*}
where we used the non-negativity of $g$ and the monotonicity of $\overline{T}$. Therefore, and because $g$ is non-negative, it holds that $g\geq \mathbb{I}_{A^c}$.

Now define the map $\overline{\mathbf{H}}:\gamblesextabove\to\gamblesextabove$, for all $f\in\gamblesextabove$, as
\begin{equation*}
\overline{\mathbf{H}}(f) \coloneqq \mathbb{I}_{A^c} + \mathbb{I}_{A^c}\hprod \overline{T}f\,.
\end{equation*}
Then, clearly,
\begin{equation*}
g = \overline{\mathbf{H}}(g)\,.
\end{equation*}
Moreover, for any $n\in\nats$, let $\overline{\mathbf{H}}^n$ denote the $n$-fold composition of $\overline{\mathbf{H}}$ with itself. 
Then it is easy to see that, for all $n\in\nats$,
\begin{equation*}
\overline{h}_A^{(n)} =
\overline{\mathbf{H}}^n(\overline{h}_{A}^{(0)}) 
= \overline{\mathbf{H}}^n(\mathbb{I}_{A^c}) \,.
\end{equation*}
Finally, it follows from the monotonicity of $\overline{T}$ that $\overline{\mathbf{H}}$ is also a monotone operator and, therefore, that $\overline{\mathbf{H}}^n$ is monotone for any $n\in\nats$. Because, as we already established, $g\geq\mathbb{I}_{A^c}$, it follows that for every $n\in\nats$,
\begin{equation*}
g=\overline{\mathbf{H}}^n(g) \geq \overline{\mathbf{H}}^n(\mathbb{I}_{A^c}) = \overline{h}_A^{(n)}\,,
\end{equation*}
from which it follows that
\begin{equation*}
g \geq \lim_{n\to+\infty} \overline{h}_A^{(n)} = \overline{h}_A^* = \overline{h}_A\,.
\end{equation*}
\end{proofof}

\section{Proofs of Statements in~\sectionref{subsec:hit_probs}}\label{appendix:hit_prob_proofs}

The following lemma states a result for any set $\mathcal{S}$ of transition matrices (and its corresponding lower and upper transition operators) that satisfies the stated requirements, instead of the specific set $\mathcal{T}$ that we use in the main text; it is stated at this level of generality so that we can reuse the result later.
\begin{lemma}\label{lemma:hitting_prob_approx_sequence_nonnegative}
Let $\mathcal{S}$ be a non-empty set of transition matrices that is closed and convex and has separately specified rows, and let $\underline{S}$ and $\overline{S}$ be its corresponding lower and upper transition operators, respectively. Consider the two sequences that are defined as $\underline{p}_A^{(0)}\coloneqq \overline{p}_A^{(0)}\coloneqq \mathbb{I}_A$ and, for all $n\in\natswith$, as
\begin{equation*}
\underline{p}_A^{(n+1)} \coloneqq \mathbb{I}_A + \mathbb{I}_{A^c}\cdot \underline{S}\underline{p}_A^{(n)}\,,
\end{equation*}
and
\begin{equation*}
\overline{p}_A^{(n+1)} \coloneqq \mathbb{I}_A + \mathbb{I}_{A^c}\cdot \overline{S}\overline{p}_A^{(n)}\,.
\end{equation*}
Then $\smash{0 \leq \underline{p}_A^{(n)} \leq \overline{p}_A^{(n)} \leq 1}$ for all $n\in\natswith$.
\end{lemma}
\begin{proof}
We give a proof by induction. First we note that by definition $\smash{\underline{p}_A^{(0)}=\overline{p}_A^{(0)}=\mathbb{I}_{A}}$. Hence, because $0 \leq \mathbb{I}_{A} \leq 1$, we have established the required induction base.

Now assume that the claim is true for $n-1$, with $n\in\nats$; we will show that the claim is also true for $n$. We start by establishing the non-negativity of $\smash{\underline{p}_A^{(n)}}$, which by definition is given by
\begin{equation*}
\underline{p}_A^{(n)} = \mathbb{I}_{A} + \mathbb{I}_{A^c}\cdot \underline{S}\,\underline{p}_A^{(n-1)}\,.
\end{equation*}
By the induction hypothesis, $\underline{p}_A^{(n-1)}$ is non-negative, which implies that also $\smash{\underline{S}\,\underline{p}_A^{(n-1)}}$ is non-negative. Hence, because $\mathbb{I}_{A^c}$ is non-negative, also the product $\smash{\mathbb{I}_{A^c}\cdot \underline{S}\,\underline{p}_A^{(n-1)}}$ is non-negative. Noting that $\mathbb{I}_{A}$ is non-negative, it follows that therefore the sum $\smash{\mathbb{I}_{A} + \mathbb{I}_{A^c}\cdot \underline{S}\,\underline{p}_A^{(n-1)}}$ is non-negative, which concludes the proof that $0\leq \underline{p}_A^{(n)}$.

To show that $\smash{\underline{p}_A^{(n)}\leq \overline{p}_A^{(n)}}$, we first use the induction hypothesis and the monotonicity of $\underline{S}$ to establish that
\begin{align*}
\underline{p}_A^{(n)} &= \mathbb{I}_{A} + \mathbb{I}_{A^c}\cdot \underline{S}\,\underline{p}_A^{(n-1)} \leq \mathbb{I}_{A} + \mathbb{I}_{A^c}\cdot \underline{S}\overline{p}_A^{(n-1)}\,.
\end{align*}
Moreover, it follows from the definitions of $\underline{S}$ and $\overline{S}$ that $\underline{S}f \leq \overline{S}f$ for any $f\in\gamblesextabove$, and hence we find
\begin{align*}
\underline{p}_A^{(n)} &\leq \mathbb{I}_{A} + \mathbb{I}_{A^c}\cdot \underline{S}\overline{p}_A^{(n-1)} \leq \mathbb{I}_{A} + \mathbb{I}_{A^c}\cdot \overline{S}\overline{p}_A^{(n-1)} = \overline{p}_A^{(n)}\,,
\end{align*}
which concludes the proof that $\smash{\underline{p}_A^{(n)} \leq \overline{p}_A^{(n)}}$.

It remains to show that $\smash{\overline{p}_A^{(n)} \leq 1}$. By definition we have that
\begin{equation*}
\overline{p}_A^{(n)} = \mathbb{I}_{A} + \mathbb{I}_{A^c}\cdot\overline{S}\,\overline{p}_A^{(n-1)}\,.
\end{equation*}
The induction hypothesis tells us that $\smash{\overline{p}_A^{(n-1)} \leq 1}$, which implies that also $\smash{\overline{S}\,\overline{p}_A^{(n-1)} \leq 1}$. Now consider any $x\in\states$. Then if $x\in A$ it holds that $\mathbb{I}_A(x)=1$ and $\mathbb{I}_{A^c}(x)=0$, and hence
\begin{align*}
\overline{p}_A^{(n)}(x) &= \mathbb{I}_{A}(x) + \mathbb{I}_{A^c}(x)\cdot\Bigl[\overline{S}\,\overline{p}_A^{(n-1)}\Bigr](x) \\
 &= 1 + 0\cdot \Bigl[\overline{S}\,\overline{p}_A^{(n-1)}\Bigr](x) = 1\,.
\end{align*}
Conversely, if $x\in A^c$ then we have $\mathbb{I}_A(x)=0$ and $\mathbb{I}_{A^c}(x)=1$, and so
\begin{align*}
\overline{p}_A^{(n)}(x) &= \mathbb{I}_{A}(x) + \mathbb{I}_{A^c}(x)\cdot\Bigl[\overline{S}\,\overline{p}_A^{(n-1)}\Bigr](x) \\
 &= 0 + 1\cdot \Bigl[\overline{S}\,\overline{p}_A^{(n-1)}\Bigr](x) \\
  &=  \Bigl[\overline{S}\,\overline{p}_A^{(n-1)}\Bigr](x) \leq 1\,.
\end{align*}
So for any $x\in\states$ we have $\overline{p}_A^{(n)}(x) \leq 1$, which means that $\overline{p}_A^{(n)} \leq 1$; this concludes the proof.
\end{proof}

The previous result is directly applicable to the corresponding sequences defined in the main text.
\begin{corollary}\label{cor:hitting_prob_approx_sequence_nonnegative}
Consider the sequences $\smash{\underline{p}_A^{(n)}}$ and $\smash{\overline{p}_A^{(n)}}$ that are defined in~\sectionref{subsec:hit_probs}. Then $\smash{0 \leq \underline{p}_A^{(n)} \leq \overline{p}_A^{(n)} \leq 1}$ for all $n\in\natswith$.
\end{corollary}
\begin{proof}
Note that $\mathcal{T}$ is non-empty, closed, convex, and has separately specified rows. Now apply~\lemmaref{lemma:hitting_prob_approx_sequence_nonnegative} to the sequences $\smash{\underline{p}_A^{(n)}}$ and $\smash{\overline{p}_A^{(n)}}$.
\end{proof}

\begin{proofof}{\lemmaref{lemma:recursion_lower_prob_is_truncated_function}}
We first give the proof for $\overline{p}_A^{(n)}$; the proof for $\underline{p}_A^{(n)}$ will follow an analogous reasoning. 

This proof is by induction, and we first establish the induction base. To this end, note that $\smash{G_A^{(0)}(X_0)}$ only depends on the state $X_0$, so it follows from~[4.] in~\propositionref{prop:vovk_satisfies_coherence} that, for all $x_0\in\states$, it holds that
\begin{equation*}
\uexpvovk\bigl[G_A^{(0)}(X_0)\,\big\vert\,X_0=x_0\bigr] = G_A^{(0)}(x_0)\,.
\end{equation*}
Furthermore, it follows immediately from the definitions of $G_A^{(0)}$ and $\overline{p}_A^{(0)}$ that
\begin{equation*}
G_A^{(0)}(x_0) = \mathbb{I}_{A}(x_0) = \overline{p}_A^{(0)}(x_0)\,,
\end{equation*}
for all $x_0\in\states$, and so we get
\begin{equation*}
\uexpvovk\bigl[G_A^{(0)}(X_0)\,\big\vert\,X_0=x_0\bigr] = \overline{p}_A^{(0)}(x_0)\,,
\end{equation*}
which is the induction base that we are after.
Note that in the above, we used a small notational trick to write $G_A^{(0)}(x_0)$ for the value of $G_A^{(0)}(X_0)$ on a path $\omega$ that satisfies $\omega(0)=x_0$.

In fact, this notational trick can be extended by noting that, for any $n\in\natswith$, the value of $G_A^{(n)}$ on any $\omega\in\Omega$ is completely determined by the values $\omega(0),\ldots,\omega(n)\in\states$. Thus, we can equivalently interpret $G_A^{(n)}$ as a function on $\states^{n+1}$. Moreover, we can then write, for any $n\in\nats$ and any $x_0,\ldots,x_n\in\states$, that
\begin{equation*}
G_A^{(n)}(x_{0:n}) = \mathbb{I}_{A}(x_0) + \mathbb{I}_{A^c}(x_0)G_A^{(n-1)}(x_{1:n})\,,
\end{equation*}
with $G_A^{(0)}(x_0)=\mathbb{I}_{A}(x_0)$. Using this observation, let us now proceed with the induction step; so, we assume that the statement is true for some $n-1$, with $n\in\nats$. We will show that the statement is also true for $n$.

First, for any $x_0\in\states$ we have that
\begin{align}\label{Eq: proof of lemma:iterate_hit_prob_is_restricted_exp}
&\quad \uexpvovk\bigl[G_A^{(n)}(X_{0:n})\,\vert\,X_0=x_0\bigr]  \nonumber \\
&= \uexpvovk\bigl[\mathbb{I}_{A}(x_0) + \mathbb{I}_{A^c}(x_0)G_A^{(n-1)}(X_{1:n})\,\vert\,X_0=x_0\bigr] \nonumber \\
&= \mathbb{I}_{A}(x_0) + \mathbb{I}_{A^c}(x_0)\uexpvovk\bigl[G_A^{(n-1)}(X_{1:n})\,\vert\,X_0=x_0\bigr] \nonumber \\
&= \mathbb{I}_{A}(x_0) + \mathbb{I}_{A^c}(x_0)\uexpvovk\Bigl[\uexpvovk\bigl[G_A^{(n-1)}(X_{1:n})\,\vert\,X_{0:1}\bigr]\,\Big\vert\,X_0=x_0\Bigr],
\end{align}
where the first step used~\lemmaref{lemma:vovk condition fixes value}, the second step used [2.] and [5.] in Proposition~\ref{prop:vovk_satisfies_coherence}, and the last step used Proposition~\ref{prop:vovk_iterated_expectation}.
Now, according to Corollary~\ref{prop: vovk time shift}, $\uexpvovk\bigl[G_A^{(n-1)}(X_{1:n})\,\vert\,X_{0:1}\bigr]$ does not depend on the initial state $X_0$, so there is a function $g_{n-1} \colon \states{} \to \realsext{}$ such that
\begin{align*}
g_{n-1}(X_1) = \uexpvovk\bigl[G_A^{(n-1)}(X_{1:n})\,\vert\,X_{0:1}\bigr].
\end{align*}
Moreover, again using Corollary~\ref{prop: vovk time shift}, we also have that 
\begin{align*}
g_{n-1}(X_0) = \uexpvovk\bigl[G_A^{(n-1)}(X_{0:n-1})\,\vert\,X_{0}\bigr].
\end{align*}
By the induction hypothesis, we therefore have that $g_{n-1}=\overline{p}_A^{(n-1)}$. Plugging this back into Equation \eqref{Eq: proof of lemma:iterate_hit_prob_is_restricted_exp}, we get that
\begin{align}\label{Equation: Eq: proof of lemma:iterate_hit_prob_is_restricted_exp_2}
&\quad \uexpvovk\bigl[G_A^{(n)}(X_{0:n})\,\vert\,X_0=x_0\bigr] \nonumber \\
&= \mathbb{I}_{A}(x_0) + \mathbb{I}_{A^c}(x_0)\uexpvovk\Bigl[\overline{p}_A^{(n-1)}(X_1)\,\Big\vert\,X_0=x_0\Bigr] \nonumber \\
&= \mathbb{I}_{A}(x_0) + \mathbb{I}_{A^c}(x_0)[\overline{T} \, \overline{p}_A^{(n-1)}] (x_0) 
= \overline{p}_A^{(n)}(x_0), 
\end{align}
for all $x_0\in\states$. 
In the expression above, the second step uses Equation \eqref{eq:def_vovk_lexp} together with \cite[Proposition 15]{natan:game_theory} which we can use because $\overline{p}_A^{(n-1)}$ is non-negative (and therefore bounded below) by~\corollaryref{cor:hitting_prob_approx_sequence_nonnegative}, and the last step uses the definition of $\overline{p}_A^{(n)}$. This concludes the first part of the proof.

Analogously, we can follow the reasoning above in order to show that the statement for $\underline{p}_A^{(n)}$ holds.
Properties [2.] and [5.] in Proposition~\ref{prop:vovk_satisfies_coherence} clearly also hold for lower expectations because of conjugacy.
The only step that requires some closer attention is the second equality in Equation~\eqref{Equation: Eq: proof of lemma:iterate_hit_prob_is_restricted_exp_2} because Equation~\eqref{eq:def_vovk_lexp} and \cite[Proposition 15]{natan:game_theory} are only given for upper expectations.
However, we can easily derive this through the use of conjugacy:
\begin{align*}
\lexpvovk\Bigl[\underline{p}_A^{(n-1)}(X_1)\,\Big\vert\,X_0=x_0\Bigr]
&= - \uexpvovk\Bigl[- \underline{p}_A^{(n-1)}(X_1)\,\Big\vert\,X_0=x_0\Bigr] \\
&= - [\overline{T} \, \bigl(-\underline{p}_A^{(n-1)}\bigr)] (x_0) \\
&= [\underline{T} \, \underline{p}_A^{(n-1)}] (x_0),
\end{align*}
for all $x_0 \in \states{}$.
Again, we used Equation \eqref{eq:def_vovk_lexp} together with \cite[Proposition 15]{natan:game_theory} in the second step, which was allowed since $\underline{p}_A^{(n-1)}$ is bounded above (and therefore $-\underline{p}_A^{(n-1)}$ is bounded below) as a consequence of~\corollaryref{cor:hitting_prob_approx_sequence_nonnegative}.
\end{proofof}

In the next couple of pages we set up the required results to prove the second part of Theorem~\ref{thm:imprecise_hit_probs_equal}. As in the previous section of the appendix, we start by isolating the troublesome states. We also again use the intuitive interpretation of $\underline{T}^n\mathbb{I}_\mathcal{B}$ (resp. $\overline{T}^n\mathbb{I}_\mathcal{B}$) encoding the lower (resp. upper) probability that the process reaches an element of any $\mathcal{B}\subseteq\states$ after $n\in\nats$ steps.
\begin{definition}
Let $\mathcal{C}\subseteq A^c$ contain all states $x\in A^c$ such that $\bigl[\overline{T}^n\mathbb{I}_A\bigr](x)=0$ for all $n\in\nats$.
\end{definition}
Thus, intuitively, $\mathcal{C}$ contains exactly those states, from which reaching $A$ has upper probability zero. We need the following result, which intuitively tells us that once the process enters $\mathcal{C}$, it will remain there with lower probability one.
\begin{lemma}\label{lemma:absorbing_states_closed}
For all $T\in\mathcal{T}$, it holds that $\bigl[T\mathbb{I}_\mathcal{C}\bigr](x)=1$ for all $x\in\mathcal{C}$.
\end{lemma}
\begin{proof}
Suppose \emph{ex absurdo} that the statement is false for some $T\in\mathcal{T}$. Then there exists some $x\in\mathcal{C}$ such that
\begin{equation*}
\bigl[T\mathbb{I}_\mathcal{C}\bigr](x) = \sum_{y\in\states} T(x,y)\mathbb{I}_\mathcal{C}(y) \neq 1\,,
\end{equation*}
which implies that there is $y\notin\mathcal{C}$ such that $T(x,y)>0$. Moreover, this $y$ is in $A^c$. To see this, suppose \emph{ex absurdo} that $y\in A$. Because $x\in\mathcal{C}$, we then find
\begin{equation*}
0=\bigl[\overline{T}\mathbb{I}_A\bigr](x) \geq \bigl[\overline{T}\mathbb{I}_{\{y\}}\bigr](x) \geq \bigl[T\mathbb{I}_{\{y\}}\bigr](x) = T(x,y) > 0\,,
\end{equation*}
a contradiction. So, indeed, $y\in A^c$, and $y\notin\mathcal{C}$.
Hence, there exists some $n\in\nats$ such that $\bigl[\overline{T}^n\mathbb{I}_A\bigr](y)>0$. But then
\begin{equation*}
\bigl[\overline{T}^{n+1}\mathbb{I}_{A}\bigr](x) \geq \bigl[T\overline{T}^{n}\mathbb{I}_{A}\bigr](x) \geq T(x,y)\bigl[\overline{T}^{n}\mathbb{I}_{A}\bigr](y) > 0\,,
\end{equation*}
where the second inequality used the fact that $T(x,z)\geq 0$ (because $T$ is a transition matrix) and $\bigl[\overline{T}^n\mathbb{I}_A\bigr](z) \geq 0$ for all $z\in\states$, due to~\lemmaref{lemma:composition_upper_dominates_composition_lower} and~\lemmaref{lemma:lower_trans_composite_basic_properties}. This implies $x\notin \mathcal{C}$, a contradiction.
\end{proof}

An immediate consequence of the previous property is that, when starting from a state in $\mathcal{C}$, the upper probability of ever hitting $A$, will be zero:
\begin{corollary}\label{cor:absorbing_dont_hit}
Let $\overline{p}_A^{(n)}$ and $\overline{p}_A^*$ be as in~\propositionref{prop:hit_prob_is_limit}. Then for all $n\in\nats$ and $x\in\mathcal{C}$, it holds that $\overline{p}_A^{(n)}(x)=0$. Moreover, therefore also $\overline{p}_A^*(x)=0$ for all $x\in\mathcal{C}$.
\end{corollary}
\begin{proof}
First note that $\overline{p}_A^{(0)}=\mathbb{I}_A$, so because $\mathcal{C}\subseteq A^c$ we immediately find that for all $x\in\mathcal{C}$, it holds that
\begin{equation*}
\overline{p}_A^{(0)}(x) = \mathbb{I}_A(x) = 0\,.
\end{equation*}
This provides the induction base for the remainder of the proof. Indeed, if we now assume, for some $n\in\natswith$, that $\overline{p}_A^{(n)}(x)=0$ for all $x\in\mathcal{C}$, then by~\lemmaref{lemma:absorbing_states_closed},
\begin{equation*}
\bigl[T\overline{p}_A^{(n)}\bigr](x) = \sum_{y\in\states} T(x,y)\overline{p}_A^{(n)}(y) = \sum_{y\in\mathcal{C}}T(x,y)\overline{p}_A^{(n)}(y) = 0\,,
\end{equation*}
for all $x\in\mathcal{C}$ and all $T\in\mathcal{T}$. Because $T\in\mathcal{T}$ was arbitrary, we find that also
\begin{equation*}
\bigl[\overline{T}\overline{p}_A^{(n)}\bigr](x) = 0\,,
\end{equation*}
for all $x\in\mathcal{C}$. Therefore, and because $\mathcal{C}\subseteq A^c$, we have for all $x\in\mathcal{C}$ that
\begin{equation*}
\overline{p}_A^{(n+1)}(x) = \mathbb{I}_A(x) + \mathbb{I}_{A^c}(x)\bigl[\overline{T}\overline{p}_A^{(n)}\bigr](x) = 0\,.
\end{equation*}
This concludes the proof of the first claim. That $\overline{p}_A^*(x) = \lim_{n\to+\infty} \overline{p}_A^{(n)}(x)=0$ for all $x\in\mathcal{C}$ is now immediate.
\end{proof}

We next prove some properties about the behaviour of the process when starting in states that \emph{can} reach $A$. Note that, for all $x\in A^c\setminus\mathcal{C}$ there is some $n\in\nats$ such that $\bigl[\overline{T}^n\mathbb{I}_A\bigr](x)>0$. We will denote the smallest such $n$ as follows.

\begin{definition}
For all $x\in A^c\setminus\mathcal{C}$, we let $n_x\in\nats$ be the smallest number such that $\bigl[\overline{T}^{n_x}\mathbb{I}_A\bigr](x)>0$.
\end{definition}

The following technical result will be useful throughout the rest of this appendix.
\begin{lemma}\label{lemma:non_absorbing_moves_to_non_absorbing}
	For all $x\in A^c\setminus \mathcal{C}$ with $n_x>1$, if for some $T\in\mathcal{T}$ it holds that $\bigl[T\overline{T}^{n_x-1}\mathbb{I}_A\bigr](x) > 0$,
	then there is some $y\in A^c\setminus\mathcal{C}$ such that $n_y=n_x-1$ and $T(x,y)>0$.
\end{lemma}
\begin{proof}
Fix any $x\in A^c\setminus \mathcal{C}$ with $n_x>1$, and suppose that for some $T\in\mathcal{T}$ it holds that
\begin{equation*}
\bigl[T\overline{T}^{n_x-1}\mathbb{I}_A\bigr](x) > 0\,.
\end{equation*}
Expanding this expression, we get
\begin{equation*}
0 < \bigl[T\overline{T}^{n_x-1}\mathbb{I}_A\bigr](x) = \sum_{y\in\states} T(x,y)\bigl[\overline{T}^{n_x-1}\mathbb{I}_A\bigr](y)\,,
\end{equation*}
which implies that there is some $y\in\states$ such that $T(x,y)\bigl[\overline{T}^{n_x-1}\mathbb{I}_A\bigr](y)>0$. Hence, we already know that $T(x,y)>0$. We will now show that $y\in A^c\setminus \mathcal{C}$. First, suppose \emph{ex absurdo} that $y\in A$. Because $T(x,y)>0$, if $y\in A$ we get
\begin{equation*}
0 < T(x,y) = \bigl[T\mathbb{I}_{\{y\}}\bigr](x) \leq \bigl[\overline{T}\mathbb{I}_{\{y\}}\bigr](x) \leq \bigl[\overline{T}\mathbb{I}_{\mathbb{I}_A}\bigr](x)\,,
\end{equation*}
which implies that $n_x=1$, a contradiction with the assumption made at the beginning of this proof. Hence $y\in A^c$. To see that $y\notin \mathcal{C}$, simply observe that $\bigl[\overline{T}^{n_x-1}\mathbb{I}_A\bigr](y)>0$, again by the selection of $y$. So, $y\in A^c\setminus\mathcal{C}$.

Moreover, because $\bigl[\overline{T}^{n_x-1}\mathbb{I}_A\bigr](y) >0$, it holds that $n_y\leq n_x-1$. Now, if $n_y=n_x-1$ then the proof is done; so, suppose \emph{ex absurdo} that $n_y<n_x-1$. Then
\begin{equation*}
\bigl[\overline{T}^{n_y+1}\mathbb{I}_A\bigr](x) \geq \bigl[T\overline{T}^{n_y}\mathbb{I}_A\bigr](x) \geq T(x,y)\bigl[\overline{T}^{n_y}\mathbb{I}_A\bigr](y) > 0\,,
\end{equation*}
where the second inequality used the fact that $T(x,z)\geq 0$ (because $T$ is a transition matrix) and $\bigl[\overline{T}^{n_y}\mathbb{I}_A\bigr](z) \geq 0$ for all $z\in\states$, due to~\lemmaref{lemma:composition_upper_dominates_composition_lower} and~\lemmaref{lemma:lower_trans_composite_basic_properties}. Since $n_x>n_y+1$, this implies that $n_x$ is not the smallest number such that $\bigl[\overline{T}^{n_x}\mathbb{I}_A\bigr](x)>0$, a contradiction. 
\end{proof}

\begin{corollary}\label{lemma:non_absorbing_steps}
 For all $x\in A^c\setminus\mathcal{C}$, if $n_x>1$, there is some $y\in A^c\setminus\mathcal{C}$ such that $n_y=n_x-1$.
\end{corollary}
\begin{proof}
Fix $x\in A^c\setminus\mathcal{C}$, and suppose that $n_x>1$. Then, by~\eqref{eq:lower_trans_reached}, there is some $T\in\mathcal{T}$ such that
\begin{equation}
\bigl[T\overline{T}^{n_x-1}\mathbb{I}_A\bigr](x) = \bigl[\overline{T}^{n_x}\mathbb{I}_A\bigr](x) > 0\,.
\end{equation}
Now apply~\lemmaref{lemma:non_absorbing_moves_to_non_absorbing}.
\end{proof}

\begin{corollary}\label{cor:non_absorbing_base}
Suppose that $A^c\setminus\mathcal{C}\neq\emptyset$. Then there is some $x\in A^c\setminus\mathcal{C}$ such that $n_x=1$.
\end{corollary}
\begin{proof}
Suppose that $A^c\setminus\mathcal{C}\neq\emptyset$, and choose any $x\in A^c\setminus\mathcal{C}$. If $n_x=1$ then we are done. Otherwise, by~\corollaryref{lemma:non_absorbing_steps}, there is some $y\in A^c\setminus\mathcal{C}$ such that $n_y=n_x-1$. Now repeat this selection until $n_y=1$.
\end{proof}

\begin{corollary}\label{cor:precise_non_absorbing_exists}
Suppose that $A^c\setminus\mathcal{C}\neq\emptyset$. Then there is some $T\in\mathcal{T}$ such that for all $x\in A^c\setminus\mathcal{C}$, 
$\bigl[T^{n_x}\mathbb{I}_A\bigr](x)>0$.
\end{corollary}
\begin{proof}
We construct the matrix $T\in\mathcal{T}$ in a row-by-row manner. First, choose any $S\in\mathcal{T}$ and, for all $x\in A\cup\mathcal{C}$, define the $x$-th row of $T$ as
\begin{equation*}
T(x,\cdot)\coloneqq S(x,\cdot)\,.
\end{equation*}
Next, for any $x\in A^c\setminus\mathcal{C}$, we can find $S_x\in\mathcal{T}$ such that
\begin{equation}\label{eq:row_select_move_to_interesting}
\bigl[S_x\overline{T}^{n_x-1}\mathbb{I}_A\bigr](x) = \bigl[\overline{T}^{n_x}\mathbb{I}_A\bigr](x) > 0\,,
\end{equation}
and we let the $x$-th row of $T$ be defined as
\begin{equation*}
T(x,\cdot) \coloneqq S_x(x,\cdot)\,.
\end{equation*}
Because $\mathcal{T}$ has separately specified rows, it then holds that $T\in\mathcal{T}$. It remains to verify the claim for this $T$.

The remainder of the proof is by induction. First consider $x\in A^c\setminus\mathcal{C}$ such that $n_x=1$; at least one such $x$ exists by~\corollaryref{cor:non_absorbing_base}. Then,
\begin{align*}
\bigl[T\mathbb{I}_A\bigr](x) &= \sum_{y\in\states} T(x,y)\mathbb{I}_A(y) \\
 &= \sum_{y\in\states} S_x(x,y)\mathbb{I}_A(y) \\
 &= \bigl[\overline{T}\mathbb{I}_A\bigr](x)>0\,.
\end{align*}
This provides the induction base. Now suppose the statement is true for all $x\in A^c\setminus\mathcal{C}$ for which $n_x\leq n$, for some $n\in\nats$. Then consider any $x\in A^c\setminus\mathcal{C}$ with $n_x=n+1$. 

By Equation~\eqref{eq:row_select_move_to_interesting}, it holds that
\begin{equation*}
0 < \bigl[S_x \overline{T}^{n_x-1}\mathbb{I}_A\bigr](x)\,.
\end{equation*}
Using~\lemmaref{lemma:non_absorbing_moves_to_non_absorbing}, this implies that there is some $y\in A^c\setminus\mathcal{C}$ such that $S_x(x,y)>0$ and $n_y=n_x-1=n$. Therefore, by the induction hypothesis, we know that
\begin{equation*}
\bigl[T^{n_x-1}\mathbb{I}_A\bigr](y) = \bigl[T^{n_y}\mathbb{I}_A\bigr](y) = \bigl[T^{n}\mathbb{I}_A\bigr](y) >0\,.
\end{equation*}
Hence, and because $T(x,\cdot) = S_x(x,\cdot)$ implies that $T(x,y) = S_x(x,y)>0$, it holds that
\begin{align*}
\bigl[T^{n_x}\mathbb{I}_A\bigr](x) &= \sum_{z\in\states} T(x,z)\bigl[T^{n_x-1}\mathbb{I}_A\bigr](z) \\
 &\geq T(x,y)\bigl[T^{n_x-1}\mathbb{I}_A\bigr](y) > 0\,.
\end{align*}
\end{proof}

We next use the construction of the matrix $T$ in the previous result, to construct sets $\mathcal{T}_\lambda\subseteq \mathcal{T}$ of transition matrices that, as we will see, are well-behaved enough to prove statements about upper hitting probabilities of imprecise Markov chains parameterised by these sets $\mathcal{T}_\lambda$.

\begin{definition}\label{definition:convex_subsets}
Fix any $T\in\mathcal{T}$ such that
$\bigl[T^{n_x}\mathbb{I}_A\bigr](x)>0$ for all $x\in A^c\setminus\mathcal{C}$; if $A^c\setminus\mathcal{C}\neq \emptyset$ then such a $T$ exists due to~\corollaryref{cor:precise_non_absorbing_exists}, otherwise any $T\in\transmatset$ will satisfy this vacuously.

Now for any $\lambda\in [0,1]$ we define the set $\mathcal{T}_\lambda$ as
\begin{equation*}
\mathcal{T}_\lambda \coloneqq \Bigl\{ \lambda T + (1-\lambda)V\,\Big\vert\,V\in\mathcal{T} \Bigr\}\,,
\end{equation*}
and let $\overline{T}_\lambda$ be the corresponding upper transition operator, defined for all $f\in\gamblesextabove$ and all $x\in\states$ as
\begin{equation*}
\bigl[\overline{T}_\lambda \,f\bigr](x) \coloneqq \sup_{S\in\mathcal{T}_\lambda} \bigl[S\,f\bigr](x)\,.
\end{equation*}
\end{definition}

For ease of exposition, we will take the selection of the matrix $T$ in ~\definitionref{definition:convex_subsets} to be fixed in the remainder of this appendix. These $\mathcal{T}_\lambda$ satisfy the following properties.

\begin{proposition}\label{prop:convex_nonabsorbing_subsets}
For any $\lambda\in[0,1]$ it holds that $\mathcal{T}_\lambda\subseteq \mathcal{T}$. Moreover, $\mathcal{T}_\lambda$ is non-empty, closed, convex, and has separately specified rows. 

Furthermore, for any $\lambda,\gamma \in[0,1]$, if $\lambda \leq \gamma$, then it holds that $\mathcal{T}_\gamma \subseteq \mathcal{T}_\lambda$.

Moreover, if $\lambda>0$ then for all $S\in \mathcal{T}_\lambda$ and all $x\in A^c\setminus\mathcal{C}$, 
it holds that $\bigl[S^{n_x}\mathbb{I}_A\bigr](x)>0$. 

Finally, $\lim_{\lambda\to 0}\overline{T}_\lambda f = \overline{T}f$ for all $f\in\gambles$.
\end{proposition}
\begin{proof}
Let $T\in\mathcal{T}$ denote the chosen matrix from~\definitionref{definition:convex_subsets}, and fix any $\lambda\in[0,1]$. First, observe that $\mathcal{T}_\lambda$ consists of convex combinations of elements of $\mathcal{T}$. Hence, the fact that $\mathcal{T}_\lambda\subseteq\mathcal{T}$ follows from the convexity of $\mathcal{T}$. 

Next, because $T\in\mathcal{T}$ and $\lambda\in[0,1]$, it holds that
\begin{equation*}
T = \lambda T + (1-\lambda) T \in \mathcal{T}_\lambda\,,
\end{equation*}
and so $\mathcal{T}_\lambda$ is non-empty.

To show that $\mathcal{T}_\lambda$ is closed, take any convergent sequence $\{S_n\}_{n\in\nats}$ in $\mathcal{T}_\lambda$ with $S_*\coloneqq \lim_{n\to+\infty} S_n$. We need to show that $S_*\in\mathcal{T}_\lambda$. If $\lambda=1$ then we trivially have $S_n=T$ for all $n\in\nats$, and hence $S_*=T\in\mathcal{T}_\lambda$. For the case that $\lambda<1$, we note that for all $n\in\nats$, because $S_n\in\mathcal{T}_\lambda$, there is some $V_n\in\mathcal{T}$ such that $S_n = \lambda T + (1-\lambda)V_n$. Therefore,
\begin{equation*}
S_* = \lim_{n\to+\infty} \Bigl(\lambda T + (1-\lambda) V_n \Bigr)= \lambda T + (1-\lambda) \lim_{n\to+\infty} V_n\,.
\end{equation*}
Hence the sequence $\{V_n\}_{n\in\nats}$ is convergent, with limit $\lim_{n\to+\infty} V_n=:V_*\in\mathcal{T}$ because $\mathcal{T}$ is closed. Thus,
\begin{equation*}
S_* = \lambda T + (1-\lambda) V_* \in \mathcal{T}_\lambda\,,
\end{equation*}
and so we conclude that $\mathcal{T}_\lambda$ is closed.

For the convexity of $\mathcal{T}_\lambda$, fix any $S_1, S_2$ in $\mathcal{T}_\lambda$, and any $\gamma\in[0,1]$, and let
\begin{equation*}
S \coloneqq \gamma S_1 + (1-\gamma) S_2\,.
\end{equation*}
We need to show that $S\in\mathcal{T}_\lambda$. Because $S_1,S_2\in\mathcal{T}_\lambda$, there exist $V_1,V_2\in\mathcal{T}$ such that
\begin{equation*}
S_1 = \lambda T + (1-\lambda) V_1\quad\text{and}\quad S_2 = \lambda T + (1-\lambda) V_2\,.
\end{equation*}
Therefore,
\begin{align*}
S &= \gamma\Bigl( \lambda T + (1-\lambda) V_1 \Bigr) + (1-\gamma) \Bigl( \lambda T + (1-\lambda) V_2 \Bigr) \\
 &= \gamma\lambda T + (1-\gamma)\lambda T + \gamma (1-\lambda) V_1 + (1-\gamma) (1-\lambda) V_2 \\
 &= \lambda T + (1-\lambda) \Bigl(\gamma V_1 + (1-\gamma) V_2\Bigr)\,.
\end{align*}
Because $\mathcal{T}$ is convex and $\gamma\in[0,1]$, there is some $V\in\mathcal{T}$ such that
\begin{equation*}
V = \gamma V_1 + (1-\gamma) V_2\,.
\end{equation*}
This implies that
\begin{equation*}
S = \lambda T + (1-\lambda) V \in \mathcal{T}_\lambda\,,
\end{equation*}
and so $\mathcal{T}_\lambda$ is convex.

To show that $\mathcal{T}_\lambda$ has separately specified rows, for all $x\in\states$, select any $S_x\in\mathcal{T}_\lambda$, and define the matrix $S$ as $S(x,\cdot)\coloneqq S_x(x,\cdot)$. We need to show that $S\in\mathcal{T}_\lambda$.

For all $x\in\states$, because $S_x\in\mathcal{T}_\lambda$, there is some $V_x\in\mathcal{T}$ such that $S_x = \lambda T + (1-\lambda) V_x$. Let $V$ be defined as $V(x,\cdot)\coloneqq V_x(x,\cdot)$ for all $x\in\states$; then $V\in\mathcal{T}$ because $\mathcal{T}$ has separately specified rows. Moreover, for any $x\in\states$, it holds that
\begin{align*}
S(x,\cdot) &= \lambda T(x,\cdot) + (1-\lambda) V_x(x,\cdot) \\
 &= \lambda T(x,\cdot) + (1-\lambda) V(x,\cdot)\,,
\end{align*}
and so
\begin{equation*}
S = \lambda T + (1-\lambda) V \in \mathcal{T}_\lambda\,,
\end{equation*}
from which we conclude that $\mathcal{T}_\lambda$ has separately specified rows.

Next, fix any $\lambda,\gamma\in[0,1]$ and suppose that $\lambda\leq \gamma$; we want to show that $\mathcal{T}_\gamma \subseteq \mathcal{T}_\lambda$. Clearly if $\lambda=1$ then also $\gamma=1$, so then the result is trivial; it remains to prove the case where $\lambda <1$. 

Fix any $S\in\mathcal{T}_\gamma$; it suffices to show that $S\in\mathcal{T}_\lambda$. Because $S\in\mathcal{T}_\gamma$, there is some $V\in\mathcal{T}$ such that $S=\gamma T + (1-\gamma)V$. Let $\delta \coloneqq \gamma-\lambda$; then  clearly $\lambda = \gamma -\delta$. Moreover,
\begin{align*}
S &= \gamma T + (1-\gamma)V \\ 
 &= \gamma T + V - \gamma V \\
 &= \gamma T - \delta T + \delta T + V - \gamma V + \delta V - \delta V \\
  &= \lambda T + V - \lambda V + \delta T - \delta V
\end{align*}
Next set
\begin{equation*}
W \coloneqq \frac{S - \lambda T}{1-\lambda}\,.
\end{equation*}
Then it holds that
\begin{equation*}
S = \lambda T + (1-\lambda) W\,,
\end{equation*}
so it follows that $S\in\mathcal{T}_\lambda$ if $W\in\mathcal{T}$. To show that this is the case, we use our previous expansion of $S$ to find
\begin{align*}
W &= \frac{S - \lambda T}{1-\lambda} \\
 &= \frac{\lambda T + V - \lambda V + \delta T - \delta V - \lambda T}{1-\lambda} \\
 &= \frac{(1-\lambda) V + \delta T - \delta V}{1-\lambda} \\
 &= V + \frac{\delta}{1-\lambda}T - \frac{\delta}{1-\lambda}V\,.
\end{align*}
Now note that $\delta=\gamma-\lambda\geq 0$ and, therefore, because $\gamma\in[0,1]$, it holds that $\Delta\coloneqq \nicefrac{\delta}{1-\lambda}\in[0,1]$. Hence,
\begin{equation*}
W = \Delta T + (1-\Delta)V \in \mathcal{T}\,,
\end{equation*}
because $\mathcal{T}$ is convex. Hence indeed $\mathcal{T}_\gamma\subseteq \mathcal{T}_\lambda$.

Next, suppose that $\lambda>0$ and fix any $S\in\mathcal{T}_\lambda$. Then there is some $V\in\mathcal{T}$ such that $S=\lambda T+(1-\lambda) V$. We will now prove by induction that $\bigl[S^{n_x}\mathbb{I}_A\bigr](x)>0$, for all $x\in A^c\setminus\mathcal{C}$. To avoid the vacuous truth of this statement, we can assume without loss of generality that $A^c\setminus\mathcal{C}\neq\emptyset$.

First, consider any $x\in A^c\setminus\mathcal{C}$ with $n_x=1$; at least one such $x$ exists due to~\corollaryref{cor:non_absorbing_base}. Then,
\begin{align*}
\bigl[S\mathbb{I}_A\bigr](x) &= \lambda\bigl[T\mathbb{I}_A\bigr](x) + (1-\lambda)\bigl[V\mathbb{I}_A\bigr](x) \\
 &\geq \lambda\bigl[T\mathbb{I}_A\bigr](x) > 0\,,
\end{align*}
where the first inequality used the fact that $V\mathbb{I}_A$ is non-negative (and that $1\geq \lambda$), and the second inequality used the selection of $T$ (and the fact that $\lambda>0$). This provides the induction base.

Now suppose the statement is true for all $x\in A^c\setminus \mathcal{C}$ for which $n_x \leq n$, for some $n\in\nats$, and consider any $x\in A^c\setminus\mathcal{C}$ with $n_x=n+1$. Then, by the selection of $T$, it holds that
\begin{equation*}
0 < \bigl[T^{n_x}\mathbb{I}_A\bigr](x) = \bigl[T T^{n_x-1}\mathbb{I}_A\bigr](x)\,.
\end{equation*}
Moreover, because $T\in\mathcal{T}$, it holds for all $y\in\states$ that
\begin{equation*}
\bigl[T^{n_x-1}\mathbb{I}_A\bigr](y) \leq \bigl[\overline{T}^{n_x-1}\mathbb{I}_A\bigr](y)\,,
\end{equation*}
and so we find that
\begin{equation*}
0 < \bigl[T T^{n_x-1}\mathbb{I}_A\bigr](x) \leq \bigl[T \overline{T}^{n_x-1}\mathbb{I}_A\bigr](x)\,.
\end{equation*}
Hence, because $T\in\mathcal{T}$ and $n_x>1$ by assumption, it follows from~\lemmaref{lemma:non_absorbing_moves_to_non_absorbing} that there is some $y\in A^c\setminus\mathcal{C}$ such that $n_y=n_x-1=n$ and $T(x,y)>0$.

For this $y$, it holds that $\bigl[S^n\mathbb{I}_A\bigr](y)>0$ by the induction hypothesis. Therefore,
\begin{align*}
\bigl[S^{n_x}\mathbb{I}_A\bigr](x) &\geq S(x,y)\bigl[S^n\mathbb{I}_A\bigr](y) \\
 &\geq \lambda T(x,y)\bigl[S^n\mathbb{I}_A\bigr](y) > 0\,,
\end{align*}
where the last inequality used the fact that $\lambda$, $T(x,y)$, and $\bigl[S^n\mathbb{I}_A\bigr](y)$ are all strictly positive. This concludes the proof that $\bigl[S^{n_x}\mathbb{I}_A\bigr](x)>0$ for all $x\in A^c\setminus\mathcal{C}$.

Let now, for any $\lambda\in[0,1]$, $\overline{T}_\lambda$ denote the upper transition operator corresponding to $\mathcal{T}_\lambda$. We will show that $\lim_{\lambda\to 0}\overline{T}_\lambda f=\overline{T}f$ for all $f\in\gambles$. So, consider any $f\in\gambles$ and any $\epsilon>0$. We then need to show that, for all $\lambda$ that are small enough, it holds that
\begin{equation}\label{prop:convex_nonabsorbing_subsets:eq:pointwise_limit}
\norm{\overline{T}_\lambda f - \overline{T}f} < \epsilon\,.
\end{equation}
Due to~\eqref{eq:lower_trans_reached}, there is some $V\in\mathcal{T}$ such that $Vf = \overline{T}f$. We will show that the inequality~\eqref{prop:convex_nonabsorbing_subsets:eq:pointwise_limit} holds for any $\lambda\in[0,1]$  for which
\begin{equation*}
\lambda \norm{Tf-Vf} < \epsilon\,.
\end{equation*}
Indeed, whenever this is the case, it holds for $\lambda T + (1-\lambda)V=:S\in\mathcal{T}_\lambda$ that
\begin{align*}
\norm{Sf - Vf} &= \norm{\lambda Tf + (1-\lambda)Vf - Vf} \\
 &= \lambda\norm{Tf - Vf} <  \epsilon\,.
\end{align*}
Moreover, because $S\in\mathcal{T}_\lambda$ and $\mathcal{T}_\lambda\subseteq\mathcal{T}$, it holds that
\begin{equation*}
Sf \leq \overline{T}_\lambda f \leq \overline{T}f = Vf\,,
\end{equation*}
from which it follows that also
\begin{equation*}
\norm{\overline{T}_\lambda f - \overline{T}f} < \epsilon\,.
\end{equation*}
\end{proof}

We next state some properties about the upper hitting probabilities for the imprecise Markov chains corresponding to the sets $\mathcal{T}_\lambda$.

\begin{lemma}\label{lemma:convex_nonabsorbing_subsets_systems}
For any $\lambda\in[0,1]$, define $\overline{p}_\lambda^{(0)}\coloneqq \mathbb{I}_A$ and, for all $n\in\natswith$, let
\begin{equation*}
\overline{p}_\lambda^{(n+1)} \coloneqq \mathbb{I}_A + \mathbb{I}_{A^c}\hprod \overline{T}_\lambda\overline{p}_\lambda^{(n)}\,.
\end{equation*}
Then, the limit
\begin{equation*}
\overline{p}_\lambda^* \coloneqq \lim_{n\to+\infty} \overline{p}_\lambda^{(n)}
\end{equation*}
exists and, moreover, $\overline{p}_\lambda^*$ is the minimal non-negative solution to
\begin{equation}\label{eq:lambda_system}
\overline{p}_\lambda^* = \mathbb{I}_A + \mathbb{I}_{A^c}\hprod \overline{T}_\lambda \overline{p}_\lambda^*
\end{equation}
in $\gambles$.
\end{lemma}
\begin{proof}
Fix any $\lambda\in[0,1]$. By~\propositionref{prop:convex_nonabsorbing_subsets}, $\mathcal{T}_\lambda$ is non-empty, closed, convex, and has separately specified rows. Therefore, it defines a game-theoretic imprecise Markov chain whose upper hitting probability satisfies
\begin{equation*}
\overline{\mathbb{E}}^\mathrm{V}_{\mathcal{T}_\lambda}\bigl[G_A\,\vert\,X_0\bigr] = \overline{p}_\lambda^* = \lim_{n\to+\infty}\overline{p}_\lambda^{(n)}\,,
\end{equation*}
due to~\propositionref{prop:hit_prob_is_limit}.

It remains to prove that $\overline{p}_\lambda^*$ is the minimal non-negative solution to~\eqref{eq:lambda_system} in $\gambles$. First, we note that it is clearly \emph{a} solution, using the argument surrounding~\eqref{eq:fixed_point_hit_prob}. 

To see that $\overline{p}_\lambda^*$ is a non-negative element of $\gambles$, we remark that, due to~\propositionref{prop:convex_nonabsorbing_subsets}, $\mathcal{T}_\lambda$ is non-empty, closed, convex, and has separately specified rows. Now apply~\lemmaref{lemma:hitting_prob_approx_sequence_nonnegative} to the sequence $\overline{p}_\lambda^{(n)}$ to find that $0 \leq \overline{p}_\lambda^{(n)}\leq 1$ for all $n\in\nats$. It follows immediately that also $\overline{p}_\lambda^* = \lim_{n\to+\infty} \overline{p}_\lambda^{(n)} \geq 0$ and, similarly, that $\overline{p}_\lambda^* \leq 1$.

To establish minimality, consider any non-negative $f\in\gambles$ such that
\begin{equation*}
f = \mathbb{I}_A + \mathbb{I}_{A^c}\hprod \overline{T}_\lambda f\,.
\end{equation*}
Note that at least one such $f$ exists because we can always choose $f=\overline{p}_\lambda^*$, which we already know to satisfy these properties.
Regardless of the exact choice of such an $f$, we will show that it always holds that $f\geq \overline{p}_\lambda^*$. First note that, for any $x\in A$, it holds that
\begin{equation*}
f(x) = \mathbb{I}_A(x) + \mathbb{I}_{A^c}(x)\bigl[ \overline{T}_\lambda f\bigr](x) = 1\,.
\end{equation*}
Therefore, and because $f$ is non-negative, it holds that $f\geq \mathbb{I}_A$.

Now we define the map $\overline{\mathbf{G}}:\gambles\to\gambles$, for all $g\in\gambles$, as
\begin{equation*}
\overline{\mathbf{G}}(g) \coloneqq \mathbb{I}_A + \mathbb{I}_{A^c}\hprod \overline{T}_\lambda g\,.
\end{equation*}
Moreover, for any $n\in\nats$, we let $\overline{\mathbf{G}}^n$ denote the $n$-fold composition of $\overline{\mathbf{G}}$ with itself. Then, clearly,
\begin{equation*}
f = \overline{\mathbf{G}}(f)\,.
\end{equation*}
Moreover, it is easy to see that, for all $n\in\nats$,
\begin{equation*}
\overline{p}_\lambda^{(n)} = \overline{\mathbf{G}}^n(\overline{p}_\lambda^{(0)}) = \overline{\mathbf{G}}^n(\mathbb{I}_A)\,.
\end{equation*}
Finally, it follows from the monotonicity of $\overline{T}_\lambda$ that $\overline{\mathbf{G}}$ is also a monotone operator and, therefore, that $\overline{\mathbf{G}}^n$ is monotone for any $n\in\nats$. Because, as we already established, $f\geq \mathbb{I}_A$, it follows that for every $n\in\nats$
\begin{equation*}
f = \overline{\mathbf{G}}^n(f) \geq \overline{\mathbf{G}}^n(\mathbb{I}_A) = \overline{p}_\lambda^{(n)}\,,
\end{equation*}
from which it follows that
\begin{equation*}
f \geq \lim_{n\to+\infty} \overline{p}_\lambda^{(n)} = \overline{p}_\lambda^*\,,
\end{equation*}
which concludes the proof.
\end{proof}

We will next show that the quantities $\overline{p}_\lambda^*$ are decreasing in $\lambda$.
\begin{lemma}\label{lemma:convex_nonabsorbing_subset_decreasing_limit}
For any $\lambda,\gamma \in[0,1]$, let $\overline{p}_\lambda^{(n)}$,  $\overline{p}_\lambda^{*}$,  $\overline{p}_\gamma^{(n)}$, and  $\overline{p}_\gamma^{*}$ be as in~\lemmaref{lemma:convex_nonabsorbing_subsets_systems}. Then if $\lambda \leq \gamma$, it holds that $\overline{p}_\lambda^* \geq \overline{p}_\gamma^*$.
\end{lemma}
\begin{proof}
Suppose that $\lambda \leq \gamma$. Due to~\propositionref{prop:convex_nonabsorbing_subsets}, this implies that $\mathcal{T}_\lambda \supseteq  \mathcal{T}_\gamma$, from which it follows that, for any $f\in\gambles$, it holds that $\overline{T}_\lambda f \geq \overline{T}_\gamma f$.

We will now show by induction that $\overline{p}_\lambda^{(n)} \geq \overline{p}_\gamma^{(n)}$ for all $n\in\natswith$. By definition it holds that $\overline{p}_\lambda^{(0)} = \mathbb{I}_A = \overline{p}_\gamma^{(0)}$, which provides the induction base that we are after. Now suppose that the statement is true for some $n\in\natswith$; we will show that it then also holds for $n+1$. To this end, observe that
\begin{align*}
\overline{p}_\lambda^{(n+1)} &= \mathbb{I}_A + \mathbb{I}_{A^c}\cdot\overline{T}_\lambda\overline{p}_\lambda^{(n)} \\
 &\geq \mathbb{I}_A + \mathbb{I}_{A^c}\cdot\overline{T}_\lambda\overline{p}_\gamma^{(n)} \\
 &\geq \mathbb{I}_A + \mathbb{I}_{A^c}\cdot\overline{T}_\gamma\overline{p}_\gamma^{(n)} = \overline{p}_\gamma^{(n+1)}\,,
\end{align*}
where the first inequality used the induction hypothesis together with the monotonicity of $\overline{T}_\lambda$ and the second inequality used the relation between $\overline{T}_\lambda$ and $\overline{T}_\gamma$ that we derived above, together with the fact that $\overline{p}_\gamma^{(n)}\in\gambles$ due to~\lemmaref{lemma:hitting_prob_approx_sequence_nonnegative}, which we can use because $\mathcal{T}_\gamma$ is non-empty, closed, convex and has separately specified rows due to~\propositionref{prop:convex_nonabsorbing_subsets}. Thus, we conclude that indeed $\overline{p}_\lambda^{(n)} \geq \overline{p}_\gamma^{(n)}$ for all $n\in\natswith$.

It is now immediate that
\begin{equation*}
\overline{p}_\lambda^* = \lim_{n\to+\infty} \overline{p}_\lambda^{(n)} \geq  \lim_{n\to+\infty} \overline{p}_\gamma^{(n)} = \overline{p}_\gamma^*\,.
\end{equation*}
\end{proof}

The next result tells us that, in order to find the hitting probabilities of an imprecise Markov chain corresponding to $\mathcal{T}$, we can look at the hitting probabilities of imprecise Markov chains corresponding to $\mathcal{T}_\lambda$ for $\lambda>0$, and take the limit $\lambda\to 0$.
\begin{lemma}\label{lemma:convex_limit_works}
Let $\overline{p}_A^*$ be as in~\propositionref{prop:hit_prob_is_limit}. Then it holds that
\begin{equation*}
\overline{p}_A^* = \overline{p}_0^* = \lim_{\lambda\to 0} \overline{p}_\lambda^*\,.
\end{equation*}
\end{lemma}
\begin{proof}
It follows immediately from the definition that $\mathcal{T}_0=\mathcal{T}$, and therefore from~\propositionref{prop:hit_prob_is_limit} and~\lemmaref{lemma:convex_nonabsorbing_subsets_systems} that $\overline{p}_A^*=\overline{p}_0^*$. It remains to show that the limit statement holds.

We remark that, for any $\lambda>0$, due to~\propositionref{prop:convex_nonabsorbing_subsets}, $\mathcal{T}_\lambda$ is non-empty, closed, convex, and has separately specified rows. Now apply~\lemmaref{lemma:hitting_prob_approx_sequence_nonnegative} to the sequence $\overline{p}_\lambda^{(n)}$ to find that $0\leq \overline{p}_\lambda^{(n)}\leq 1$ for all $n\in\nats$. It follows immediately that also $0 \leq \overline{p}_\lambda^* \leq 1$ for any $\lambda>0$.

Now, let
\begin{equation*}
\overline{p} \coloneqq \lim_{\lambda\to 0}\overline{p}_\lambda^*\,.
\end{equation*}
Note that $\overline{p}$ exists because the sequence $\overline{p}_\lambda^*$ is uniformly bounded in the interval $[0,1]$, and non-decreasing as $\lambda$ goes to zero; this last property follows from~\lemmaref{lemma:convex_nonabsorbing_subset_decreasing_limit}. Moreover, we immediately see that also $0\leq \overline{p}\leq 1$.

Also because of~\lemmaref{lemma:convex_nonabsorbing_subset_decreasing_limit}, it holds that $\overline{p}_\lambda^* \leq \overline{p}_0^*$ for all $\lambda\in[0,1]$, from which it follows that $\overline{p}\leq \overline{p}_0^*$.

Next, fix any $\lambda\in (0,1)$ and use~\eqref{eq:lambda_system} to derive
\begin{align*}
 &\quad \norm{\overline{p} - \mathbb{I}_A - \mathbb{I}_{A^c}\hprod \overline{T}_0\overline{p}} \\
  &= \norm{\overline{p} - \overline{p}_\lambda^* + \mathbb{I}_A + \mathbb{I}_{A^c}\hprod \overline{T}_\lambda\overline{p}_\lambda^* - \mathbb{I}_A - \mathbb{I}_{A^c}\hprod \overline{T}_0\overline{p}} \\
 &\leq \norm{\overline{p} - \overline{p}_\lambda^*} + \norm{\mathbb{I}_{A^c}\hprod \overline{T}_\lambda\overline{p}_\lambda^* - \mathbb{I}_{A^c}\hprod \overline{T}_0\overline{p}} \\
 &\leq \norm{\overline{p} - \overline{p}_\lambda^*} + \norm{\overline{T}_\lambda\overline{p}_\lambda^* -  \overline{T}_0\overline{p}} \\
 &\leq \norm{\overline{p} - \overline{p}_\lambda^*} + \norm{\overline{T}_\lambda\overline{p}_\lambda^* - \overline{T}_\lambda \overline{p}} + \norm{\overline{T}_\lambda \overline{p} -  \overline{T}_0\overline{p}} \\
 &\leq \norm{\overline{p} - \overline{p}_\lambda^*} + \norm{\overline{p}_\lambda^* - \overline{p}} + \norm{\overline{T}_\lambda \overline{p} -  \overline{T}_0\overline{p}}\,,
\end{align*}
where we used~\lemmaref{lemma:norm_properties_transop} in the last step; we can use this because, by~\propositionref{prop:convex_nonabsorbing_subsets}, $\mathcal{T}_\lambda$ is non-empty, closed, convex, and has separately specified rows, which means that the corresponding upper transition operator $\overline{T}_\lambda$ satisfies all the assumptions on which~\lemmaref{lemma:norm_properties_transop} (implicitly) relies. Now note that all summands vanish because $\lim_{\lambda\to 0}\overline{p}_\lambda^* = \overline{p}$ by definition, and $\lim_{\lambda\to 0}\overline{T}_\lambda \overline{p} = \overline{T}\overline{p} = \overline{T}_0\overline{p}$ due to~\propositionref{prop:convex_nonabsorbing_subsets}.

We conclude that
\begin{equation*}
\overline{p} = \mathbb{I}_A + \mathbb{I}_{A^c}\hprod \overline{T}_0\overline{p}\,.
\end{equation*}
This implies $\overline{p}\geq \overline{p}_0^*$ due to~\lemmaref{lemma:convex_nonabsorbing_subsets_systems} (which we can use since $0\leq \overline{p}\leq 1$ as we showed above, i.e. $\overline{p}$ is non-negative and real-valued) and, since we already established that $\overline{p}\leq\overline{p}_0^*$, we find that $\overline{p}=\overline{p}_0^*=\overline{p}_A^*$.
\end{proof}

\begin{lemma}\label{lemma:nonabsorbing_matrix_has_inverse}
Suppose that $A^c\setminus\mathcal{C}\neq\emptyset$. Fix any $\lambda\in (0,1)$, choose any $S\in\mathcal{T}_\lambda$, and define the $\lvert A^c\setminus\mathcal{C}\rvert \times \lvert A^c\setminus\mathcal{C}\rvert$ matrix $F$, for all $x,y\in A^c\setminus\mathcal{C}$, as $F(x,y)\coloneqq S(x,y)$. Let $I$ denote the identity matrix. Then $(I-F)^{-1}$ exists.
\end{lemma}
\begin{proof}
We will start by proving that $\lim_{n\to+\infty} F^n = 0$, where the right hand side is a zero matrix. Let $1$ denote the vector with constant value one. 

First we note that the linear map $F: \reals^{\lvert A^c\setminus\mathcal{C}\rvert} \to \reals^{\lvert A^c\setminus\mathcal{C}\rvert}$ is monotone. To see this, fix any $f,g\in  \reals^{\lvert A^c\setminus\mathcal{C}\rvert}$ such that $f\leq g$. Then, for any $x\in A^c\setminus\mathcal{C}$,
\begin{align*}
\bigl[Ff\bigr](x) &= \sum_{y\in A^c\setminus \mathcal{C}} F(x,y)f(y) \\
 &= \sum_{y\in A^c\setminus \mathcal{C}} S(x,y)f(y) \\
 &\leq \sum_{y\in A^c\setminus \mathcal{C}} S(x,y)g(y) \\
 &= \sum_{y\in A^c\setminus \mathcal{C}} F(x,y)g(y) = \bigl[Fg\bigr](x)\,,
\end{align*}
where the inequality used the fact that $S\in\mathcal{T}_\lambda$ is a transition matrix, which implies that $S(x,y)\geq 0$ for all $y\in A^c\setminus\mathcal{C}$. Because $x\in A^c\setminus\mathcal{C}$ was arbitrary, this implies that $Ff \leq Fg$, which is to say, $F$ is monotone. Note that this immediately implies that also $F^n$ is monotone for any $n\in\nats$. Moreover, clearly also the entries of $F^n$ are non-negative for any $n\in\nats$ (because the entries of $F$ are, as mentioned above).

Now note that $F1 \leq 1$. Therefore, by the monotonicity of $F$, for any $n\in\nats$ it holds that $F^n1 \leq 1$. Moreover, recall that for any $x\in A^c\setminus\mathcal{C}$, by~\propositionref{prop:convex_nonabsorbing_subsets}, it holds that $\bigl[S^{n_x}\mathbb{I}_A\bigr](x)>0$. 

Next consider any $x\in A^c\setminus\mathcal{C}$ such that $n_x=1$; at least one such $x$ exists by~\corollaryref{cor:non_absorbing_base} because we are assuming that $A^c\setminus\mathcal{C}\neq\emptyset$. Then $\bigl[S\mathbb{I}_A\bigr](x) >0$, which implies that there is some $y\in A$ such that $S(x,y)>0$. Therefore, since $S(x,\cdot)$ is a probability mass function, we have that
\begin{equation*}
\bigl[F1\bigr](x) = \sum_{z\in A^c\setminus\mathcal{C}} F(x,z) = \sum_{z\in A^c\setminus\mathcal{C}} S(x,z) < 1\,.
\end{equation*}
Hence, by the monotonicity of $F$ and the fact that $F^{n-1}1 \leq 1$, it holds for all $n\in\nats$ that
\begin{equation}\label{lemma:nonabsorbing_matrix_has_inverse:eq:eventually_contractive}
\bigl[F^n1\bigr](x) = \bigl[FF^{n-1}1\bigr](x) \leq \bigl[F1\bigr](x) < 1\,.
\end{equation}
We will now use this result to prove by induction that 
$\big[F^{n_x}1\bigr](x)<1$ for all $x\in A^c\setminus \mathcal{C}$. To establish the induction base, consider any $x\in A^c\setminus \mathcal{C}$ with $n_x=1$. Then $\big[F^{n_x}1\bigr](x)<1$ by~\eqref{lemma:nonabsorbing_matrix_has_inverse:eq:eventually_contractive}. Now suppose that there is some $n\in\nats$ such that this claim is true for all $x\in A^c\setminus \mathcal{C}$ with $n_x\leq n$, and then consider any $x\in A^c\setminus \mathcal{C}$ with $n_x=n+1$. Then it holds that
\begin{equation*}
\bigl[S^{n_x}\mathbb{I}_A\bigr](x) > 0\,.
\end{equation*}
Moreover, because $S\in\mathcal{T}_\lambda\subseteq\mathcal{T}$, it holds that 
\begin{equation*}
S^{n_x-1}\mathbb{I}_A \leq \overline{T}_\lambda^{n_x-1}\mathbb{I}_A \leq \overline{T}^{n_x-1}\mathbb{I}_A
\end{equation*}
and therefore
\begin{equation*}
\bigl[S\overline{T}^{n_x-1}\mathbb{I}_A\bigr](x) > 0\,.
\end{equation*}
By~\lemmaref{lemma:non_absorbing_moves_to_non_absorbing}, this implies that there is some $y\in A^c\setminus\mathcal{C}$ with $n_y=n_x-1=n$ and $S(x,y)>0$. Because $n_y=n$ it holds that $\bigl[F^{n_y}1\bigr](y)<1$ by the induction hypothesis. Furthermore, $F(x,y)=S(x,y)>0$, and therefore it holds that
\begin{align*}
\bigl[F^{n_x}1\bigr](x) &= \sum_{z\in A^c\setminus\mathcal{C}} F(x,z)\bigl[F^{n_y}1\bigr](z) \\
&= F(x,y)\bigl[F^{n_y}1\bigr](y) \\
 &\quad\quad\quad + \sum_{z\in A^c\setminus\bigl\{\mathcal{C}\cup\{y\}\bigr\}} F(x,z)\bigl[F^{n_y}1\bigr](z) \\
 &< F(x,y)+ \sum_{z\in A^c\setminus\bigl\{\mathcal{C}\cup\{y\}\bigr\}} F(x,z)\bigl[F^{n_y}1\bigr](z) \\
 &\leq F(x,y)+ \sum_{z\in A^c\setminus\bigl\{\mathcal{C}\cup\{y\}\bigr\}} F(x,z) \\
 &= \sum_{z\in A^c\setminus\mathcal{C}} S(x,z) \leq 1\,,
\end{align*}
where in the final inequality we used the fact that $F^{n_y}1\leq 1$ as proved above, and where we also used that $F(x,z)=S(x,z)\geq 0$ because $S$ is a transition matrix. This finishes the induction step. 

Now let $n\coloneqq \max_{x\in A^c\setminus\mathcal{C}}n_x$, and note that $\bigl[F^n1\bigr](x)<1$ for all $x\in A^c\setminus\mathcal{C}$; this follows because, for each $x\in A^c\setminus\mathcal{C}$, as we have just shown, $\bigl[F^{n_x}1\bigr](x)<1$, and therefore by the monotonicity of $F$, the fact that $n_x\leq n$, and $F^{n-n_x}1 \leq 1$ (as discussed above), we get
\begin{equation*}
\bigl[F^n1\bigr](x) = \bigl[F^{n_x}F^{n-n_x}1\bigr](x) \leq \bigl[F^{n_x}1\bigr](x) < 1\,,
\end{equation*}
for all $x\in A^c\setminus\mathcal{C}$. Let now
\begin{equation*}
p \coloneqq \max_{x\in A^c\setminus\mathcal{C}} \bigl[F^n1\bigr](x)\,.
\end{equation*}
Then $p<1$ because $A^c\setminus\mathcal{C}$ is finite. Fix any $m\in\nats$. To be explicit in what follows, let $\mathbf{p}$ be the vector with constant value $p$. Then, for all $x\in A^c\setminus\mathcal{C}$, it holds that
\begin{align*}
\bigl[F^{mn}1\bigr](x) &= \bigl[F^{(m-1)n}F^n1\bigr](x) \\
 &\leq \bigl[F^{(m-1)n}\mathbf{p}\bigr](x) \\
 &= p\bigl[F^{(m-1)n}1\bigr](x) \cdots \leq p^m\,,
\end{align*}
where the inequalities use the monotonicity of $F$ and its powers, as discussed at the beginning of this proof.

Because $p<1$ this implies that
\begin{equation*}
\lim_{m\to+\infty} \bigl[F^{mn}1\bigr] = 0\,,
\end{equation*}
and it easily follows from the monotonicity of $F$ that therefore also
\begin{equation*}
\lim_{m\to+\infty} \bigl[F^m1\bigr] = 0\,.
\end{equation*}
Therefore, and because the entries of each $F^m$ are non-negative, it holds that $\lim_{m\to+\infty} F^m=0$, where the right hand side is a zero matrix.

We are now ready to prove that $(I-F)^{-1}$ exists. To this end, choose any real vector $f$ and suppose that $(I-F)f=0$ or, in other words, that $f=Ff$. Iterating this, for any $n\in\nats$ we have that $f=F^nf$. Passing to the limit and using the above result, we find that
\begin{equation*}
f = \lim_{n\to+\infty} F^nf = 0\,.
\end{equation*}
This implies that the kernel of $I-F$ is trivial,
\begin{equation*}
\mathrm{ker}(I-F) = \{0\}\,,
\end{equation*}
whence $(I-F)^{-1}$ exists.
\end{proof}

The next result tells us that the upper hitting probability for an imprecise Markov chain corresponding to $\mathcal{T}_\lambda$, is obtained by a homogeneous Markov chain $P\in\imchom$.
\begin{lemma}\label{lemma:convex_set_reaches_upper_prob}
For any $\lambda\in(0,1)$ there is some $P\in \imchom$ such that $p_A^P = \overline{p}_\lambda^*$.
\end{lemma}
\begin{proof}
We first recall that, due to~\corollaryref{cor:absorbing_dont_hit}, $\overline{p}_A^*(x)=0$ for all $x\in\mathcal{C}$. By~\lemmaref{lemma:convex_nonabsorbing_subsets_systems}, $\overline{p}_\lambda^*$ is non-negative, and, furthermore, $\overline{p}_\lambda^*\leq \overline{p}_0^*$ by~\lemmaref{lemma:convex_nonabsorbing_subset_decreasing_limit} and $\overline{p}_0^*=\overline{p}_A^*$ by \lemmaref{lemma:convex_limit_works}. This implies that also $\overline{p}_\lambda^*(x)=0$ for all $x\in \mathcal{C}$. Moreover, by~\lemmaref{lemma:convex_nonabsorbing_subsets_systems}, $\overline{p}_\lambda^*$ satisfies
\begin{equation*}
\overline{p}_\lambda^* = \mathbb{I}_A + \mathbb{I}_{A^c}\hprod \overline{T}_\lambda\overline{p}_\lambda^*\,.
\end{equation*}
This implies that, for all $x\in A$,
\begin{equation*}
\overline{p}_\lambda^*(x) = \mathbb{I}_A(x) + \mathbb{I}_{A^c}(x)\bigl[ \overline{T}_\lambda\overline{p}_\lambda^*\bigr](x) = 1\,.
\end{equation*}
Next, because by~\propositionref{prop:convex_nonabsorbing_subsets} $\mathcal{T}_\lambda$ is a non-empty, closed, and convex set of transition matrices that has separately specified rows and upper transition operator $\overline{T}_\lambda$, using~\lemmaref{lemma:trans_op_continuous}, we can find some $S\in\mathcal{T}_\lambda$ such that
\begin{equation}\label{eq:upper_prob_reached_lin_sys}
\overline{p}_\lambda^* = \mathbb{I}_A + \mathbb{I}_{A^c}\hprod S\overline{p}_\lambda^*\,.
\end{equation}
Because $S\in \mathcal{T}_\lambda\subseteq\mathcal{T}$, this implies that there is some homogeneous Markov chain $P\in\imchom$ with transition matrix $S$ and hitting probability $p_A^P$. By~\lemmaref{lemma:homogen_hitting_prob_is_minimal_system_solution}, $p_A^P$ is also a solution of~\eqref{eq:upper_prob_reached_lin_sys}. We will show that $p_A^P=\overline{p}_\lambda^*$.

To this end, we will first show that also $p_A^P(x)=0$ for all $x\in \mathcal{C}$, and that $p_A^P(x)=1$ for all $x\in A$. First, by~\lemmaref{lemma:homogen_hitting_prob_is_minimal_system_solution}, $p_A^P$ is non-negative. Now note that, because  $P\in\imchom\subseteq \imcirr$, by~\propositionref{prop:vovk_imprecise_dominates_compatible_precise_vovk} and~\lemmaref{lemma:precise_vovk_hit_prob_equal_precise_measure_hit_prob}, it holds that
\begin{align*}
\uexpvovk\bigl[G_A\,\vert\,X_0\bigr] &\geq \sup_{Q\in\imcirr} \overline{\mathbb{E}}_Q^\mathrm{V}\bigl[G_A\,\vert\,X_0\bigr] \\ 
&= \sup_{Q\in\imcirr} \mathbb{E}_Q\bigl[G_A\,\vert\,X_0\bigr] \\
&= \uexpirr\bigl[G_A\,\vert\,X_0\bigr] \\
&\geq \uexphom\bigl[G_A\,\vert\,X_0\bigr] \geq \mathbb{E}_P[G_A\,\vert\,X_0]\,.
\end{align*}
We have already seen that, for $x\in \mathcal{C}$, it holds that $0=\overline{p}_A^*(x) = \uexpvovk\bigl[G_A\,\vert\,X_0=x\bigr]$ (using~\propositionref{prop:hit_prob_is_limit} for the second equality), so, because $p_A^P$ is non-negative, this implies that also $p_A^P(x)=\mathbb{E}_P[G_A\,\vert\,X_0=x]=0$.

For the values of $p_A^P$ in $A$, note that by~\lemmaref{lemma:homogen_hitting_prob_is_minimal_system_solution}, for all $x\in A$, $p_A^P$ satisfies
\begin{equation*}
p_A^P(x) = \mathbb{I}_A(x) + \mathbb{I}_{A^c}(x)\bigl[Sp_A^P\bigr](x) = 1\,.
\end{equation*}
We conclude that therefore $p_A^P$ agrees with $\overline{p}_\lambda^*$ on $\mathcal{C}$ and on $A$. Finally, it remains to prove that also $p_A^P(x)=\overline{p}_\lambda^*(x)$ for all $x\in A^c\setminus\mathcal{C}$. If $A^c\setminus\mathcal{C}=\emptyset$ then this holds vacuously; so, we can assume without loss of generality that $A^c\setminus\mathcal{C}\neq\emptyset$.

To this end, consider any non-negative $f\in\gambles$ such that $f(x)=0$ for all $x\in\mathcal{C}$, and such that
\begin{equation}\label{eq:upper_prob_reached_lin_sys_general}
f = \mathbb{I}_A + \mathbb{I}_{A^c}\cdot Sf\,.
\end{equation}
Before we move on, let us show that both $\overline{p}_\lambda^*$ and $p_A^P$ are valid choices of such $f$. First, for $\overline{p}_\lambda^*$, we note that it is non-negative and in $\gambles$ due to~\lemmaref{lemma:convex_nonabsorbing_subsets_systems}, satisfies $\overline{p}_\lambda^*(x)=0$ for all $x\in\mathcal{C}$ as argued at the beginning of this proof, and satisfies~\eqref{eq:upper_prob_reached_lin_sys_general} because it satisfies~\eqref{eq:upper_prob_reached_lin_sys}. For $p_A^P$, we remark that it is non-negative, in $\gambles$, and satisfies~\eqref{eq:upper_prob_reached_lin_sys_general} due to~\lemmaref{lemma:homogen_hitting_prob_is_minimal_system_solution}; the fact that $p_A^P(x)=0$ for all $x\in\mathcal{C}$ was argued above.

Because both $\overline{p}_\lambda^*$ and $p_A^P$ are valid choices for such an $f$, to finish the proof it suffices to find an expression for $f$ whose values on $A^c\setminus\mathcal{C}$ do not depend on the specific choice of $f$; this is what we will derive below.

First, for all $x\in A^c\setminus\mathcal{C}$, it holds that
\begin{align}
f(x) &= \mathbb{I}_A(x) + \mathbb{I}_{A^c}(x)\bigl[ S\,f\bigr](x) \nonumber \\
&= \bigl[S\,f\bigr](x) \label{eq:lemma:convex_set_reaches_upper_prob:compact_form} \\
&= \bigl[S \bigl(\mathbb{I}_A + \mathbb{I}_{A^c}\hprod S\,f\bigr)\bigr](x) \nonumber \\
&= \bigl[S\mathbb{I}_A\bigr](x) + \bigl[S\bigl(\mathbb{I}_{A^c}\hprod S\,f\bigr)\bigr](x)\,.\label{eq:lemma:convex_set_reaches_upper_prob:important_expansion}
\end{align}
Next, we will show that, for all $x\in\states$,
\begin{equation*}
\bigl[\mathbb{I}_{A^c}\hprod Sf\bigr](x) = \left\{\begin{array}{ll}
0 & \text{if $x\in A$, and} \\
0 & \text{if $x\in\mathcal{C}$, and} \\
f(x) & \text{otherwise}.
\end{array}\right.
\end{equation*}
This first case is trivial due to the occurrence of the indicator $\mathbb{I}_{A^c}$. The second case follows from $\bigl[S\mathbb{I}_\mathcal{C}\bigr](x)=1$ (this is~\lemmaref{lemma:absorbing_states_closed}) and the assumption that $f=0$ for all $x\in\mathcal{C}$; this implies
\begin{equation*}
\bigl[S f\bigr](x) = \sum_{y\in\states} S(x,y) f(y) = \sum_{y\in\mathcal{C}} S(x,y)f(y) = 0\,,
\end{equation*}
where the second equality used the fact that $S$ is a transition matrix. The final case considers that $x\in A^c\setminus\mathcal{C}$ and so follows from the identity~\eqref{eq:lemma:convex_set_reaches_upper_prob:compact_form}.

Define the $\lvert A^c\setminus\mathcal{C}\rvert\times \lvert A^c\setminus\mathcal{C}\rvert$ matrix $F$, for all $x,y\in A^c\setminus\mathcal{C}$, as $F(x,y)\coloneqq S(x,y)$. Let $I$ denote the identity matrix. Then by~\lemmaref{lemma:nonabsorbing_matrix_has_inverse}, and because we are working under the assumption that $A^c\setminus\mathcal{C}\neq\emptyset$, $(I-F)^{-1}$ exists. Moreover, for any $x\in A^c\setminus\mathcal{C}$,
\begin{align*}
\bigl[S\bigl(\mathbb{I}_{A^c}\hprod S\,f\bigr)\bigr](x) &= \sum_{y\in\states} S(x,y) \bigl[\mathbb{I}_{A^c}\hprod S\,f\bigr](y) \\
&= \sum_{y\in A^c\setminus \mathcal{C}} S(x,y)f(y) \\
&= \sum_{y\in A^c\setminus \mathcal{C}} F(x,y)f(y) \\
&= \left[F\left(f\big\vert_{A^c\setminus\mathcal{C}}\right)\right](x)\,,
\end{align*}
where $\cdot\big\vert_{A^c\setminus\mathcal{C}}$ denotes the restriction of $\cdot$ to $A^c\setminus\mathcal{C}$. Plugging this back into Equation~\eqref{eq:lemma:convex_set_reaches_upper_prob:important_expansion}, we can write
\begin{equation*}
f\big\vert_{A^c\setminus\mathcal{C}} = \bigl[S\mathbb{I}_A\bigr]\big\vert_{A^c\setminus\mathcal{C}} + F\left(f\big\vert_{A^c\setminus\mathcal{C}}\right)\,,
\end{equation*}
or in other words, using the fact that $f$ (and therefore $f\big\vert_{A^c\setminus\mathcal{C}}$) is real-valued by assumption,
\begin{equation*}
(I-F)(f\big\vert_{A^c\setminus\mathcal{C}}) = \bigl[S\mathbb{I}_A\bigr]\big\vert_{A^c\setminus\mathcal{C}}\,,
\end{equation*}
from which we get
\begin{equation*}
f\big\vert_{A^c\setminus\mathcal{C}} = (I-F)^{-1}\bigl[S\mathbb{I}_A\bigr]\big\vert_{A^c\setminus\mathcal{C}}\,,
\end{equation*}
using the existence of $(I-F)^{-1}$. 

Using the aforementioned fact that both $\overline{p}_\lambda^*$ and $p_A^P$ are valid choices for $f$, we find that
\begin{equation*}
\overline{p}_\lambda^*\big\vert_{A^c\setminus\mathcal{C}} = (I-F)^{-1}\bigl[S\mathbb{I}_A\bigr]\big\vert_{A^c\setminus\mathcal{C}} = p_A^P\big\vert_{A^c\setminus\mathcal{C}}\,,
\end{equation*}
which concludes the proof.
\end{proof}

We now have all the required results to prove the second part of~\theoremref{thm:imprecise_hit_probs_equal}.
\newline
\begin{proofof}{\theoremref{thm:imprecise_hit_probs_equal}}
The proof for the lower hitting probability was already given in the main text, so we here only give the proof for the upper hitting probability.

Fix any $x\in\states$. We note that $\imchom\subseteq \imcirr$ and, therefore,  by~\propositionref{prop:vovk_imprecise_dominates_compatible_precise_vovk} and~\lemmaref{lemma:precise_vovk_hit_prob_equal_precise_measure_hit_prob}, it holds that
\begin{align}
\uexpvovk\bigl[G_A\,\vert\,X_0=x\bigr] &\geq \sup_{Q\in\imcirr} \overline{\mathbb{E}}_Q^\mathrm{V}\bigl[G_A\,\vert\,X_0=x\bigr] \nonumber \\ 
&= \sup_{Q\in\imcirr} \mathbb{E}_Q\bigl[G_A\,\vert\,X_0=x\bigr] \nonumber \\
&= \uexpirr\bigl[G_A\,\vert\,X_0=x\bigr] \nonumber \\
&\geq \uexphom\bigl[G_A\,\vert\,X_0=x\bigr]\,. \label{eq:thm:imprecise_hit_probs_equal:vovk_bounds_homogen}
\end{align}
Hence, it suffices to prove that $\uexpvovk\bigl[G_A\,\vert\,X_0=x\bigr]\leq \uexphom\bigl[G_A\,\vert\,X_0=x\bigr]$ for all $x\in\states$, from which it then follows that also $\uexpvovk\bigl[G_A\,\vert\,X_0=x\bigr] = \uexpirr\bigl[G_A\,\vert\,X_0=x\bigr]$.

To this end, we first remark that $\uexpvovk\bigl[G_A\,\vert\,X_0=x\bigr] = \overline{p}_A^*(x)$ due to~\propositionref{prop:hit_prob_is_limit}. 
Now fix any $\epsilon>0$.
By~\lemmaref{lemma:convex_limit_works} it holds that $\lim_{\lambda\to 0}\overline{p}_\lambda^*=\overline{p}_0^*=\overline{p}_A^*$ which, because $\overline{p}_\lambda^*$ and $\overline{p}_0^*$ are real-valued due to~\lemmaref{lemma:convex_nonabsorbing_subsets_systems}, implies that there is some $\lambda \in(0,1)$ such that
\begin{equation*}
\overline{p}_\lambda^* \geq \overline{p}_A^* - \epsilon = \uexpvovk\bigl[G_A\,\vert\,X_0=x\bigr] -\epsilon\,.
\end{equation*}
Moreover, by~\lemmaref{lemma:convex_set_reaches_upper_prob}, there is some $P^*\in\imchom$ such that $p_A^{P^*}=\overline{p}_\lambda^*$. 
It then holds that
\begin{align*}
\uexphom\bigl[G_A\,\vert\,X_0=x\bigr] &= \sup_{P\in\imchom} p_A^P(x) \\
 &\geq p_A^{P^*}(x) = \overline{p}_\lambda^*(x) \\
  &\geq \uexpvovk\bigl[G_A\,\vert\,X_0=x\bigr] -\epsilon\,.
\end{align*}
Because $\epsilon>0$ was arbitrary this concludes the proof.
\end{proofof}

\begin{proofof}{\corollaryref{cor:imprecise_hitting_prob_is_minimal_system_solution}}
The proof for $\underline{p}_A$ was already given in the main text, so we here only give the proof for $\overline{p}_A$.

Due to~\theoremref{thm:imprecise_hit_probs_equal}, the upper hitting probability $\overline{p}$ is the same for any type of imprecise Markov chain; using~\propositionref{prop:hit_prob_is_limit}, let $\overline{p}=\overline{p}_A^*$ be this hitting probability. Then, by~\lemmaref{lemma:convex_limit_works}, $\overline{p}_0=\overline{p}_A^*$ and, by~\lemmaref{lemma:convex_nonabsorbing_subsets_systems}, $\overline{p}_0$ is the minimal non-negative solution to
\begin{equation*}
\overline{p}_0 = \mathbb{I}_A + \mathbb{I}_{A^c}\hprod \overline{T}\overline{p}_0\,,
\end{equation*}
which concludes the proof.
\end{proofof}

\end{document}